\documentclass[11pt]{article}

\usepackage{amsmath, amsfonts, amssymb,latexsym}

\usepackage{amsthm}
\usepackage{bbm}

\input epsf.sty

\newtheorem{Lemma}{Lemma}[section]
\newtheorem{Proposition}[Lemma]{Proposition}

\newtheorem{Definition}[Lemma]{Definition}

\newcommand{\BEQ}{\begin{equation}}     
\newcommand{\BEA}{\begin{eqnarray}}
\newcommand{\BD}{\begin{displaymath}}
\newcommand{\EEQ}{\end{equation}}       
\newcommand{\EEA}{\end{eqnarray}}
\newcommand{\ED}{\end{displaymath}}

\newcommand{\del}{\delta}

\newcommand{\eps}{\varepsilon}          

\newcommand{\supp}{{\mathrm{supp}}}

\newcommand{\Tr}{{\mathrm{Tr}}}

\newcommand{\Vol}{{\mathrm{Vol}}}

\newcommand{\R}{\mathbb{R}}
\newcommand{\C}{\mathbb{C}}

\newcommand{\N}{\mathbb{N}}

\newcommand{\Id}{{\mathrm{Id}}}

\def\proba{{\mathbb{P}}}

\def\esper{{\mathbb{E}}}

\def\sgn{{\mathrm{sgn}}}

\newcommand{\eop}{\hfill $\Box$}        

\newcommand{\II}{{\rm i}}               
\renewcommand{\Re}{{\rm Re\ }}          
\renewcommand{\Im}{{\rm Im\ }}          
\newcommand{\half}{{1\over 2}}          
\renewcommand{\vec}[1]{\boldsymbol{#1}} 

\catcode`\@=11
\def\numberbysection{\@addtoreset{equation}{section}
        \def\theequation{\thesection.\arabic{equation}}}
\numberbysection

\begin{document}

\vspace*{1.5cm}
\begin{center}
{\Large \bf Global fluctuations for 1D log-gas dynamics}

\end{center}

\vspace{2mm}
\begin{center}
{\bf  J\'er\'emie Unterberger$^a$}
\end{center}

\vskip 0.5 cm
\centerline {$^a$Institut Elie Cartan,\footnote{Laboratoire 
associ\'e au CNRS UMR 7502} Universit\'e de Lorraine,} 
\centerline{ B.P. 239, 
F -- 54506 Vand{\oe}uvre-l\`es-Nancy Cedex, France}
\centerline{jeremie.unterberger@univ-lorraine.fr}

\vspace{2mm}
\begin{quote}

\renewcommand{\baselinestretch}{1.0}
\footnotesize
{We study in this article the hydrodynamic limit in the macroscopic regime of the coupled system of stochastic
differential equations,
\BEQ d\lambda_t^i=\frac{1}{\sqrt{N}} dW_t^i - V'(\lambda_t^i) dt+ \frac{\beta}{2N}
\sum_{j\not=i} \frac{dt}{\lambda^i_t-\lambda^j_t}, \qquad i=1,\ldots,N, \EEQ
with $\beta>1$, sometimes called {\em generalized Dyson's Brownian motion}, describing the dissipative
dynamics of a log-gas of $N$ equal charges with equilibrium measure corresponding to a $\beta$-ensemble, with
sufficiently regular convex potential $V$. The limit $N\to\infty$ is known to satisfy a mean-field
Mac-Kean-Vlasov equation.  We prove that, for suitable initial conditions,  fluctuations around the limit are Gaussian and
satisfy an explicit PDE.

The proof is very much indebted to the harmonic potential case treated in
Israelsson
\cite{Isr}. Our key
argument consists in showing that the time-evolution generator may be written in the
form of  a transport operator on the upper half-plane, plus a bounded non-local operator
 interpreted in terms of a signed
 jump process.  }

\end{quote}

\vspace{4mm}
\noindent

 \medskip
 \noindent {\bf Keywords:} random matrices, Dyson's Brownian motion, log-gas, beta-ensembles, hydrodynamic limit, Stieltjes transform,
 entropy

\smallskip

\noindent
{\bf Mathematics Subject Classification (2010):}  60B20; 60F05; 60G20; 60J60; 60J75; 60K35

\newpage

\tableofcontents



\section{Introduction and statement of main results}


\subsection{Introduction}


Let $\beta\ge 1 $ be a fixed parameter, and  $N\ge 1$ an integer. We consider the following system of coupled stochastic differential equations
driven by $N$ independent standard Brownian motions $(W_t^1,\ldots,W_t^N)_{t\ge 0}$,

\BEQ d\lambda_t^i=\frac{1}{\sqrt{N}} dW_t^i - V'(\lambda_t^i) dt+ \frac{\beta}{2N}
\sum_{j\not=i} \frac{dt}{\lambda^i_t-\lambda^j_t}, \qquad i=1,\ldots,N  \label{eq:SDE} \EEQ

Letting
\BEQ {\cal W}(\{\lambda^i\}_i):=\sum_{i=1}^N V(\lambda^i)-\frac{\beta}{4N} \sum_{i\not=j}
\log(\lambda^i-\lambda^j),\EEQ
we can rewrite (\ref{eq:SDE}) as $d\lambda_t^i=\frac{1}{\sqrt{N}} dW_t^i-\nabla_i {\cal W}
(\lambda_t^1,\ldots,\lambda_t^N) dt$.
Thus the corresponding equilibrium measure,
\BEQ d\mu^N_{eq}(\{\lambda^i\}_i) = \frac{1}{Z^N_V} e^{-2N{\cal W}(\{\lambda^i\}_i)}=\frac{1}{Z^N_V} \left(\prod_{j\not=i}  |\lambda^j-\lambda^i|\right)^{\beta/2} \exp\left(-2N\sum_{i=1}^N V(\lambda^i)\right)
\ d\lambda^1\cdots d\lambda^N \label{eq:eq-measure} \EEQ
is that of a $\beta$-log gas with confining potential $V$. 

\medskip

\noindent Let us start with a historical overview of the subject as a motivation for our study. This system of equations was originally considered in a particular case by Dyson \cite{Dys} who wanted to
describe the Markov evolution of a Hermitian matrix $M_t$ with i.i.d. increments $dG_t$ taken from the 
Gaussian unitary ensemble (GUE). In Dyson's idea, this matrix-valued process was to be
a matrix analogue of Brownian motion. The latter time-evolution being invariant through conjugation
by unitary matrices, we may project it onto a time-evolution of the set of eigenvalues 
$\{\lambda_t^1,\ldots,\lambda_t^N\}$ of the matrix, and obtain (\ref{eq:SDE}) with
$\beta=2$ and $V\equiv 0$. Keeping $\beta=2$, it is easy to prove that (\ref{eq:SDE}) 
is  equivalent to a generalized matrix Markov evolution, $dM_t=dG_t-V'(M_t) dt$.  The Gibbs measure
$${\cal P}^N_V(M)=\frac{1}{Z_N} e^{-N\Tr V(M)} dM, \qquad dM=\prod_{i=1}^N dM_{ii}
\prod_{1\le i<j\le n} d\Re M_{ij} \ d\Im M_{ij} $$
 can then be proved to be an  equilibrium measure. Such measures, together with their
projection onto the eigenvalue set, $\mu^N_{eq}(\{\lambda^1,\ldots,\lambda^N\})$, are the main
object of random matrix theory, see e.g. \cite{Meh},\cite{And}, \cite{Pas}. The {\em equilibrium
eigenvalue distribution} can be studied by various means, in particular using orthogonal
polynomials with respect to the weight $e^{-NV(\lambda)}$. The scaling in $N$ (called
{\em macroscopic scaling} in random matrix theory) ensures
the convergence of the random point measure $X^N:=\frac{1}{N} \sum_{i=1}^N \del_{\lambda^i}$ to  a deterministic  measure $\mu_V$ with {\em compact} support and density $\rho$ when $N\to\infty$ (see
e.g. \cite{Joh}, Theorem 2.1). One finds e.g. the well-known semi-circle law, $\rho(x)=\frac{1}{\pi}\sqrt{2-x^2}$, when $V(x)=x^2/2$. Looking more closely at the limit of the point
measure, one finds for arbitrary {\em polynomial} $V$  (Johansson \cite{Joh}) Gaussian fluctuations of order $O(1/N)$, contrasting with the
$O(1/\sqrt{N})$ scaling of fluctuations for the means of $N$ independent random variables,
typical of the central limit theorem. Assuming that the support of the measure is
connected (this essential "one-cut" condition holding in particular for $V$ {\em convex}), Johansson proves that the {\em covariance} of the limiting law depends on $V$ only through the support of the measure -- it is thus
{\em universal} up to a scaling coefficient --, while the means is
equal to $\rho$, plus an apparently non-universal correction in $O(1/N)$. 

\medskip

Then Rogers and Shi \cite{RogShi}, disregarding the random matrix background,  studied directly for its sake the system (\ref{eq:SDE}) in the case where
$V$ is harmonic (i.e. quadratic) and $\beta=2$, which we call {\em Hermite case} henceforth (by reference to the
corresponding class of equilibrium orthogonal polynomials), proving in particular the following two facts:
\begin{itemize}
\item[(i)]  two arbitrary eigenvalues never collide, which implies the non-explosion of (\ref{eq:SDE});

\item[(ii)] the random point process $X^N_t:=\frac{1}{N} \sum_{i=1}^N \del_{\lambda^i_t}$
satisfies in the limit $N\to\infty$ a deterministic hydrodynamic equation of Mc-Kean Vlasov type, namely, the asymptotic density 
\BEQ \rho_t \equiv X_t:={\mathrm{w\!-\!\!\lim}}_{N\to\infty} X_t^N \EEQ
 satisfies the PDE
\BEQ \frac{\partial \rho_t(x)}{\partial t}=\frac{\partial}{\partial x}\left( \left( x-
p.v. \int \frac{dy}{x-y}\rho_t(y) \right)\rho_t(x)\right), \label{eq:McKV0} \EEQ
where $p.v. \int \frac{dy}{x-y}\rho_t(y)$ is a principal value integral. 
\end{itemize}

In the case studied by Rogers and Shi, explicit formulas for finite $N$ are known for the Markov
 generator in the form of a determinant, called {\em extended kernel}, see e.g. \cite{For},
 chapter XI, or  \cite{ForNag}), whose 
 asymptotics for $N\to\infty$ may in principle be used to study the macroscopic limit. This is accomplished by noting that a simple conjugation trick turns the generator of the process into an  $N$-particle Hamiltonian with a {\em one-body} potential only, whose eigenfunctions are deduced from those of one-particle Hamiltonians (actually, harmonic
 oscillators).  {\em For general
 $V$, however}, in marked contrast with respect to the equilibrium case,  no explicit formulas are known for finite $N$, even for $\beta=2$, since the  conjugation trick produces a supplementary two-body potential making the spectral problem unsolvable. To be more precise, Macedo and Macedo \cite{MacMac} classified all random matrix dynamics
which are unitary equivalent to imaginary-time evolution under a  Calogero-Sutherland
type Hamiltonian, providing explicit determinantal solutions in connection to classical
orthogonal polynomials when $\beta=2$; however, restricting to SDE's with  {\em additive} noise, the latter class contains only the Hermite case. Related models of diffusions
conditioned on non-intersecting, solvable in terms of classical orthogonal polynomials,
have been considered in Duits \cite{Dui}, who showed convergence of fluctuation field
to inhomogeneous Gaussian free field.   Then for $\beta\not=2$, 
 the finite $N$ equilibrium measure is not fully understood, even in the harmonic case,
 see \cite{ValVir}. 
 
\medskip
\noindent This makes the direct study of (\ref{eq:SDE}) for general $V$ and $\beta$ all the
more interesting. Whereas the PDE appearing in the hydrodynamical limit is known \cite{LiLiXie}, the law of fluctuations is not known in general, and it is the purpose of
this study, and of the forthcoming article \cite{Unt-law}, to fill this gap. 
 S. Li, X.-D. Li and Y.-X. Xie \cite{LiLiXie}, generalizing properties (i) and (ii)
 above, prove that the
random point process $X^N_t$ satisfies in the limit $N\to\infty$  a generalization of
the above Mc-Kean Vlasov equation, namely,
\BEQ \frac{\partial \rho_t(x)}{\partial t}=\frac{\partial}{\partial x} \left(
\left(V'(x)- \frac{\beta}{2} \ p.v. \int \frac{dy}{x-y}\,  \rho_t(y) \right) \rho_t(x)\right)
\label{eq:McKV} \EEQ
The equilibrium measure $\rho$, defined as the solution of the integral equation
$\frac{\beta}{2} p.v. \int \frac{dy}{x-y} \rho(y)=V'(x)$,  cancels the right-hand side 
of (\ref{eq:McKV}). Replacing  in (\ref{eq:SDE}) $\frac{\beta}{2N} \sum_{j\not=i} \frac{1}{\lambda_t^i-\lambda_t^j}$ with\\ $-\frac{1}{N} \sum_{j\not=i} \nabla U(\lambda_t^i-\lambda_t^j)$ where
$U$ is some convex two-body potential satisfying some very general properties of regularity
and growth at infinity, one may show that there appears in the same limit an
equation similar to (\ref{eq:McKV}),
\BEQ \frac{\partial \rho_t(x)}{\partial t}=\frac{\partial}{\partial x} \left(
\left(V'(x)+ \frac{\beta}{2} \ \int dy \, U'(x-y)  \rho_t(y) \right) \rho_t(x)\right) \EEQ
Solutions of this type of equations, common in plasma theory  and the study of
granular media \cite{GolRut,BenCag} and in particular,
the rate of convergence of these to equilibrium, have been studied
in details using Otto's infinite dimensional differential calculus \cite{Ott} in a series of papers,
see e.g.  \cite{Car}, \cite{OttVil}, \cite{Vil}. However, as already noted by S. Li, X.-D. Li
and Y.-X. Xie, the range of applicability of these papers, written under the assumption
that $U$ be {\em Lipschitz}, does not seem to extend to our case when $U(x)=-c\log|x|$. 
Since formally the law of fluctuations is obtained by linearizing the system of
equations (\ref{eq:SDE}) around its macroscopic limit $\rho$, it is clear that one must
find some way to deal with (\ref{eq:McKV}). 

\medskip
\noindent
Rogers' and Shi's approach to (\ref{eq:SDE}) has been successfully generalized to the
case of a harmonic potential with
arbitrary $\beta$ by Israelsson \cite{Isr} and Bender \cite{Ben}. The present study owes
very much to these two articles, so let us describe to some extent their contents. 
There are two main ideas. Let $Y^N_t:=N(X^N_t-X_t)$ be the rescaled fluctuation process
for finite $N$; we want to prove that $Y^N_t \overset{law}{\to} Y_t$ when $N\to\infty$ and
identify the law of the process $(Y_t)_{t\ge 0}$. First, It\^o's formula implies that
\BEQ d\langle Y_t^N,f_t\rangle=\half(1-\frac{\beta}{2}) \langle X_t^N,f''_t\rangle dt+
\frac{1}{\sqrt{N}} \sum_{i=1}^N f'_t(\lambda_t^i) dW_t^i \label{eq:Ito} \EEQ
if the test functions $(f_t)_{0\le t\le T}$, $f_t:\R\to\R$ solve the following linear PDE
  
\BEQ \frac{\partial f_t}{\partial t}(x)=V'(x)f'_t(x)-\frac{\beta}{4} \int \frac{f'_t(x)-f'_t(y)}{x-y} (X_t^N(dy)+X_t(dy)) \label{eq:intro-PDE-f} 
\EEQ
(see Proposition \ref{prop:Israelsson} below). Eq. (\ref{eq:intro-PDE-f}) is a dualized, linearized
version of (\ref{eq:McKV}).  Second,
eq. (\ref{eq:intro-PDE-f}) may be integrated in the harmonic case by means of a {\em Stieltjes transform} (see Definition \ref{def:Stieltjes}). Namely, the family of functions
$\{\frac{c}{\cdot-z}\}_{c\in\C,z\in\C\setminus\R}$ is preserved by (\ref{eq:intro-PDE-f}). The solution $\frac{c^N_t}{\cdot-z^N_t}$ at time $t$ with terminal condition $\frac{c^N_T}{\cdot-z^N_T}=\frac{c}{\cdot-z}$ is obtained by solving two coupled ordinary differential
equations for $c^N_t$ and $z^N_t$ depending on $X$ and the random point measure $X^N$ (see \cite{Isr}, Lemma 2). Substituting to $X^N$ its deterministic limit $X$ in the
r.-h.s. of (\ref{eq:intro-PDE-f}), one gets in a natural way a system of two coupled ordinary differential
equations for $(z_t)_{0\le t\le T}$, $(c_t)_{0\le t\le T}$ that  describes a solution of
the asymptotic limit of (\ref{eq:intro-PDE-f}) in the limit $N\to\infty$, namely,
\BEQ \frac{\partial f_t}{\partial t}(x)=V'(x)f'_t(x)-\frac{\beta}{2} \int \frac{f'_t(x)-f'_t(y)}{x-y} X_t(dy) \label{eq:intro-PDE-f-asymptotic} 
\EEQ

Bender interprets these equations as {\em characteristic equations} for a generalized
transport operator (see section \ref{sec:transport}) which is never stated explicitly.
Then (at least formally), It\^o's formula (see \cite{Isr}, p. 29) makes it possible to
find explicitly the Markov kernel in the limit $N\to\infty$. Namely, consider
a finite number of points $(z^k)_k$ in $\C\setminus\R$, and the solutions $(z^k_t)_{t\le T}$
of the corresponding characteristic equations with terminal condition $(z_T^k)_k$. Letting $f_t(x):=\sum_k
\frac{c^k_t}{x-z^k_t}$ be the solution of (\ref{eq:intro-PDE-f-asymptotic}), and $\phi_{f_t}(Y_t^N):=e^{\II\langle Y_t^N,f_t\rangle}$,
  \BEQ \esper[\phi_{f_T}(Y_T)\big| {\cal F}_t]=\esper[\phi_{f_t}(Y_t)] \exp\left( \half
\int_t^T  \left[ \II (1-\frac{\beta}{2}) \langle X_s,f''_s\rangle -
\langle X_s,(f'_s)^2\rangle \right] ds  \right) \label{eq:intro-Gaussian0} \EEQ
Since functions $f$ of the above form are dense in some appropriate Sobolev space, 
formula (\ref{eq:intro-Gaussian0}) allows to conclude that the limit process is Gaussian.
Then Bender solves explicitly the characteristic equations, which take on a particularly
simple form in the harmonic case, and deduces first the covariance of the Stieltjes
transform of the fluctuation process, Cov$(U_{t_1}(z_1),U_{t_2}(z_2))$,
$U_t(z):=\langle Y_t,\frac{1}{\cdot-z}\rangle$, and then (taking boundary values and
using the Plemelj formula, see section \ref{sec:Stieltjes}), the (distribution-valued)
covariance kernel Cov$(Y_{t_1}(x_1),Y_{t_2}(x_2))$. 

\bigskip
\noindent Our approach for the case of a general potential has exactly the
same starting point, but dealing with eq. (\ref{eq:intro-PDE-f}) turns out to be
more complicated than in the harmonic case.  The reason is that the family of 
functions $\left\{\frac{c}{\cdot-z}\right\}_{c\in\C,z\in\C\setminus\R}$ is no more preserved
by (\ref{eq:intro-PDE-f}): this is easily seen if $V$ is a polynomial or extends analytically
to a strip around the real axis, since
\BEQ -V'(x)\partial_x \left(\frac{1}{x-z}\right)=(V'(z)\partial_z+V''(z))
\left( \frac{1}{x-z}\right) + \left( \frac{V^{(3)}(z)}{2!} + \frac{V^{(4)}(z)}{3!} (x-z)+
\ldots \right)\EEQ
features extra unwanted polynomial terms.
In practice we need only assume that $V$ is sufficiently regular, and  (letting $z=:a+\II b$)
write instead, for $a$ in a neighbourhood of the support of the random point measure
\BEQ V'(x)=V'(a)+V''(a)(x-a)+V'''(a) \frac{(x-a)^2}{2} +(x-a)^3 W_a(x-a), \label{eq:intro-W} \EEQ
 and find for the first three terms,
\BEQ -V'(a)\partial_x \left(\frac{1}{x-z}\right)=V'(a)\partial_a \left(\frac{1}{x-z}\right),
\ \  -V''(a)(x-a) \partial_x \left(\frac{1}{x-z}\right)=V''(a)\left(1+b\partial_b\right)
\frac{1}{x-z}, \label{eq:1.14} \EEQ
\BEQ -V'''(a) \frac{(x-a)^2}{2} \partial_x\left(\frac{1}{x-z}\right)=\half V'''(a)+
\half V'''(a) (2\II b+b^2\partial_a) \frac{1}{x-z}, \label{eq:1.15} \EEQ
defining a generalized transport operator 
\BEQ -V'(a)\partial_a-V''(a)(1+b\partial_b)-
\half V'''(a) (2\II b+b^2\partial_a). \label{eq:generalized-transport-operator} 
\EEQ
The new piece is the last (Taylor's remainder) term in (\ref{eq:intro-W}). We must give up at this point the
 idea that 
the time-evolution is a simple characteristic evolution, and prove that the Taylor
remainder produces instead a {\em non-local kernel}.  Let us highlight the main points
while avoiding technicalities. The main tool here is the use
of {\em Stieltjes decompositions} of order $\kappa$ (see Definition
 \ref{def:Stieltjes-decomposition}): 
for any $b_{max}>0$ and $\kappa=0,1,2,\ldots$, any sufficient regular, integrable
 function $f:\R\to\R$ may be written as an integral over the strip
$\Pi_{b_{max}}:=\{a\pm\II b\ |\  0<|b|<b_{max} \}$
\BEQ f(x)=\int_{-\infty}^{+\infty} da \int_{-b_{max}}^{b_{max}} db\, (-\II b) \frac{|b|^{\kappa}}{(1+\kappa)!} \,  \,  {\mathfrak{f}}_z(x)
\, h(a,b) \label{eq:intro-Stieltjes-decomposition} \EEQ 
\BEQ {\mathrm{where}}\qquad  {\mathfrak{f}}_z(x):=\frac{1}{x-z}. \EEQ
The mapping $f\mapsto h$ is clearly not one-to-one. Explicit Stieltjes decompositions
 are produced in \cite{Isr}, Lemma 9;  part of the job consists precisely in {\em choosing
} Stieltjes decompositions with good properties.  Let $\kappa'\ge\kappa\ge 0$. Inserting the time-evolution operator (\ref{eq:intro-PDE-f}) into  (\ref{eq:intro-Stieltjes-decomposition}), we prove
that:

-- the $(\frac{1}{x})$-potential and the second-order Taylor expansion of the operator
$V'(x)\partial_x$, see (\ref{eq:1.14},\ref{eq:1.15}), act together as a
{\em transport operator} ${\cal H}^{\kappa}:L^1(\Pi_{b_{max}})\to L^1(\Pi_{b_{max}})$;

-- the {\em Taylor remainder term} (see \S 3.7), to which one must add an inessential
off-support contribution (see \S 3.8) and  boundary terms (see \S 3.9), may be realized as a non-local operator
$|b|^{\kappa'-\kappa}{\cal H}^{\kappa';\kappa}_{nonlocal}(t)$ acting on the coefficient function $h$,
\BEQ |b|^{\kappa'-\kappa}({\cal H}^{\kappa';\kappa}_{nonlocal}(t))(h)(a,b):=|b|^{\kappa'-\kappa}\int_{-\infty}^{+\infty} da_T \int_{-b_{max}}^{b_{max}} db_T
\, g_{nonlocal}^{\kappa';\kappa}(a,b;a_T,b_T)h(a_T,b_T) \EEQ
such that 
\BEQ {\cal H}^{\kappa+1;\kappa}_{nonlocal}(t) :L^1(\Pi_{b_{max}})\to L^1(\Pi_{b_{max}})  \label{eq:regularizing} \EEQ 
are {\em bounded}. From (\ref{eq:regularizing}) we get: $b{\cal H}^{\kappa+1;\kappa}_{nonlocal}(h)(a,b)= \tilde{h}(a,b)$ with $\tilde{h}\in L^1(\Pi_{b_{max}},|b|^{-1}\, da\,
 db)$.  In other words, the non-local part of the time evolution  is  (in some weak sense) {\em regularizing} near the real axis, and acts therefore as a bounded perturbation of 
 ${\cal H}^{\kappa}_{transport}$.

\subsection{Notations and basic facts}


In this paragraph, we simply assume that $V$ is convex.
The filtration $({\cal F}_t)_{t\ge 0}$ is the filtration of the Brownian $(W^i_t)_{t\ge 0,i=1,\ldots,N}$.

\begin{Definition} \label{def:empirical-fluctuation}

\begin{enumerate}
\item
Let 
\BEQ X_t^N:=\frac{1}{N} \sum_{i=1}^N \del_{\lambda_t^i} \EEQ
be the empirical measure process.
\item
Call
\BEQ  Y_t^N:=N(X_t^N-X_t) \EEQ
 the finite $N$ {\em fluctuation process}.
\end{enumerate}

\end{Definition}

\noindent

The case developed in \cite{Isr} and \cite{Ben} is the {\em harmonic case}, 
$V(x)=\half x^2$ (up to normalization), to which we shall often refer. 
Apart from this very particular case, classical examples include the Landau-Ginzburg
potential  $V(x)=\half x^2+ \frac{\lambda}{4}x^4$, $\lambda\ge 0$, for which the support of the equilibrium measure is connected, and the density
is the product of the  (rescaled) semi-circle law  by some explicit polynomial
of degree 2 (see e.g. \cite{Joh}, p. 164).

It is proved in (Li-Li-Xie \cite{LiLiXie}, Theorem 1.3) that, provided
$X_0^N\overset{N\to\infty}{\to}\rho_0$, a deterministic density,  the empirical measure process $(X_t^N)$
converges in law to a deterministic measure process $(X_t)_{t\ge 0}$ with density $\rho_t$
solution of the non-linear Fokker-Planck equation,
\BEQ \frac{\partial \rho_t(x)}{\partial t}=\frac{\partial}{\partial x} \left(
\left(V'(x)- \frac{\beta}{2} \ p.v. \int \frac{dy}{x-y}\,  \rho_t(y) \right) \rho_t(x)\right)
\label{eq:rho-t} \EEQ
with initial condition $\rho_0$.
Equivalently, for any test function $f$,
\BEA &&  \frac{d}{dt} \langle X_t,f\rangle=\frac{d}{dt} \int \rho_t(x) f(x) \, dx=
-\int V'(x) f'(x) X_t(dx) + \nonumber\\
&&\qquad \qquad +  \frac{\beta}{4} \int \int \frac{f'(x)-f'(y)}{x-y}\,  X_t(dx)\,  X_t(dy). \label{eq:Xt-t} \EEA

The equilibrium measure $\mu^N_{eq}$, see (\ref{eq:eq-measure}) converges weakly when
$N\to\infty$ to  the stationary, deterministic  solution $X_t(dx)=\rho_{eq}(x)dx$ of (\ref{eq:rho-t}), where
$\rho_{eq}$ is the solution of the following integral equation, called {\em cut equation},
\cite{Joh}
\BEQ p.v.\ \int \frac{\rho_{eq}(x)\, dx}{x-y}=-\frac{2}{\beta} V'(y)  \label{eq:cut-eq1} \EEQ

\medskip

Formula (\ref{eq:Xt-t}) is formally obtained as in \cite{RogShi} by taking the limit $N\to\infty$ in
the finite $N$ It\^o formula (eq. (3) in \cite{Isr}),
\BEA && d\langle X_t^N,f\rangle = \left(\frac{\beta}{4} \int \int \frac{f'(x)-f'(y)}{x-y}
 X^N_t(dx)X^N_t(dy)  - \int V'(x) f'(x) X^N_t(dx) \right) dt \nonumber\\
 && \qquad\qquad \qquad  + \half(1-\frac{\beta}{2}) \frac{1}{N}
  \langle X_t^N,f''\rangle \, dt +
 \frac{1}{N\sqrt{N}} \sum_{i=1}^N f'(\lambda_t^i)dW_t^i. \label{eq:Ito0} \EEA

Roughly speaking, both terms in the second line of (\ref{eq:Ito0}) are $O(\frac{1}{N})$
(the argument for the martingale term  relies on an $L^2-$bound based on the independence of the $(W^i)_{1\le i\le N}$).

\begin{Definition}[Stieltjes transform] \label{def:Stieltjes}
Fix $z\in\C\setminus\R$.
\begin{itemize}
\item[(i)] Let ${\mathfrak{f}}_z(x):=\frac{1}{x-z}$ ($x\in\R$).
\item[(ii)]
Let, for $z\in\C\setminus\R$, 
\BEQ M_t^N:=\langle X_t^N,{\mathfrak{f}}_z\rangle=\sum_{i=1}^N \frac{1}{\lambda^i_t-z}\EEQ
and 
\BEQ M_t:=\langle X_t,{\mathfrak{f}}_z\rangle=\int \frac{\rho_t(x)}{x-z}\, dx\EEQ
be the Stieltjes transform of $X_t^N$, resp. $X_t$.
\end{itemize}
\end{Definition}

Starting from the cut equation (\ref{eq:cut-eq1}) and applying
Plemelj's formula (see section \ref{sec:Stieltjes}), one finds {\em at equilibrium}
\BEQ M(x+\II 0)=-\frac{2}{\beta} V'(x)+\II \pi\rho_{eq}(x)   \label{eq:cut-eq2} \EEQ
\BEQ M(x+\II 0)-M(x-\II 0)=2\II \pi\rho_{eq}(x),\qquad M(x+\II 0)+M(x-\II 0)=-\frac{4}{\beta}V'(x).\label{eq:cut-eq3} \EEQ

\medskip

A PDE for the Stieltjes transform of $X_t$ is determined easily from (\ref{eq:Xt-t}),

\BEQ \frac{\partial M_t}{\partial t}=\frac{\partial}{\partial z} \left(
\frac{\beta}{4} (M_t(z))^2+ V'(z)M_t(z)+T_t(z) \right),  \label{eq:T}\EEQ
 where   
\BEQ T_t(z):=\int \frac{V'(x)-V'(z)}{x-z} \, X_t(dx).\EEQ
{\em In the harmonic case}, $T$ is simply a constant, hence $M$ is the solution of
a complex Burgers equation on $\C\setminus\R$, see e.g. \cite{Isr}, eq. (6). However, this is no more the case in our general setting, and
the non-local term $T$ in the right-hand side prevents any explicit solution of
the equation. Yet
the Stieltjes transform will turn out to be a very convenient technical tool in the
computations.

\bigskip

\noindent The first idea, coming from \cite{Isr}, is to {\em transfer the drift in the time-evolution of $Y_t^N$ to the test function $f$}. This is done through a {\em straightforward generalization} of (Israelsson \cite{Isr}, Lemma 1):

\begin{Proposition}(see Israelsson \cite{Isr}) \label{prop:Israelsson}
Assume the following event holds for some constant $R>0$,
\BEQ \Omega_R: \ \sup_{0\le t\le T} \max_{i=1,\ldots,N} |\lambda^i|\le R; \qquad \forall t\le T,\ \supp(X_t)\subset[-R,R], \EEQ
i.e. that the support of the random point measure $X_t^N$ and of the measure
$X_t$ is $\subset[-R,R]$ for $0\le t\le T$. 
Let $(f_t)_{0\le t\le T}$, $f_t:\R\to\R$ be such that  
\BEQ \frac{\partial f_t}{\partial t}(x)=V'(x)f'_t(x)-\frac{\beta}{4} \int \frac{f'_t(x)-f'_t(y)}{x-y} (X_t^N(dy)+X_t(dy)) \label{eq:PDE-f}  \EEQ
for all $|x|\le R$.
Then
\BEQ d\langle Y_t^N,f_t\rangle=\half(1-\frac{\beta}{2}) \langle X_t^N,f''_t\rangle dt+
\frac{1}{\sqrt{N}} \sum_{i=1}^N f'_t(\lambda_t^i) dW_t^i.\EEQ
\end{Proposition}

As emphasized in the above Proposition, (\ref{eq:PDE-f}) need only hold on
$[-R,R]$, because $\langle X_t^N,f\rangle$ and $\langle X_t,f\rangle$ do {\em not} depend  on the values of $f$
on $\C\setminus[-R,R]$.

\medskip\noindent
The above Proposition is a direct consequence of It\^o's formula applied to the
fluctuation process (just subtract (\ref{eq:Ito0}) from (\ref{eq:Xt-t})),
\BEA &&  d\langle Y_t^N,f\rangle=\Big(\frac{\beta}{4} \int \int \frac{f'(x)-f'(y)}{x-y}
\left[X_t^N(dx)+X_t(dx)\right] Y_t^N(dy)-\int V'(x)f'(x)Y_t^N(dx) 
\Big) \, dt \nonumber\\
&& \qquad \qquad +\half(1-\frac{\beta}{2}) \langle X_t^N,f''\rangle \, dt + \frac{1}{\sqrt{N}}
\sum_{i=1}^N f'(\lambda_t^i)dW_t^i.  \label{eq:dY_tnf} \EEA

\medskip
Then Israelsson solves eq. (\ref{eq:PDE-f}) in the harmonic case by using as test
functions the $\frac{c}{\cdot-z}$, $c\in\C,z\in\C\setminus\R$, on which the
generator of time-evolution acts in a particularly simple way. We do not reproduce their
results here however, since they do not separate the analysis of the term due to the
harmonic potential from that due to the two-body logarithmic potential. We shall analyze
(\ref{eq:PDE-f}) in section 3 after we have introduced Stieltjes decompositions. 

\bigskip

{\em Normalization:} the reader willing to compare our results with those of
Israelsson \cite{Isr} or Bender \cite{Ben} should take into account the different
choices of normalization.  Compared to \cite{Isr}, we fix $\sigma=1$ and let $\alpha=\frac{\beta}{2}$,
$\gamma=\half(1-\frac{\beta}{2})$.
After  rescaling the $\lambda^i$'s by a factor $\beta^{-1/2}$, we obtain for $V$ quadratic   \cite{Ben}
with $\sigma=\half$.


\subsection{Main result and outline of the article}


 {\bf Assumptions on $V$.}
\medskip

\noindent {\em We assume $V$ to be a convex function in ${\cal C}^{11}$.}

\bigskip

\noindent Main examples are convex polynomials, or suitable, smooth perturbations
thereof.

\bigskip\noindent
Under our assumptions (see e.g. \cite{Joh}, Theorem 2.1 and Proposition 3.1) the equilibrium measure
$\rho_{eq}$ is well-defined and compactly supported, its support $[a,b]$ is connected, and $\rho_{eq}$ is a solution of the cut-equation (\ref{eq:cut-eq1}).

\medskip

\noindent  

\medskip\noindent {\bf Assumptions on the initial  measure.}

\medskip\noindent Let $\mu^N_0=\mu_0(\{\lambda^i_0\}_i)$ be the initial measure of the stochastic process
$\{\lambda^i_t\}_{t\ge 0, i=1,\ldots,N}$, and $X_0^N:=\frac{1}{N}\sum_{i=1}^N
\del_{\lambda^i_0}$ be the initial empirical measure. Since $N$ varies,
we find it useful here to add an extra upper index $(\lambda^{N,i}_0)_{
i=1,\ldots,N}$ to denote the initial condition of the process for a given value of $N$.   {\em We assume that:

\medskip
\begin{itemize}
 \item[(i)]  (large deviation estimate for the initial support) there exist some constants $C_0,c_0,R_0>0$  such that, for every $N\ge 1$,
\BEQ \proba[ \max_{i=1,\ldots,N} |\lambda_0^{N,i}| > R_0] \le 
C_0 e^{-c_0 N }. \label{eq:support-bound0} \EEQ
 \item[(ii)]  $X_0^N\overset{law}{\to}
\rho_0(x)\, dx$ when $N\to\infty$, where $\rho_0(x)$ is a deterministic measure;
\item[(iii)] (rate of convergence)  
\BEQ \left(\esper[|M_0^N(z)-M_0(z)|^2]\right)^{1/2}=O(\frac{1}{N b}) \label{hyp:initial-measure} \EEQ
 for $z=a+\II b\in\C\setminus\R$, where $M_0(z):=\int dx\, \frac{\rho_0(x)}{x-z}$ is the
Stieltjes transform of $\rho_0$.

\end{itemize}}

\medskip\noindent

 Lemma \ref{lem:support} proves  that
the initial large deviation estimate (i) implies a uniform-in-time large deviation estimate for the
support of the random point measure, which is essential for our main result.

\medskip
\begin{Definition}[Sobolev spaces] \label{def:Sobolev}
Let $H_n:=\{f\in L^2(\R)\ |\ ||f||_{H_n}<\infty\}$ $(n\ge 0)$, where 
$||f||_{H_n}:=\left(\int d\xi\, (1+|\xi|^2)^{n} |{\cal F}f(\xi)|^2\right)^{1/2}$, and 
$H_{-n}:=(H_n)'$ its dual. 
\end{Definition}

\medskip
\noindent The measure-valued process $Y^N$ may be shown to converge in $C([0,T],H_{-14})$:

\bigskip

\noindent {\bf Main Theorem (Gaussianity of limit fluctuation process).}

\noindent {\em  Let $Y_t^N$ be the finite $N$ fluctuation process (see Definition \ref{def:empirical-fluctuation}). Then:
\begin{enumerate}
\item $Y^N \overset{law}{\to} Y$ when $N\to\infty$, where $Y$ is a  Gaussian process.
More precisely, $Y^N$ converges to $Y$ weakly in $C([0,T],H_{-14})$; 
\item  let $\phi_{f_T}(Y_t^N):=e^{\II \langle Y_t^N,f_T\rangle}$, with 
$f_T$ compactly supported function in $C^2$.
Then 
\BEQ \esper[\phi_{f_T}(Y_T)\big| {\cal F}_t]=\phi_{h_t}(Y_t) \ 
 \exp\left( \half
\int_t^T  \left[ \II (1-\frac{\beta}{2}) \langle X_s,f''_s\rangle -
\langle X_s,(f'_s)^2\rangle \right] ds  \right)  \label{eq:intro-Gaussian} \EEQ
where $(f_s)_{0\le s\le T}$ is the solution of the asymptotic equation (\ref{eq:intro-PDE-f-asymptotic}).
 
\end{enumerate} }

\medskip\noindent {\em Scheme of proof.}
As in Israelsson \cite{Isr}, the main task is to prove a uniform in N Sobolev bound,
called "$H_8$-bound", see (\ref{eq:bd-H4}),
\BEQ \esper[\sup_{0\le s\le T} |\langle Y_s^N,\phi\rangle|] \le C_T ||\phi||_{H_8}
\EEQ
implying in particular tightness in some Sobolev space with negative index. Representing $\phi$ in terms of its
standard Stieltjes decomposition of order 5, $\phi={\cal C}^5 h$, this is  shown (by
technical arguments developed in \cite{Isr})
to hold provided

\BEQ   \esper[|N(M_t^N(z)-M_t(z))|^2]\le C |b|^{-12} \EEQ
(see (\ref{eq:4.14bis}))
or equivalently $\esper[|\langle Y_t^N,{\mathfrak{f}}_z\rangle|^2]  \le C |b|^{-12}$. Apply
(\ref{eq:intro-PDE-f}): start from terminal condition $f_T:={\mathfrak{f}}_{z_T}$ and
integrate in time, $\langle Y_T^N,{\mathfrak{f}}_{z_T}\rangle=\Big\langle Y_0^N,f_0 +
\half(1-\frac{\beta}{2}) \int_0^T dt\, \Big(\langle X_t^N,f_t''\rangle + \frac{1}{\sqrt{N}}\sum_{i=1}^N
f'_t(\lambda_t^i)dW^i_t \Big) \Big\rangle$. Terms in the r.-h.s. are bounded in section 4 using
a control over $(f_t)_{0\le t\le T}$, solution of the evolution equation (\ref{eq:intro-PDE-f}). The above equation is solved  in the following way: it is proved to be compatible
with the Stieltjes decomposition of order $\kappa$, 
\BEQ f_t(x)\equiv ({\cal C}^{\kappa} h_t)(x)=\int_{-\infty}^{+\infty} da \int_{-b_{max}}^{b_{max}} db\, (-\II b) \frac{|b|^{\kappa}}{(1+\kappa)!} \,  \,  {\mathfrak{f}}_z(x)
\, h_t(a,b)  \EEQ 
see Definition \ref{def:Stieltjes-decomposition}, if $\frac{\partial h_t}{\partial t}(a,b)=
{\cal H}_t h_t(a,b)$ for a certain time-dependent operator ${\cal H}_t$ -- a
"Stieltjes transform" of the evolution operator featuring in (\ref{eq:intro-PDE-f}) --  acting on 
$L^1(\Pi_{b_{max}})$, which is written down explicitly and analyzed in great details
in section 3. 


\bigskip


The article is organized as follows. We first introduce a family of {\em Stieltjes
decompositions} ${\cal C}^{\kappa}$ depending on a regularity index $\kappa=0,1,2,\ldots$ (see section 2). The main technical section is section 3, where we rewrite the
r.-h.s. of 
eq. (\ref{eq:PDE-f}) using Stieltjes decompositions as a sum of
linear operators which we call {\em generators}; these are of  two types:
generalized transport operators, including $V$-dependent terms sketched above in  (\ref{eq:generalized-transport-operator}), summing up to ${\cal H}_{transport}$, and bounded operators summing up to ${\cal H}_{nonlocal}$. We prove our Main Theorem in section
4. Since (\ref{eq:intro-Gaussian}) is formally just a consequence of It\^o's formula, and most of the technical arguments used in Israelsson's paper to justify this formula hardly
depend on $V$,   section 4 really revolves around a fundamental estimate,
Lemma \ref{lem:fundamental-estimate}, which is based on properties of the characteristics, hence is strongly $V$-dependent. The analysis of the generators
made in details in section 3 allows one to prove the latter estimate. To conclude, one uses as input {\em large deviation
estimates for the support of the measure} proved in section 5.  Finally, sections 6 and 7 are appendices, where we collected some
well-known facts and formulas about transport equations and Stieltjes transforms.

\medskip\noindent In an article in preparation \cite{Unt-law}, we solve (\ref{eq:intro-Gaussian}) and obtain the Gaussian kernel of the limiting fluctuation process $Y$. 


\section{Stieltjes decompositions}


In (Israelsson \cite{Isr}, Lemma 9) one finds the following decomposition of an $L^1$ function $f:\R\to\R$ living in the Sobolev space $H_2$
as a sum of functions of the type 
\BEQ {\mathfrak{f}}_z: x\mapsto \frac{1}{x-z}, \EEQ
 where $z=a+\II b$, $b\not=0$ (see section \ref{sec:Stieltjes}),

\BEQ f(x)=\int_{-\infty}^{+\infty} da  \, \int_{-\infty}^{+\infty} (-\II b) \,db\,   {\mathfrak{f}}_z(x) h(a), \qquad
h(a)=-f''(a). \label{eq:Stieltjes-decomposition0} \EEQ

The above "reproducing kernel type" decomposition is clearly not unique. The proof is based on the fact that, for $\kappa=0,1,2,\ldots$ (see (\ref{eq:8.11}))
\BEA \int_{-\infty}^{+\infty}  db\, (-\II b)\,  |b|^{\kappa}\,   \, \cdot\, {\cal F}({\mathfrak{f}}_{\II b})(s) &=& 2 \int_0^{+\infty} db\, |b|^{1+\kappa} {\cal F}(\Im {\mathfrak{f}}_{\II b})(s) \nonumber\\ &=&2\pi\int_{-\infty}^{+\infty} db\, |b|^{1+\kappa}\, e^{-b|s|}= 2\pi \, \cdot\, (1+\kappa)!  \
|s|^{-2-\kappa}, \label{eq:q} \nonumber\\
 \EEA
where $\cal F$ is the Fourier transform (see section 7 for normalization), from which we also  get the
following family of decompositions, valid for $f\in L^1\cap H_{2+\kappa}$, $\kappa\in\N$,
\BEQ f(x)=\int_{-\infty}^{+\infty} da  \, \int_{-\infty}^{+\infty} (-\II b)\, db\, \frac{|b|^{\kappa}}{(1+\kappa)!}\,  {\mathfrak{f}}_z(x) h(a), \ 
h(a)=  {\cal F}^{-1}(|s|^{2+\kappa} {\cal F}(f))(a), \label{eq:2Nd-Stieltjes-decomposition0} \EEQ
a straightforward generalization of (\ref{eq:Stieltjes-decomposition0}) obtained by
choosing some arbitrary value of $\kappa$ instead of $\kappa=0$ in (\ref{eq:q}).
Note that $h$ is {\em real} since ${\cal F}^{-1}(|s|\, \cdot)$ is given by a 
{\em real-valued} convolution kernel (see section \ref{sec:Stieltjes}).

The reason for introducing this $\kappa$-dependent family of decompositions is  that  the coefficient
of ${\mathfrak{f}}_z(x)$ now {\em vanishes to order $1+\kappa$ instead of $1$ on the real axis}, a property
inherited from the assumed supplementary regularity of $f$.
 Note that, for
$\kappa$ {\em even}, ${\cal F}^{-1}(|s|^{2+\kappa}{\cal F}(\cdot))$ is the differential operator $(-\partial_s^2)^{1+\kappa/2}$. For $\kappa$ {\em odd}, on the other hand, one gets
derivatives of the $(\frac{1}{x})$-kernel (see section \ref{sec:Stieltjes}).

\medskip\noindent 
Since all interesting phenomena appear when $|b|$ is small, and we  want
to avoid artificial problems arising when $|b|$ is not bounded, we shall
actually
use
analogous decompositions in which $|b|$ ranges from $0$ to some maximal value $b_{max}>0$. 
This introduces the following changes. First, instead of (\ref{eq:Stieltjes-decomposition0}), we get
\BEQ f(x)=\int_{-\infty}^{+\infty} da  \, \int_{-b_{max}}^{b_{max}} (-\II b)\, db\, \frac{|b|^{\kappa}}{(1+\kappa)!}\,  {\mathfrak{f}}_z(x) h(a), \ 
h(a)= ({\cal F}^{-1}(K^{\kappa}_{b_{max}})\ast f)(a), \label{eq:2Nd-Stieltjes-decomposition} \EEQ
where
\BEQ K^{\kappa}_{b_{max}}(s):=\left( 2\int_{0}^{b_{max}}  db\, |b|^{1+\kappa}\,   \, \cdot\, {\cal F}(\Im({\mathfrak{f}}_{\II b}))(s) \right)^{-1}  \label{eq:2.6} \EEQ
(note that the above integral is $>0$ by (\ref{eq:8.11})).
We now study the convolution operator
\BEQ  {\cal K}^{\kappa}_{b_{max}}: f\mapsto {\cal F}^{-1}(K^{\kappa}_{b_{max}})\ast f, \label{eq:calK-kappa} \EEQ  
depending on the parity
of $\kappa$:

\begin{itemize}
\item[(i)] For $\kappa$ {\em even},
\BEQ  \int_{0}^{b_{max}}  db\, |b|^{1+\kappa}\,   \, \cdot\, {\cal F}(\Im {\mathfrak{f}}_{\II b})(s) =\pi\int_{0}^{b_{max}} db\, |b|^{1+\kappa}\, e^{-b|s|}=\pi (-\partial_s^2)^{\kappa/2}
(k^0_{b_{max}}(|s|)), \label{eq:bmax-q-even} \EEQ
where
\BEQ k^0_{b_{max}}(|s|)=\frac{1}{s^2}\left(1-(1+b_{max}|s|)e^{-b_{max}|s|}\right).\EEQ
When $|s|\to \infty$, $k^0_{b_{max}}(|s|)\sim s^{-2}$; on the other hand, $k^0_{b_{max}}(|s|)
\overset{s\to 0}{\sim}  b_{max}^2 \sum_{k\ge 0} \frac{(-1)^k}{(k+2)!} (k+1) (b_{max}|s|)^{k}$.
Thus $(-\partial_s^2)^{\kappa/2}(k^0_{b_{max}}(|s|))\overset{|s|\to\infty}{\sim} (-1)^{\kappa/2} \frac{s^{-(2+\kappa)}}{(\kappa+1)!} $, and 
$(-\partial_s^2)^{\kappa/2}(k^0_{b_{max}}(|s|))\overset{|s|\to 0}{\sim}  
(-1)^{\kappa/2}\frac{b_{max}^{2+\kappa}}{2+\kappa}.$ 
It is a simple exercise to prove the following: let
\BEQ \underline{K}^{\kappa}_{b_{max}}(s):=(-1)^{\kappa/2} \left( (\kappa+1)!\  s^{2+\kappa} +  (2+\kappa)  b_{max}^{-(2+\kappa)} \right)^{-1} K^{\kappa}_{b_{max}}(s)-1. \label{eq:underlineK-kappa-even}
\EEQ
Then   $(\underline{K}^{\kappa}_{b_{max}})^{(j)}(s)$, $j=0,1,2$  is  $O(\frac{b^j_{max}}{1+b_{max}^2 s^2})$ uniformly in $s$ and $b_{max}$.   Hence {\em the convolution operator
\BEQ \underline{\cal K}^{\kappa}_{b_{max}}: f\mapsto {\cal F}^{-1}(\underline{K}^{\kappa}_{b_{max}})\ast f\EEQ
is a bounded operator from  $L^1(\R)\cap L^{\infty}(\R)$ to $L^1(\R)\cap L^{\infty}(\R)$.} Indeed, $$||| \underline{\cal K}^{\kappa}_{b_{max}}|||_{L^1(\R)\to L^1(\R)},\ 
||| \underline{\cal K}^{\kappa}_{b_{max}}|||_{L^{\infty}(\R)\to L^{\infty}(\R)} \le || {\cal F}^{-1}(\underline{K}^{\kappa}_{b_{max}})||_{L^1(\R)}$$ and
\BEA |{\cal F}^{-1}(\underline{K}^{\kappa}_{b_{max}})(x)| &\le&  \min\left(
 ||\underline{K}^{\kappa}_{b_{max}}||_{L^1}, \frac{1}{x^2} ||(\underline{K}^{\kappa}_{b_{max}})''||_{L^1}\right) \nonumber\\
 &=&
O\left( \int \frac{ds}{1+b_{max}^2 s^2} \right) \ \cdot\ \min(1,(\frac{b_{max}}{x})^2)  \nonumber\\
&=& \frac{1}{b_{max}}  \ \cdot\ \min(1,(\frac{b_{max}}{x})^2) ,
\label{eq:2.11} 
\EEA
from which $||| \underline{\cal K}^{\kappa}_{b_{max}}|||_{L^1(\R)\to L^1(\R)},
||| \underline{\cal K}^{\kappa}_{b_{max}}|||_{L^{\infty}(\R)\to L^{\infty}(\R)}=O(1).$
We may therefore write
\BEQ {\cal K}^{\kappa}_{b_{max}}=(1+\underline{\cal K}^{\kappa}_{b_{max}})\left(-(\kappa+1)! \ \partial_x^{2+\kappa}+  (-1)^{\kappa/2} (2+\kappa)b_{max}^{-(2+\kappa)} \right). \label{eq:calK-kappa-even} \EEQ

\item[(ii)] For $\kappa$ {\em odd},
\BEQ  \int_{0}^{b_{max}}  db\, |b|^{1+\kappa}\,   \, \cdot\, {\cal F}(\Im {\mathfrak{f}}_{\II b})(s) =\pi\int_0^{+\infty} db\, |b|^{1+\kappa}\, e^{-b|s|}=\pi (-\partial_s^2)^{(\kappa+1)/2}
(k^1_{b_{max}}(|s|)), \label{eq:bmax-q-odd} \EEQ
where
\BEQ k^1_{b_{max}}(|s|)=\frac{1}{|s|} (1-e^{-b_{max}|s|}).\EEQ
When $|s|\to \infty$, $k^1_{b_{max}(|s|)}\sim |s|^{-1}$; on the other hand, $k^1_{b_{max}}(|s|)
\overset{s\to 0}{\sim}  b_{max} \sum_{k\ge 0} \frac{(-1)^k}{(k+1)!}  (b_{max}|s|)^{k}$.
Thus $(-\partial_s^2)^{(\kappa+1)/2}(k^1_{b_{max}}(|s|))\overset{|s|\to\infty}{\sim} (-1)^{(\kappa+1)/2} \frac{|s|^{-1} s^{-(1+\kappa)}}{(\kappa+1)!} $, and 
$(-\partial_s^2)^{(\kappa+1)/2}(k^1_{b_{max}}(|s|))\overset{|s|\to 0}{\sim}  
(-1)^{(\kappa+1)/2}\frac{b_{max}^{2+\kappa}}{2+\kappa}.$ 
Therefore the previous analysis, from (\ref{eq:underlineK-kappa-even}) to the line before
(\ref{eq:calK-kappa-even}), remains valid, with $(-1)^{\kappa/2}$, resp. $s^{2+\kappa}$, replaced with
$(-1)^{(\kappa+1)/2}$, resp. $|s|s^{1+\kappa}$, and one obtains using (\ref{eq:C1}):
\BEQ {\cal K}^{\kappa}_{b_{max}}=(1+\underline{\cal K}^{\kappa}_{b_{max}})\left(-(\kappa+1)! \ \II H \partial_x^{2+\kappa}+ (-1)^{(\kappa+1)/2} (2+\kappa)b_{max}^{-(2+\kappa)} \right) \label{eq:calK-kappa-odd} \EEQ
where $H$ is the Hilbert transform, defined by the principal value integral
\BEQ Hf(x):=\frac{1}{\pi} p.v. \int_{-\infty}^{+\infty} \frac{1}{x-y}f(y)\, dy \EEQ
and $\underline{\cal K}^{\kappa}_{b_{max}}$, defined  by analogy with (i) as the
convolution with respect to the inverse Fourier transform of 
\BEQ \underline{K}^{\kappa}_{b_{max}}(s)=(-1)^{(\kappa+1)/2}
\left( (\kappa+1)!\  |s|^{2+\kappa} + (2+\kappa) b_{max}^{-(2+\kappa)} \right)^{-1}
K^{\kappa}_{b_{max}}(s)-1, \EEQ
 has operator norm
$O(1)$ on $L^1(\R)$ and on $L^{\infty}(\R)$. 
\end{itemize}

\medskip

The above formulas are  particular instances of {\em Stieltjes decompositions},
{\em where $h=h(a,b)$ is allowed to be complex-valued and to depend on $b$}.

\begin{Definition}[upper half-plane]
\begin{enumerate}
\item  Let $\Pi^+:=\{z\in\C\ |\ \Im(z)>0\}$.
\item For $b_{max}>0$, let $\Pi^+_{b_{max}}:=\{z\in \C\ |\ 0<\Im (z)<b_{max}\}$.
\item Let $\Pi^-:=-\Pi^+$, $\Pi^-_{b_{max}}:=-\Pi^+_{b_{max}}$ and $\Pi:=\Pi^+ \uplus \Pi^-,\ \Pi_{b_{max}}:=\Pi^+_{b_{max}}\uplus \Pi^-_{b_{max}}$.
\end{enumerate}
\end{Definition}

\begin{Definition} Let, for $p\in[1,+\infty]$ and $b_{max}>0$,
\BEQ L^p(\Pi_{b_{max}}):=\{h:\Pi_{b_{max}}\to\C\ |\  h(\bar{z})=\overline{h(z)}\ 
 (z\in \Pi^+_{b_{max}}) \, {\mathrm{\ and\ }} ||h||_{L^p(\Pi_{b_{max}})}<\infty\}, \EEQ
 where 
\BEQ ||h||_{L^p(\Pi_{b_{max}})}:= \left(\int_{-\infty}^{+\infty} da\,  \int_{-b_{max}}^{b_{max}} db\,  |h(a,b)|^p\right)^{1/p} \ (p<\infty), \qquad
||h||_{L^{\infty}(\Pi_{b_{max}})}:=\sup_{z\in\Pi_{b_{max}}} |h(z)|. \EEQ
\end{Definition}

We will be mostly interested in the extreme cases $p=1$, $p=\infty$. 
Letting formally $b_{max}\to +\infty$ one obtains in an obvious way the space $L^p(\Pi)$ with its norm $||\ ||_{L^p(\Pi)}$. However, we shall actually fix
some finite value of $b_{max}$, say (for convenience only), $0\le b_{max}\le \half$, implying: $\ln(1/|b|)\ge \ln 2>0$.

\begin{Definition}[{\bf Stieltjes decomposition}] \label{def:Stieltjes-decomposition}
Let $\kappa=0,1,2,\ldots$ 

\begin{enumerate}
\item
Let $h\in L^1(\Pi_{b_{max}})$. We say that $f:\R\to\R$ {\em has  Stieltjes decomposition $h$ of order $\kappa$ and cut-off $b_{max}$ on $[-R,R]$
} if, for all $|x|\le R$, 
\BEQ f(x)=({\cal C}^{\kappa}h)(x):=\int_{-\infty}^{+\infty} da \, \int_{-b_{max}}^{b_{max}} (-\II b)\, db\, \frac{|b|^{\kappa}}{(1+\kappa)!}\,  \  
{\mathfrak{f}}_z(x) h(a,b). \label{eq:Stieltjes-decomposition}  \EEQ
\item
The function $h:(a,b)\mapsto  {\cal K}^{\kappa}_{b_{max}}(f)(a)$
is called the {\em standard  Stieltjes decomposition of order $\kappa$} and {\em cut-off} $b_{max}$ of $f$.
\end{enumerate}
\end{Definition}

Thanks to the symmetry condition, $h(\bar{z})=\overline{h(z)}$, 
(\ref{eq:Stieltjes-decomposition})  may be rewritten in the form
\BEQ ({\cal C}^{\kappa}h)(x)=2 \int_{-\infty}^{+\infty} da \, \int_{0}^{b_{max}}  db\, \frac{|b|^{1+\kappa}}{(1+\kappa)!}\,  \  
\Im\left[{\mathfrak{f}}_z(x) h(a,b)\right],   \EEQ
from which it is apparent that $f$ is indeed real-valued.

\bigskip\noindent As already emphasized before, Stieltjes decompositions are not unique. In fact, it
turns out to be useful to introduce a larger family of decompositions depending on a 
further scale parameter, $\rho>0$. We shall give less details since we apply these to
the functions $\Im {\mathfrak{f}}_{z_T}$ with  $\Im z_T>0$  only, and concentrate
on the case $\kappa=0$ for computations. For $0<\rho\le b_{max}$, we write
\BEQ K^{\kappa}_{b_{max},\rho}(s):=
\left( \int_0^{b_{max}}  db\, b^{1+\kappa}\,  e^{-b/\rho}  \, \cdot\, {\cal F}(\Im({\mathfrak{f}}_{\II b}))(s) \right)^{-1}\EEQ
(compare with (\ref{eq:2.6}))  and let as before
\BEQ  {\cal K}^{0}_{b_{max},\rho}: f\mapsto {\cal F}^{-1}(K^{\kappa}_{b_{max},\rho})\ast f, \label{eq:calK-rho} \EEQ  
so that, for $b_T\not=0$, 
\BEQ \Im[{\mathfrak{f}}_{\II b_T}](x)=\int_{-\infty}^{+\infty} da\, \int_{0}^{b_{max}} db \, b^{\kappa+1}\, 
e^{-b/\rho} \, \Im[{\mathfrak{f}}_z(x)]\, {\cal K}^0_{b_{max},\rho}(\Im {\mathfrak{f}}_{\II b_T})(a).\EEQ

Then (compare with (\ref{eq:bmax-q-even})) $(K^{0}_{b_{max},\rho})^{-1}=\pi k^0_{b_{max}}(\frac{1}{\rho}+|s|)$. Thus (emphasizing only the differences with
$K^0_{b_{max}}=\lim_{\rho\to +\infty} K^0_{b_{max},\rho}$) 
$k^0_{b_{max},\rho}(|s|)\overset{s\to 0}{\sim} C\rho^2$ instead of $b_{max}^2$,
with $C:=1-(1+b_{max}/\rho)e^{-b_{max}/\rho}$ bounded away from $0$. Hence

\BEA \underline{K}^0_{b_{max},\rho}(s) &:=& (|s|+\rho^{-1})^{-2} K^0_{b_{max},\rho}-1 \nonumber\\
&=& - \frac{(1+b_{max}(|s|+\rho^{-1}))e^{-b_{max}(|s|+\rho^{-1})}}{1-(1+b_{max}(|s|+\rho^{-1}))e^{-b_{max}(|s|+\rho^{-1})}} \EEA
is a $O(e^{-\half b_{max}(|s|+\rho^{-1})})$, defining a bounded operator on $L^1\cap L^{\infty}$ with
operator norm
\BEQ O\left( \int ds\, e^{-\half b_{max}(|s|+\rho^{-1})}\right) \, \cdot\,
\int dx\, \min(1,\frac{b_{max}}{x})^2 
=O(e^{-\half b_{max}/\rho}) \EEQ
(compare with (\ref{eq:2.11})), so that ${\cal K}^0_{b_{max},\rho}=
{\cal F}^{-1}((|s|+\rho^{-1})^2)\ast$, times (1+plus  bounded perturbation).  More generally,
it may be proved that (for some constant $c_{\kappa}$)  ${\cal K}^{\kappa}_{b_{max},\rho}=c_{\kappa}
{\cal F}^{-1}((|s|+\rho^{-1})^{2+\kappa})\ast$, times (1+  bounded perturbation). On the other hand, 
 letting $f= \Im[{\mathfrak{f}}_{z_T}]$, with
 $z_T=a_T+\II b_T$, $b_T>0$, we find using (\ref{eq:8.4}), (\ref{eq:8.11})
\BEQ {\cal F}^{-1}((|s|+\rho^{-1})^{2+\kappa} {\cal F} \Im[{\mathfrak{f}}_{z_T}])(a)=
 \left( \rho^{-1} -\frac{\partial}{\partial b_T}\right)^{2+\kappa}  \left(\frac{b_T}{(a-a_T)^2+|b_T|^2} \right) \label{eq:fz=Ch} \EEQ
 
Hence, letting 
\BEQ h(a,b):= e^{-b/\rho} {\cal K}^{\kappa}_{b_{max},\rho}(\Im {\mathfrak{f}}_{z_T})(a), \label{eq:2.29-0}
\EEQ
 we get:

\BEA  \int_{-\infty}^{+\infty} da\, \int_{-b_{max}}^{b_{max}} db\, h(a,b) &= &  O(1) \rho \left\{
\frac{1}{\rho^{2+\kappa}} + \frac{1}{|b_T|^{2+\kappa}} \right\} \nonumber\\
&=& O\left(
\frac{1}{\rho^{1+\kappa}}(1+ O( (\frac{\rho}{b_T})^{2+\kappa} ) ) \right)
\EEA
which is minimal, of order
\BEQ ||h||_{L^1(\Pi_{b_{max}})}\approx b_T^{-1-\kappa},  \label{eq:hL1b-1} \EEQ
 when $\rho\approx b_T$. Choosing $\rho=b_T/C$ for some
 large enough absolute constant $C>0$, we further obtain -- specifically in the case $\kappa=0$   -- a {\em positive} function $h$,
 which can hence be interpreted as a {\em density}. Also, one easily checks that, still
 with $\rho\approx b_T$,
\BEQ ||(a,b)\mapsto \ln(1/|b|) h(a,b)||_{L^1(\Pi_{b_{max}})}\approx \ln(1/b_T)b_T^{-1-\kappa} 
\label{eq:log-estimate} \EEQ
if $0<b_T\le \half$.  On the other hand,
 \BEQ ||h||_{L^{\infty}(\Pi_{b_{max}})}\approx b_T^{-3-\kappa} \label{eq:hLinftyb-3}. \EEQ

\medskip\noindent 
In section 4, we consider the time-evolution of ${\cal C}^{\kappa}h_T$, where $h_T$
is essentially as in (\ref{eq:2.29}),
\BEQ h_T(a,b):= e^{-b/\rho} {\cal K}^{\kappa}_{b_{max},\rho}(\chi_R\, \Im {\mathfrak{f}}_{z_T})(a), \label{eq:2.29}
\EEQ
for some cut-off function (see Definition \ref{def:cut-off} below) essentially supported
on the ball $B(0,R)$ for some fixed radius $R>0$. An easy adaptation of the above arguments, and a use of
(\ref{eq:bound-vp2}) when $\kappa$ is odd in order to deal with the Hilbert transform,
show that the above estimates (\ref{eq:hL1b-1},\ref{eq:log-estimate},\ref{eq:hLinftyb-3}) remain correct, while now $h_T$ is $O(1)$, independently of $b_T$, far from
the support, e.g. on $\Big( (B(0,2R))^c\times [-b_{max},b_{max}] \Big) \cup
\Big( \R\times ([-b_{max},b_{max}]\setminus[-\half b_{max},\half b_{max}]) \Big)$.


\section{Generators}


The general purpose of the section is the following: for $\kappa=0,1,2,\ldots$ {\em
fixed}, we want to write down an explicit 
time-dependent operator ${\cal H}(t)$   such that the right-hand side of
(\ref{eq:PDE-f}) for $f_t$ decomposed as 
\BEQ  f_t(x)=\int_{-\infty}^{+\infty} da  \, \int_{-b_{max}}^{b_{max}}  db\,  (-\II b)\,  \frac{|b|^{\kappa}}{(1+\kappa)!} \,  \  f_{z}(x) \, h_t(a,b) \label{eq:generator} \EEQ
(see Definition \ref{def:Stieltjes-decomposition}) is equal to  
\BEQ \int_{-\infty}^{+\infty} da\, \int_{-b_{max}}^{b_{max}}  db\, (-\II b)\,  \frac{|b|^{\kappa}}{(1+\kappa)!}\, 
 {\mathfrak{f}}_z(x)\,  {\cal H}(h_t)(t;a,b) \EEQ
 where ${\cal H}(h_t)(t;a,b)\equiv ({\cal H}(t)(h_t))(a,b).$

\noindent Given the characteristic evolution in the $z$-coordinate found by Israelsson 
(recalled in \S 3.1 below) -- which one may view as a {\em deterministic}
Markov process -- it is natural to think of the function $h$ as  a density $h(a,b)da\, db$, and then to interpret ${\cal H}(t)$ as a {\em Fokker-Planck operator}, whence (by
duality) ${\cal L}(t):=({\cal H}(t))^{\dagger}$ as the generator of a random process (see
section \ref{sec:transport} for more). However, contrary
to the harmonic case studied by Israelsson and Bender, in general we obtain a truly
{\em random} process, furthermore, a {\em signed} process, with $h$ a signed function. 
For lack of references on these notions, we shall refrain from developing this signed Markov process
interpretation, and solve instead the evolution equation
\BEQ \frac{dh_t}{dt}(a,b)={\cal H}(h_t)(t;a,b)\EEQ
using {\em semi-group theory}. 

\medskip
\noindent We have been voluntarily been vague up to this point about the double 
dependency of ${\cal H}$ on the integer index $\kappa$ and the cut-off scale $b_{max}$.
Why couldn't one just set $\kappa=0$ and let $b_{max}\to +\infty$, as does Israelsson ?
The reason is, we cannot handle properly the potential-dependent part of the generator
(save when $V$ is order $\le 2$, which is the case considered in \cite{Isr}) without 
introducing various cut-offs and perturbative arguments -- unless maybe if $V$ is analytic 
(or even better, polynomial), where another strategy is perhaps possible. Since we do not want to make this assumption, we shall:

\begin{itemize}
\item[(1)] ({\em support issues}) in practice replace $V(x)$ by its {\em Taylor expansion to order 2 around all points in the support of the measures $(X_t^N)_{0\le t\le T}$ and
$(X_t)_{0\le t\le T}$} and {\em treat the Taylor remainder as a perturbation}. Since
$V$ is not  bounded at infinity, it is important {\em not to} Taylor expand
around any point on the real line, but only on a {\em compact interval}; choosing
as compact interval the support  is natural because the singular kernel
term in (\ref{eq:PDE-f}) vanishes outside. {\em For brevity, we shall
henceforth call {\em support of the measure} the random set $\left(\cup_{0\le t\le T}
\ \supp(X_t^N) \right) \cup \left( \cup_{0\le t\le T}\  \supp(X_t)\right)$ and
denote by $R>0$ any number such that the support of the measure is $\subset[-R,R]$.} We rely on the  bounds developed in
section \ref{sec:entropy} to argue that the probability of the support not to be
included in $[-R,R]$ for some large enough $R$ is exponentially small in $N$ when
$R$ is large enough.  Then we naturally decompose $h\in L^1(\Pi_{b_{max}})$ as $h^{int}+h^{ext}$ where
$\supp(h^{int})\subset[-3R,3R]\times[-b_{max},b_{max}]$ and $\supp(h^{ext})\subset (\R\setminus[-2R,2R])
\times[-b_{max},b_{max}]$, 
for instance by writing $h(a,b)=\bar{\chi}_R(a) h(a,b)+(1-\bar{\chi}_R(a))h(a,b)$,
where $\bar{\chi}_R:\R\to[0,1]$ is some smooth function such that $\supp(\bar{\chi}_R)\subset[-3R,3R]$
and $\supp(1-\bar{\chi}_R)\subset\R\setminus[-2R,2R]$.  The action of  ${\cal H}$ on
$h^{ext}$ is very simple and can be added to  the action of the remainder term
${\cal H}_{nonlocal}$
discussed in (2).

\item[(2)] ({\em varying $\kappa$}) in order to be able to treat the part  (thereafter denoted by ${\cal H}_{nonlocal}$) of the generator due to the remainder term as a perturbation, it is
important to see that ${\cal H}_{nonlocal}h(a,b)=O(|b|)$ when $b\to 0$. This being the case,
we may also consider ${\cal H}_{nonlocal}$ as an operator {\em intertwining
a Stieltjes decomposition of order $\kappa$ with a Stieltjes decomposition of order $\kappa+1$}, leading to a modification of the scheme developed around (\ref{eq:generator}): namely,
we want the right-hand side of
(\ref{eq:PDE-f}) for $f_t$ decomposed as  (\ref{eq:generator}) for some integer index
 $\kappa$ to be equal to  
\BEQ \int_{-\infty}^{+\infty} da\, \int_{-b_{max}}^{b_{max}} (-\II b)\,  db\, \frac{|b|^{1+\kappa}}{(2+\kappa)!}\, 
 {\mathfrak{f}}_z(x)\,  {\cal H}^{\kappa+1,\kappa}(h_t)(t;a,b). \EEQ
Thus, instead of a single operator ${\cal H}^0$, we  deal simultaneously
with a family of operators ${\cal H}^{\kappa}$ and a family of intertwining operators
${\cal H}_{nonlocal}^{\kappa+1,\kappa}$, for $\kappa=0,1,2,\ldots.$   These 
intertwinings result in an expansion of the Green kernel explicited in (\ref{eq:2Nd-Green-function-expansion}). This strategy yields optimal bounds in section 4.

\end{itemize}

\noindent Let us collect by anticipation all the terms which will come out of the
computations in \S 3.1 through \S 3.9.  As a general rule, if ${\cal H}: L^1([-3R,3R]\times
[-b_{max},b_{max}])\oplus L^1( (\R\setminus[-2R,2R])\times
[-b_{max},b_{max}]) \to  L^1([-3R,3R]\times
[-b_{max},b_{max}])\oplus L^1( (\R\setminus[-2R,2R])\times
[-b_{max},b_{max}])$ is an (unbounded) operator, we denote by
\BEQ {\cal H}^{int}:={\cal H}\big|_{L^1([-3R,3R]\times
[-b_{max},b_{max}])\oplus 0}, \qquad 
{\cal H}^{ext}:={\cal H}\big|_{0 \oplus L^1([-3R,3R]\times
[-b_{max},b_{max}])}
\EEQ
its restrictions to either factor, and ${\cal H}^{(int,int)},{\cal H}^{(int,ext)},{\cal H}^{(ext,int)},{\cal H}^{(ext,ext)}$ its four block-components, so that
\BEQ {\cal H}=\left( \begin{array}{ccc} {\cal H}^{int} & | & {\cal H}^{ext} \end{array}
\right)=\left( \begin{array}{cc} {\cal H}^{(int,int)} & {\cal H}^{(ext,int)} \\
{\cal H}^{(int,ext)} & {\cal H}^{(ext,ext)} \end{array} \right). \EEQ

By explicit computation we show that
\BEQ {\cal H}^{\kappa}=\left(\begin{array}{cc} {\cal H}^{\kappa}_{transport} & 0 \\
0 & 0 \end{array}\right) 
+{\cal H}^{\kappa}_{nonlocal}; \EEQ


\BEQ {\cal H}^{\kappa}_{transport}:={\cal H}^{\kappa}_0+\sum_{k=0}^2 {\cal H}_{pot}^{\kappa,(k)} \EEQ 
 is a sum of 4 transport
operators -- with ${\cal H}^{\kappa}_0$ coming from the  $\left(\frac{1}{x}\right)$-kernel
 part (see \S 3.1), and ${\cal H}_{pot}^{\kappa,(k)}$, $k=0,1,2$ 
 (see \S 3.3, \S 3.4, \S 3.5)
from the Taylor expansion of the potential --,  which are {\em unbounded operators}; 

\BEQ 
{\cal H}^{\kappa}_{nonlocal}=
 \left(\begin{array}{ccc} {\cal H}_{pot}^{\kappa,(3)}  + {\cal H}^{\kappa}_{\mathrm{bdry}} 
 & | 
 & {\cal H}_{pot}^{\kappa,ext} \end{array}\right),  \label{eq:H-dec3} \EEQ
 is a sum of 
{\em bounded} operators. The operator ${\cal H}^{\kappa,(3)}_{pot}$ -- coming
from the third order Taylor remainder for the potential -- is introduced in 
\S 3.7. The operator $ {\cal H}_{\mathrm{bdry}}\equiv 
{\cal H}^{\kappa}_{h-bdry}+{\cal H}^{\kappa}_{v-bdry}$ is itself
a sum of boundary terms (see \S 3.9): contributions coming from horizontal boundary 
$[-3R,3R]\times \{\pm b_{max}\}$, collected in ${\cal H}^{\kappa}_{\mathrm{h-bdry}}$, and
contributions coming from vertical boundary 
$\{\pm 3R\}\times [-b_{max},b_{max}]$, collected in ${\cal H}^{\kappa}_{\mathrm{v-bdry}}$.  The eight operators ${\cal H}_0^{\kappa},
{\cal H}_{pot}^{\kappa,(0)}, {\cal H}_{pot}^{\kappa,(1)}, {\cal H}_{pot}^{\kappa,(2)},
{\cal H}_{pot}^{\kappa,(3)}, {\cal H}^{\kappa}_{\mathrm{h-bdry}}, {\cal H}^{\kappa}_{\mathrm{v-bdry}} 
, {\cal H}^{\kappa,ext}_{pot}$ are
defined resp. in\\ (\ref{eq:H0},\ref{eq:HV0},\ref{eq:HV1},\ref{eq:HV2},\ref{eq:HV3},\ref{eq:Hhbdry}, \ref{eq:Hvbdry},
\ref{eq:Hext}).

\medskip\noindent We also write down expressions for ${\cal H}^{\kappa+1,\kappa}_{nonlocal}$, namely, ${\cal H}_{pot}^{\kappa+1,\kappa,(3)},  {\cal H}^{\kappa+1,\kappa}_{\mathrm{h-bdry}}, {\cal H}^{\kappa+1,\kappa}_{\mathrm{v-bdry}} 
, {\cal H}^{\kappa+1,\kappa,ext}_{pot}$, to be found resp. in (\ref{eq:HV3+1},
\ref{eq:Hhbdry+1},\ref{eq:Hvbdry+1},\ref{eq:Hext+1}).

\medskip\noindent The kernels of these operators are denoted by the letter $g$, for
instance,
\BEQ {\cal H}_0^{\kappa}(h)(a,b)=\int_{-\infty}^{+\infty} da_T \int_{-b_{max}}^{b_{max}}
db_T\ g_0^{\kappa}(a,b;a_T,b_T)h(a_T,b_T) \EEQ
and similarly for the other operators.

\medskip\noindent Let us simply state as as general remark that the dependence on $R$
of the bounds of the present section will never be made explicit, since for the proof
of our Main Theorem (see section 4), we shall simply fix $R=R(T)$, where $R(T)$ is
a fixed radius depending only on the ptoential and on the time horizon $T$, defined
in section 5.

\bigskip

\subsection{The $\left(\frac{1}{x}\right)$-kernel part}

\medskip

(\ref{eq:PDE-f}) is easily solved by the characteristic method in the  case $V\equiv 0$
for test functions $f$ of the form $f(x)=\frac{c}{x-z}$. Up to conjugation
we may assume that $b:=\Im z>0$. This is done in (Israelsson \cite{Isr}, Lemmas 2-4) -- we need only subtract the trivial contribution of the
harmonic potential --\ :

\begin{Proposition}(see Israelsson \cite{Isr}) \label{prop:lemma-2}
Assume $V\equiv 0$. Then eq. (\ref{eq:PDE-f}) with terminal condition $f_T(x)=\frac{c}{x-z}$ $(\Im z>0)$ is solved as
\BEQ f_t^T(x)=\frac{C_t}{x-Z_t},\EEQ
where $(C_t)_{0\le t\le T}, (Z_t)_{0\le t\le  T}$ solve the following ode's,
\BEQ \frac{dZ_t}{dt}=-\frac{\beta}{4} (M_t^N(Z_t)+M_t(Z_t)), \qquad 
\frac{dC_t}{dt}=-\frac{\beta}{4} ((M_t^N)'(Z_t)+M_t'(Z_t)) C_t  \label{eq:dc/dt} \EEQ
with terminal conditions $Z_T=z,C_T=c$.
In particular,
\BEQ \Im \frac{dZ_t}{dt}=-\frac{\beta}{4} \langle X_t^N+X_t, \Im({\mathfrak{f}}_{z_T})\rangle \le 0.
\label{eq:Im-dz/dt<0} \EEQ
\end{Proposition}

Obviously, the solution of eq. (\ref{eq:PDE-f}) with terminal condition $f_T(x)=
\frac{c}{x-\bar{z}}$ is now $f_t^T(x)= \frac{\bar{C}_t}{x-\bar{Z}_t}$.

\medskip

The last inequality (\ref{eq:Im-dz/dt<0}), a simple consequence of (\ref{eq:Im-fz>0}),
implies that $Z_t$ moves {\em away} from the real axis as $t$ decreases, hence away
from singularities. From the above and from general properties of the Stieltjes transform of $\rho_t$ (see section \ref{sec:Stieltjes}), one can deduce bounds for the displacement, as in
\cite{Isr}. First  
 $|M_{t}^N(z)|$, $|M_{t}(z)|\le 1/b_t$,  hence (letting
 $Z_t=:A_t+\II B_t$, $B_t>0$), for some large enough constant $C>0$,
 
\BEQ B_T\le B_t\le \sqrt{|B_T|^2+C(T-t)}. \EEQ
 Similarly, $|A_t-A_T|=O\left(\frac{T-t}{B_T}\right)$.
Finally,
$\frac{\beta}{4}\left(|(M_t^N)'(Z_t)|+|(M_t)'(Z_t)|\right)\le \frac{|dB_t/dt|}{B_t}=\left|\frac{d}{dt} (\ln(B_t))\right|$,
whence $C_t\le \frac{B_t}{B_T}\le \sqrt{|B_T|^2+C(T-t)} \ \, / \, B_T.$
Summarizing: for a 
given final condition $z$, $|\frac{dA_t}{dt}|,|\frac{dB_t}{dt}|\le \frac{\beta}{2} \frac{1}{B_t}\le \frac{\beta}{2} \frac{1}{b}$ is
 bounded along the characteristics, but may become arbitrarily large when
 $b\to 0$; $|B_t-B_T|\le \sqrt{C(T-t)}$ is bounded independently of $b$, while $|A_t-A_T|$ is {\em not}. On the other hand, starting from $Z_T$ far enough from the support
 of $(X^N_t)_{0\le t\le T},(X_T)_{0\le t\le T}$, e.g. $|\Re Z_T|\ge CR$, where $C>1$  and (by assumption) $\supp(X^N_t),\supp(X_t)\subset B(0,R)$, then (by
 (\ref{eq:2Nd-bound-Mt})) $|A_t|>C' R$ $(1<C'<C)$ for all $t\in[0,T]$, whence $
 |\dot{B}_t|,|\dot{A}_t|\le \frac{\beta}{R}
 <\infty$, provided $T<\frac{R^2}{\beta}(C-C')$.

\begin{Definition}
Let ${\cal L}_0={\cal L}_0(t)$ be the time-dependent operator defined by
\BEA &&  ({\cal L}_{0}(t) \phi)(a,b)=-\frac{\beta}{4} \left( \Re \left[(M_{t}^N+M_{t})(z)\right] \partial_a\phi(a,b) +\Im \left[(M_{t}^N+M_{t})(z)\right] \partial_b\phi(a,b)  \right. \nonumber\\
&&  \qquad \qquad \qquad \qquad \left.  + ((M_{t}^N)'+M_{t}')(z) \phi(a,b)\right)  \label{eq:L0} \EEA
\end{Definition}

The motivation for this definition is the following. Let $f_t^T(x)=\frac{C_t}{x-Z_t}$ as in Proposition \ref{prop:lemma-2}. Then $\frac{\partial f_t^T}{\partial t}\big|_{t=T}=({\cal L}_0 f_T)(T;z):=({\cal L}_0(T) f_T)(z)$. Take, more generally, a terminal condition for
(\ref{eq:PDE-f}) of
the form
\BEQ f_T(x)=\int da_T \int db_T \, \frac{h_T(a_T,b_T)}{x-z_T}.\EEQ
Then
\BEA \frac{\partial f_t^T}{\partial t}\big|_{t=T} &=& \int da_T db_T\,  ({\cal L}_0( \frac{1}{x - \cdot})(z_T)  h_T(a_T,b_T) \nonumber\\
&=&  \int da_T db_T\, {\mathfrak{f}}_{z_T}(x) ({\cal L}_0^{\dagger} h_T)(T,z_T),\EEA
where $({\cal L}_0^{\dagger}h_T)(T;z_T)=(({\cal L}_0(T))^{\dagger}h_T)(z_T)$
is obtained from the adjoint of ${\cal L}_0(T)$. Thus ${\cal L}_0$, resp.
${\cal L}_0^{\dagger}$ may be considered as the generator of a -- here
{\em deterministic} -- {\em generalized} Markov process $Z_{\cdot}=A_{\cdot}+\II B_{\cdot}$,
resp. the associated Fokker-Planck generator, where {\em generalized} refers
to the supplementary order $0$ term in ${\cal L}_0(t)$ (second line of (\ref{eq:L0})), which has on top of that the nasty property of not being
even real-valued.

\bigskip
We now need some very general development, which we apply to ${\cal L}_0$ in this paragraph.
{\em Let $w:\Pi\to\R_+^*$ be some weight function.}
 The relation $\frac{d}{dt} (wh)= 
 w ({\cal L}_0^w(t))^{\dagger}h$ defines an operator $({\cal L}_0^w(t))^{\dagger}$,
\BEA ({\cal L}_0^w)^{\dagger}(t;a,b)&:=&(w {\cal L}_0(t) w^{-1})^{\dagger}(a,b) \nonumber\\
&=& -\frac{\beta}{4}  w(a,b)^{-1}  \left(-\partial_a\  \Re \left[(M_t^N+M_{t})(z)\right]  - \partial_b\  \Im \left[(M_{t}^N+M_{t})(z)\right] \right. \nonumber\\
&&  \qquad \qquad \qquad \qquad \left.  + ((M_{t}^N)'+M_{t}')(z) \right)   w(a,b)
\label{eq:w} \EEA
which is the adjoint of  ${\cal L}_0$ with respect to the measure $w(a,b)da\, db$ on $\Pi$. In other words,  we see that the
solution $(h_t)_{0\le t\le T}$ of the equation 
\BEQ \frac{\partial h_t}{\partial t}=({\cal L}_0^w)^{\dagger}(t) h_t \EEQ
with terminal condition $h_T$ is the {\em density} of $Z_{\cdot}$ with respect to the
measure $w(a,b) da\, db$. 

\medskip\noindent
Let us consider specifically the cases 
\BEQ w(a,b):=b\, |b|^{\kappa}, \qquad \kappa=0,1,2,\ldots\EEQ
for which we write
\BEQ {\cal L}_0^w\equiv {\cal L}_0^{\kappa}.\EEQ
These cases allow a direct connection to Stieltjes decompositions of order $\kappa$ with $b_{max}=+\infty$, namely: 
{\em if 
\BEQ f_T(x)={\cal C}^{\kappa}h_T(x)=\int da \int_0^{+\infty}   db\, (-\II b)\, \frac{|b|^{\kappa}}{(1+\kappa)!}\,
h_T(a,b)\ {\mathfrak{f}}_z(x) \EEQ
 then 
\BEQ f_t(x):={\cal C}^{\kappa} h_t(x),\EEQ
where $(h_t)_{0\le t\le T}$ is the  solution of 
\BEQ \frac{\partial h_t}{\partial t}=({\cal L}_0^{\kappa})^{\dagger}(t) h_t,
\label{eq:dagger} \EEQ
 }

\medskip
\noindent  The generator
\BEQ {\cal L}_0^{\kappa}= \left((b\, |b|^{\kappa})^{-1}  {\cal L}_0^{\dagger}
b\, |b|^{\kappa} \right)^{\dagger}=|b|^{1+\kappa}  {\cal L}_0 (|b|^{1+\kappa})^{-1} 
\label{eq:sharp} \EEQ
  is obtained from (\ref{eq:L0}) by replacing the
multiplicative term $ \frac{\beta}{4}((M_{t}^N)'+M_{t}')(z) \phi(a,b)$ with
$\frac{\beta}{4} \left( ((M_{t}^N)'+M_{t}')(z)-\frac{1+\kappa}{b}\Im[(M_{t}^N+M_{t})(z)] \right) \phi(t;a,b)$, from which

\BEA && {\cal H}_0^{\kappa}(h)(t;a,b) :=({\cal L}_0^{\kappa})^{\dagger}(t)(h_t)(a,b)
\nonumber\\ && \qquad \qquad =\frac{\beta}{4}
\left(
\partial_a\left[   \Re \left((M_{t}^N+M_{t})(z)\right) h(t;a,b)\right] +
\partial_b\left[ \Im \left((M_{t}^N+M_{t})(z)\right) h(t;a,b) \right] \right. \nonumber\\
&&  \qquad \qquad  \qquad \left.  + \left[
\frac{1+\kappa}{b} \Im\left((M_{t}^N+M_{t})(z)\right) -  ((M_{t}^N)'+M_{t}')(z)
 \right]  h(t;a,b)\right)  \label{eq:H0} \EEA

Consider the extended characteristics $(z_t,c_t):=(a_t+\II b_t,c_t^{\kappa})$ of 
the operator ${\cal H}_0^{\kappa}$ as in section \ref{sec:transport}:
 they
are defined as the solution of 

\BEQ \frac{da_t}{dt}=\frac{\beta}{4} \Re(M_{t}^N+M_{t})(a_t+\II b_t),  
\qquad  \frac{db_t}{dt}=\frac{\beta}{4} \Im(M_{t}^N+M_{t})(a_t+\II b_t)
\EEQ

\BEQ
\frac{dc^{\kappa}_{t}}{dt}=\frac{\beta}{4} \left(
\frac{1+\kappa}{b_t} \Im\left[(M_{t}^N+M_{t})(a_t+\II b_t)\right] +
\left((\bar{M}_{t}^N)'+\bar{M}_{t}'\right)(a_t+\II b_t)   \right) c^{\kappa}_{t}  \label{eq:char-L0} \EEQ

Here we used the fact that $\partial_a\, \Re (M_t^N+M_t)(z)=\partial_b\,  \Im
(M_t^N+M_t)(z)=\Re ((M_t^N)'+M'_t)(z)$ since $z\mapsto (M_t^N +M_t)(z)$ is holomorphic. {\em Mind the sign changes} with respect to Proposition 
\ref{prop:lemma-2}: characteristics of the dual operator ${\cal H}_0^{\kappa}$
have  reversed velocities with respect to those of ${\cal L}$, with characteristic curves $(Z_t)_{0\le t\le T}$ running
backwards in time (see section \ref{sec:transport} for more details). By convention, 
characteristics are killed upon touching the real axis.

Now it follows from (\ref{eq:M'M}) that {\em $\Re (-\tau(t,z_t))\equiv \Re\left( -(c^{\kappa}_t)^{-1} 
\frac{dc_{t}^{\kappa}}{dt} \right)\le 0$ for $\kappa\ge 0$, whence (see section \ref{sec:transport})
${\cal H}_0^{\kappa}$  is a generator
of a semi-group of $L^{\infty}$-contractions. Because the first line of
(\ref{eq:H0}) is in divergence form, and the second line has positive
real part, ${\cal H}_0^{\kappa}$ is also the generator of a semi-group of   $L^1$-contractions. }

\bigskip\noindent

\bigskip
\noindent Let us now see what happens for $R$ and $b_{max}$ {\em finite}. It is clear
that the cut-off does not change the characteristic equations for $(a_t+\II b_t, c_t)$. On the other hand, we get two supplementary 
{\em boundary terms}. This can be proved as follows. Assume  
\BEQ f_T(x)={\cal C}^{\kappa} h_T(x)=\int_{-3R}^{3R} da_T \int_{-b_{max}}^{b_{max}} db_T\,
 (-\II b_T)\, \frac{|b|^{\kappa}}{(1+\kappa)!}\,
h_T(a_T,b_T)\ {\mathfrak{f}}_{z_T}(x). \EEQ
Then (coming back directly to the characteristics equations of Proposition
\ref{prop:lemma-2})
\BEA &&-\frac{\beta}{4} \int \int \frac{f'_T(x)-f'_T(y)}{x-y}( X^N_t(dy)+ X_t(dy))
\nonumber\\
&&\ = -\frac{\beta}{4} \int_{-3R}^{+3R} da_T \int_{-b_{max}}^{b_{max}}
(-\II b_T)\,  db_T\,
\frac{|b_T|^{\kappa}}{(1+\kappa)!}\,  
 \left\{ \Im[(M^N_T+M_T)(z_T)]\,  (\partial_{b_T} {\mathfrak{f}}_{z_T})(x) +
 \right. \nonumber\\
&& \qquad \qquad \left.  \Re[(M^N_T+M_T)(z_T)]\,  (\partial_{a_T} {\mathfrak{f}}_{z_T})(x) \right\} \  h_T(a_T,b_T)+\cdots \nonumber\\
&&\  = \frac{\beta}{4} \int_{-3R}^{+3R} da_T \int_{-b_{max}}^{b_{max}}
(-\II b_T)\,  db_T\,  
\frac{|b_T|^{\kappa}}{(1+\kappa)!}\,
\left\{ \partial_{b_T}\left( \Im[(M^N_T+M_T)(z_T)]\, h_T(a_T,b_T)\right) +
\right.\nonumber\\
&& \qquad\qquad\qquad \left.
\partial_{a_T}\left( \Re[(M^N_T+M_T)(z_T)]\, h_T(a_T,b_T)\right) \right\} \ {\mathfrak{f}}_{z_T}(x)
\ \ 
+\cdots\nonumber\\
&& \qquad\qquad + {\mathrm{bdry}}
\EEA
where "$\cdots$" denote the contribution of the $c$-characteristics
(which we can ignore), so
(by integration by parts) we get a boundary term bdry $\equiv$ h-bdry$_0+$
 v-bdry$_0$ on the support of the measure, with
\BEA &&  {\mathrm{h\!-\!bdry}}_0=-\frac{\beta}{4} \frac{b_{max}^{1+\kappa}}{(1+\kappa)!}  \int_{-\infty}^{+\infty} da_T
\,  \ \cdot \chi_R(x)\cdot\   \left( \Im[(M^N_T+M_T)(a_T+\II b_{max})] \ \cdot
\right.\nonumber\\ && \qquad \left. \cdot\  f_{a_T+\II b_{max}}(x) h(a_T,b_{max}) 
-\Im[(M^N_T+M_T)(a_T-\II b_{max})] f_{a_T-\II b_{max}}(x) h(a_T,-b_{max}) \right)  \label{eq:h-bdry0} \nonumber\\
 \EEA
and
\BEA && {\mathrm{v\!-\!bdry}}_0=-\frac{\beta}{4} \int_{-b_{max}}^{b_{max}} (-\II b_T)\,  db_T \, 
\frac{|b_T|^{\kappa}}{(1+\kappa)!}   \ \cdot  \chi_R(x)\ \cdot \  \left\{\Re [(M_T^N+M_T)(3R+\II b_T)] \ \cdot\ \right.\nonumber\\
&&\qquad \left.\cdot\   f_{3R+\II b_T}(x)\,  h(3R,b_T) -
\Re [(M_T^N+M_T)(-3R+\II b_T)] f_{-3R+\II b_T}(x)\,  h(-3R,b_T) \right\} \nonumber\\
\label{eq:v-bdry0}  
\EEA


\subsection{The potential-dependent part: general introduction}


Consider now the potential-dependent part in (\ref{eq:PDE-f}). 
{\em Without further mention we fix for the discussion $R\ge 1$ and some arbitrary $b_{max}\in(0,\half]$.} Generally speaking we want to {\em write the action of the operator $V'(x)\frac{\partial}{\partial x}$ on
a function $f\equiv f_T$ with  Stieltjes decomposition of order $\kappa$} on $[-R,R]$
\BEQ f(x)=\int_{-\infty}^{+\infty} da_T \, \int_{-b_{max}}^{b_{max}} 
 db_T\,  (-\II b_T)\,\frac{|b_T|^{\kappa}}{(1+\kappa)!}\,  \  
{\mathfrak{f}}_{z_T}(x) h(a_T,b_T), \qquad |x|\le R \label{eq:3.30} \EEQ
(see Definition \ref{def:Stieltjes-decomposition}) .  In principle (see below though)
it may be done in the following way.

\begin{Definition}[cut-offs] \label{def:cut-off}
Let $\chi_R:\R\to\R$ be a smooth cut-off function such that: 
\begin{itemize}
\item[(i)] $\chi_R\equiv 1$ on $[-R,R]$;
\item[(ii)] $\chi_R\big|_{B(0,\frac{3}{2} R)^c}\equiv 0$.
\end{itemize}
and $\bar{\chi}_R$ be the function  $x\mapsto \chi_R(\frac{x}{2})$.
\end{Definition}

\begin{Definition}[$g$-kernel] \label{def:g-kernel}
 Let, for $\kappa,\kappa'=0,1,2,\ldots$
\BEQ g^{\kappa';\kappa}_{pot}(a,b;a_T,b_T):={\bf 1}_{|b|<b_{max}} \frac{(-\II b_T) |b_T|^{\kappa}}{(1+\kappa)!} \, {\cal K}^{\kappa'}_{b_{max}} \left(x\mapsto \chi_R(x)V'(x) f'_{z_T}(x) \right)(a) \label{eq:g-kernel} \EEQ
with ${\cal K}^{\kappa'}_{b_{max}}$ as in (\ref{eq:calK-kappa}).
\end{Definition}

In practice we are only interested in couples of indices $(\kappa,\kappa'=\kappa)$ and $(\kappa,\kappa'=\kappa+1)$, and
let $g^{\kappa}_{pot}\equiv g^{\kappa;\kappa}_{pot}$. 
Let $f$ be a function as in (\ref{eq:3.30}).
Then, {\em using the standard  order $\kappa'$   Stieltjes decomposition with
cut-off $b_{max}$  of $V'(x)\frac{\partial}{\partial x}({\mathfrak{f}}_{z_T}(x))$}, we get  
\BEA &&  V'(x)\frac{\partial}{\partial x} f(x)=\int_{-\infty}^{+\infty} da  \, \int_{-b_{max}}^{b_{max}}   db\,  (-\II b)\, \frac{|b|^{\kappa'}}{(1+\kappa')!} \,  \  {\mathfrak{f}}_z(x) \nonumber\\
&& \qquad\qquad  \int_{-\infty}^{+\infty} da_T \, \int_{-b_{max}}^{b_{max}} db_T \,    g^{\kappa';\kappa}_{pot}(a,b;a_T,b_T) h(a_T,b_T),\EEA
thus defining (unbounded) operators ${\cal H}^{\kappa}_{pot},{\cal H}^{\kappa+1;\kappa}_{pot}: L^1(\Pi_{b_{max}},da_T \, db_T)\to L^1(\Pi_{b_{max}},da\, db)$,
$$ {\cal H}^{\kappa}_{pot}(h)(a,b)=\int_{-\infty}^{+\infty} da_T\, \int_{-b_{max}}^{b_{max}} db_T\, g^{\kappa}_{pot}(a,b;a_T,b_T)h(a_T,b_T) $$
$$ {\cal H}^{\kappa+1;\kappa}_{pot}(h)(a,b)=\int_{-\infty}^{+\infty} da_T\, \int_{-b_{max}}^{b_{max}} db_T\, g^{\kappa+1;\kappa}_{pot}(a,b;a_T,b_T)h(a_T,b_T).$$

\medskip\noindent However, these Stieltjes representations of the vector field $V'(x)\frac{\partial}{\partial x}$ 
will be  used directly only  when $h=0\oplus h^{ext}$ has support in $\R\setminus[-2R,2R]$,
producing the ${\cal H}^{ext}$-term. When $h=h^{int}\oplus 0$ has support 
$\subset[-3R,3R]$, we separate first the Taylor
expansion of order 2 of $V'(x)$ around $a_T$, as explained in (\ref{eq:intro-W})
or (\ref{eq:W}) below, which is directly analyzed without further Stieltjes
decomposition.

\bigskip\noindent
Let us discuss in details  how we proceed when $h=h^{int}$.
As we have just mentioned, the first step is to use a second-order Taylor-expansion of  $V'$ around $a_T$ for
$|a_T|\le 3R$ and $x\in$ supp$(\chi_R)$, i.e. $|x|\le \frac{3}{2}R$,
\BEQ  V'(x)=V'(a_T)+V''(a_T)(x-a_T)+ V'''(a_T) \frac{(x-a_T)^2}{2}+ (x-a_T)^3  W_{a_T}(x-a_T), \EEQ
where $W_{a_T}$ is $C^7$. Thus 
\BEQ V'(x)\frac{\partial}{\partial x}=V'(a_T)\partial_x+V''(a_T)(x-a_T)
\partial_x + V'''(a_T)\frac{(x-a_T)^2}{2} \partial_x+ (x-a_T)^3 W_{a_T}(x-a_T) \partial_x  \label{eq:W} \EEQ
 
makes  four  different contributions to the generator, ${\cal H}_{pot}^{\kappa}$, $\kappa=0,1,2,3$, resp. called {\em constant, linear, quadratic and
Taylor remainder term}, plus some boundary contributions. Computations
show the following: ${\cal H}_{pot}^{\kappa,(k)}$, $k=0,1,2$, are directly of the adjoint form 
\BEQ h\mapsto \Big( (a_T,b_T)\mapsto 
\partial_{a_T}(v_{hor}(a_T,b_T) h(a_T,b_T))+\partial_{b_T}(v_{vert}(a_T,b_T) h(a_T,b_T)) - \tau(a_T,b_T) h(a_T,b_T) \Big) \EEQ
with $\Re \tau(\cdot)\le 0$, implying that {\em they generate $L^1$-contractions} (see section \ref{sec:transport}, and recall that we go {\em backwards} in time). 
As a nice feature of this problem, 
${\cal H}_{pot}^{\kappa,(k)}$, $k=0,1,2$ {\em also generate $L^{\infty}$-contractions}. 
 Replacing the
vector field $\chi_R(x)V'(x)\partial_x$ in (\ref{eq:g-kernel}) above by $\chi_R(x) (x-a_T)^3
W_{a_T}(x-a_T)\partial_x$ produces a kernel $g^{(3)}_{pot}(a,b;a_T,b_T)$ discussed
in \S 3.8. As for the first three terms, they are directly shown to be equivalent to
the action of a transport operator (see \S 3.3, 3.4, 3.5). We sum up in \S 3.6 the contributions
of the transport operators introduced in \S 3.1, 3.3, 3.4 and  3.5. The operators
${\cal H}^{\kappa,(3)}_{pot}$ and ${\cal H}^{\kappa,ext}_{pot}$ are studied in \S 3.7 and 3.8. Finally,
the contribution of the boundary terms is analyzed in \S 3.9.

\bigskip\noindent The terms in ${\cal H}_{nonlocal}$, on the other hand, {\em do not} generate
neither $L^{\infty}$- nor $L^1$-contractions. 
Thanks to the horizontal and vertical cut-offs however, they are bounded,
hence generate by integration 
some exponentially increasing time factor $e^{C_R t}$, with $C_R$ depending on $R$.


\bigskip

\subsection{Constant term}

 Inserting the constant operator $V'(a_T)\frac{\partial}{\partial x}$ inside the  Stieltjes decomposition $$f_T(x)=\int_{-3R}^{+3R} da_T \, \int_{-b_{max}}^{b_{max}}(-\II b_T)\,  db_T\, |b_T|^{\kappa} \, 
{\mathfrak{f}}_{z_T}(x) h(a_T,b_T),$$  we get 
\BEA && \int_{-3R}^{3R} da_T\ \int_{-b_{max}}^{b_{max}} (-\II b_T)\, db_T\, |b_T|^{\kappa}  V'(a_T) 
f'_{z_T}(x)  h(a_T,b_T) \nonumber\\
&& \qquad =- \int_{-3R}^{3R} da_T\ \int_{-b_{max}}^{b_{max}} (-\II b_T)\, db_T\, |b_T|^{\kappa}\,  
V'(a_T)   \frac{\partial}{\partial a_T} 
( {\mathfrak{f}}_{z_T}(x))  h(a_T,b_T) \nonumber\\
&& \qquad =\int_{-3R}^{3R} da_T\ \int_{-b_{max}}^{b_{max}} (-\II b_T)\, db_T\, |b_T|^{\kappa}\,   {\mathfrak{f}}_{z_T}(x) \frac{\partial}{\partial a_T}
\left( V'(a_T)h(a_T,b_T) \right) + {\mathrm{bdry}}. \nonumber\\
\EEA

where
\BEA && {\mathrm{bdry}}\equiv{\mathrm{v\!-\!bdry}}_{pot}^{(0)}=- \int_{-b_{max}}^{b_{max}} (-\II b_T)\,  db_T \, |b_T|^{\kappa}  \ \cdot  \chi_R(x)\ \cdot \nonumber\\
&& \qquad \cdot \,  \left(V'(3R) f_{3R+\II b_T}(x)\,  h(3R,b_T) -
V'(-3R) f_{-3R+\II b_T}(x)\,  h(-3R,b_T) \right)
\label{eq:v-bdry-0}  
\EEA
is a vertical boundary term coming from integration by parts, which we shall not discuss till \S 3.9.
 
This defines a new operator ${\cal H}^{(0)}_{pot}$ in divergence form,
\BEQ {\cal H}^{\kappa,(0)}_{pot}(h)(a_T,b_T)=\frac{\partial}{\partial a_T}
\left( V'(a_T)h(a_T,b_T) \right), \label{eq:HV0} \EEQ
with corresponding extended characteristics 
\BEQ \frac{da_t}{dt}=V'(a_t), \qquad \frac{db_t}{dt}=0, \qquad 
\frac{dc^{\kappa}_{t}}{dt}=V''(a_t). \label{eq:char-LV0} \EEQ


\subsection{Linear term}   Proceeding as in the {\em constant term} case, we
insert the operator $V''(a_T)(x-a_T)\frac{\partial}{\partial x}$ inside the Stieltjes decomposition
of $f_T$, getting
\BEA && \int da_T \int_{-b_{max}}^{b_{max}} (-\II b_T)\, db_T \, |b_T|^{\kappa}\,   V''(a_T) \left[ (x-a_T)\partial_x
{\mathfrak{f}}_{z_T}(x)\right] h(a_T,b_T) \nonumber\\
&& \qquad  =-  \int da_T \int_{-b_{max}}^{b_{max}} (-\II b_T)\, db_T \, |b_T|^{\kappa}\,   V''(a_T)  \ \cdot\partial_{b_T}(b_T {\mathfrak{f}}_{z_T}(x)) \ h(a_T,b_T) \nonumber\\
&& \qquad = {\mathrm{bdry}} \, +\ \int da_T \int_{-b_{max}}^{b_{max}} (-\II b_T)\, db_T \, |b_T|^{\kappa}\,  {\mathfrak{f}}_{z_T}(x) \,  V''(a_T) 
(b_T\partial_{b_T}+1+\kappa)
h(a_T,b_T), \nonumber\\
\EEA
where
\BEA && {\mathrm{bdry}}\equiv{\mathrm{h-bdry}}_{pot}^{(1)}=-b_{max}^{2+\kappa} \int da_T \, V''(a_T) \ \cdot  \chi_R(x)\ \cdot \nonumber\\
&& \qquad \cdot \,  \left(f_{a_T+\II b_{max}}(x)\,  h(a_T,b_{max}) -
f_{a_T-\II b_{max}}(x)\,  h(a_T,-b_{max}) \right)
\label{eq:h-bdry-1}
\EEA
is a horizontal boundary term coming from integration by parts.
This defines an operator, ${\cal H}_{pot}^{\kappa,(1)}$,
\BEQ {\cal H}_{pot}^{\kappa,(1)}(h)(a_T,b_T)
=\left(\partial_{b_T} b_T+\kappa\right)(V''(a_T)h(a_T,b_T)).
\label{eq:HV1} \EEQ
with associated characteristics,
\BEQ \frac{da_t}{dt}=0, \qquad  \frac{db_t}{dt}=V''(a_t)b_t, \qquad 
\frac{dc^{\kappa}_{t}}{dt}= (1+\kappa) V''(a_t)c^{\kappa}_{t} \label{eq:char-LV1} \EEQ

\subsection{Quadratic term}


Proceeding as in the previous paragraph, we  insert the operator
$\half V'''(a_T) (x-a_T)^2 \frac{\partial}{\partial x}$ inside the Stieltjes decomposition. Since $(x-a_T)^2 \frac{\partial}{\partial x}=\partial_x (x-a_T)^2 -
2 (x-a_T)$, and
\BEQ \partial_x \left[ (x-a_T)^2 {\mathfrak{f}}_{z_T}(x)\right]=1+b_T^2 \partial_{a_T} {\mathfrak{f}}_{z_T}(x),
\qquad  -2(x-a_T){\mathfrak{f}}_{z_T}(x)=-2-2\II b_T {\mathfrak{f}}_{z_T}(x)  \label{eq:1-2} \EEQ
this term produces
a new operator

\BEQ {\cal H}_{pot}^{\kappa,(2)}(h)(a_T,b_T)=\half  \left(- \partial_{a_T} b_T^2 -2\II
b_T\right) (V'''(a_T) h(a_T,b_T)).  \label{eq:HV2} \EEQ

with associated  characteristics
\BEQ \frac{da_t}{dt}
=-\half V'''(a_t) b_t^2, 
 \qquad  \frac{db_t}{dt}=0, \qquad 
\frac{dc^{\kappa}_{t}}{dt}=-\II V'''(a_t) b_t\, c^{\kappa}_{t} \label{eq:char-LV2} \EEQ

plus a vertical boundary term,
\BEA && {\mathrm{bdry}}\equiv{\mathrm{v-bdry}}_{pot}^{(2)}=- \int_{-b_{max}}^{b_{max}} (-\II b_T)\,  db_T \, |b_T|^{\kappa+2}  \ \cdot  \chi_R(x)\ \cdot \nonumber\\
&& \qquad \cdot \,  \left(V'''(3R) f_{3R+\II b_T}(x)\,  h(3R,b_T) -
V'''(-3R) f_{-3R+\II b_T}(x)\,  h(-3R,b_T) \right)
\label{eq:v-bdry-2}
\EEA

The remaining term due to the constant $\half V'''(a_T) (1-2)$ in (\ref{eq:1-2}), integrated
with respect to the measure $|b_T|^{1+\kappa}  da_T\,  db_T$, adds a time-independent constant $C$ to $f_t$, 
which however disappears from the computations since: (i) the right-hand side of 
(\ref{eq:PDE-f})  features only $f'_t$; (ii) $d\langle Y_t^N,C\rangle=dC=0$ in
(\ref{eq:Ito}).

\medskip


\subsection{Recapitulating: the transport contribution}

We can now write down the  action of the sum of our three  transport operators,
\BEQ {\cal H}^{\kappa}_{transport}(t):={\cal H}^{\kappa}_0(t)+{\cal H}_{pot}^{\kappa,(0)}
+{\cal H}_{pot}^{\kappa,(1)}+{\cal H}_{pot}^{\kappa,(2)} \EEQ
(note that only ${\cal H}_0^{\kappa}$ depends on the time variable), defined
respectively in (\ref{eq:H0}), (\ref{eq:HV0}), (\ref{eq:HV1}) and (\ref{eq:HV2}).
Summing the contributions
in (\ref{eq:char-L0}), (\ref{eq:char-LV0}), (\ref{eq:char-LV1}) and (\ref{eq:char-LV2}), we obtain the following equation for the 
characteristics on $[-3R,3R]\times[-b_{max},b_{max}]$,
\BEQ \frac{da_t}{dt}=\frac{\beta}{4} \Re(M_{t}^N+M_{t})(a_t+\II b_t) +V'(a_t)-\half V'''(a_t) b_t^2   \label{eq:da(tau)/dtau} \EEQ
\BEQ \frac{db_t}{dt}=\frac{\beta}{4} \Im(M_{t}^N+M_{t})(a_t+\II b_t)+ V''(a_t)b_t \label{eq:db(tau)/dtau} \EEQ 
\BEA &&  \frac{dc^{\kappa}_{t}}{dt}= \left[\frac{\beta}{4} \left(
\frac{1+\kappa}{b_t}\Im(M_{t}^N+M_{t})(a_t+\II b_t) +  ((\bar{M}_{t}^N)'+(\bar{M}_{t})')(a_t+\II b_t)  \right) \right. \nonumber\\
&& \qquad\qquad\qquad \left. +(2+\kappa) V''(a_t) 
-\II V'''(a_t) b_t\right] c^{\kappa}_{t}. \label{eq:dEsharp(tau)/dtau} \EEA

\noindent Let us make the following observations
 (see \S 3.1):  
 
\begin{itemize}
\item[(i)]  $t\mapsto |b_t|$ increases, whence $|b_t|\le |b_T|$ $(0\le t\le T)$.

\item[(ii)] velocities are $O(1)$ far from the support, e.g. on $\{|a|>CR\}\cup
\{|b|>b_{max}/2\}$ (see discussion below Proposition \ref{prop:lemma-2});
\item[(iii)] as is already true of each individual generator ${\cal H}^{\kappa}_0$,
${\cal H}^{\kappa,(i)}_{pot}$ ($i=0,1,2$), the sum
${\cal H}^{\kappa}_{transport}(t)$ generates a semi-group  of $L^1$-contractions. The same holds if one
considers
$L^{\infty}$-
norms, since the real part of (\ref{eq:dEsharp(tau)/dtau}) is positive.
\end{itemize}

\noindent As a side remark, we may choose $T$ small enough so that characteristics
$(a_t+\II b_t)_{0\le t\le T}$ started at time $T$
on the boundary $\Big(\{\pm 3R\}\times[-b_{max},b_{max}]\Big)\cup\Big([-3R,3R]\times \{\pm
b_{max}\}\Big)$  always remains far from the support, in the sense of (ii).

\medskip\noindent Section \ref{sec:transport} therefore implies:

\begin{Lemma} \label{lem:U-transport}
Let $u_T\in (L^1\cap L^{\infty})([-3R,3R]\times(0,b_{max}])$ and $\kappa=0,1,2\ldots$ 
Then the backward evolution equation $\frac{du}{dt}={\cal H}^{\kappa}_{transport}(t)u(t), 
u\big|_{t=T}=u_T$ $(0\le t\le T)$ has a unique solution $u(t):=U_{transport}(t,T)u_T$, such that 
\BEQ ||u_t||_{L^1([-3R,3R]\times(0,b_{max}])}\le  ||u_T||_{L^1([-3R,3R]\times(0,b_{max}])} \EEQ
\BEQ ||u_t||_{L^{\infty}([-3R,3R]\times(0,b_{max}])}\le  ||u_T||_{L^{\infty}([-3R,3R]\times(0,b_{max}])} \EEQ
\end{Lemma}

\medskip\noindent In section 4, we will separate (\ref{eq:dEsharp(tau)/dtau}) into its
real and imaginary parts. Solving the $(a,b)$-characteristics coupled with
the real characteristic $\tilde{c}^{\kappa}$,
\BEA &&  \frac{d\tilde{c}^{\kappa}_{t}}{dt}= \Re \left[\frac{\beta}{4} \left(
\frac{1+\kappa}{b_t}\Im(M_{t}^N+M_{t})(a_t+\II b_t) +  ((\bar{M}_{t}^N)'+(\bar{M}_{t})')(a_t+\II b_t)  \right) \right. \nonumber\\
&& \qquad\qquad\qquad \left. +(2+\kappa) V''(a_t) 
-\II V'''(a_t) b_t\right] \tilde{c}^{\kappa}_{t}, \label{eq:ctilde} \EEA
one gets a backward evolution equation equation, 
\BEQ \frac{d\tilde{u}}{dt}=[\Re{\cal H}^{\kappa}_{transport}(t)]\tilde{u}(t),  \ \ 
\tilde{u}\big|_{t=T}=\tilde{u}_T\equiv u_T \qquad (0\le t\le T), \label{eq:utilde}
\EEQ
 solved
as $\tilde{u}(t):=\tilde{U}_{transport}(t,T)\tilde{u}_T$, such that
\BEQ ||\tilde{u}_t||_{L^1([-3R,3R]\times(0,b_{max}])}\le  ||\tilde{u}_T||_{L^1([-3R,3R]\times(0,b_{max}])} \EEQ
\BEQ ||\tilde{u}_t||_{L^{\infty}([-3R,3R]\times(0,b_{max}])}\le  ||\tilde{u}_T||_{L^{\infty}([-3R,3R]\times(0,b_{max}])} \EEQ
{\em Assuming $u_T\ge 0$, $\tilde{u}$ is simply the modulus of $u$,}
\BEQ \tilde{u}_t(a,b)=|u_t(a,b)|.\EEQ
In particular, $\tilde{u}_t\ge 0$. The {\em adjoint evolution} with
generator $(\Re {\cal H}^{\kappa}_{transport})^{\dagger}$ is
{\em sub-Markovian}, i.e. $\tilde{u}_t$ is the density at time $t$ of a
(deterministic)
time-reversed  Markov $(\tilde{A}_t,\tilde{B}_t)_{0\le t\le T}$
process (whose trajectories run backwards w.r. to $(a_t,b_t)$) with kernel $p(t,\tilde{a}_t,\tilde{b}_t;s,\tilde{a}_s,\tilde{b}_s)\ge 0$ ($t\le s$) such that
$\int d\tilde{a} \int d\tilde{b} \, p(t,\tilde{a},\tilde{b};s,\tilde{a}_s,\tilde{b}_s)\le 1$.  

\bigskip\noindent
{\em Remark.}
It is instructive to look at terms that would be produced by continuing the Taylor
expansion to infinity. Note that the boundary value of the  operator ${\cal H}^{\kappa,(2)}_{pot}(h)$ vanishes (${\cal H}_{pot}^{\kappa,(2)}(h)(a_T,b_T=0)\equiv 0$ vanishes to order $\ge 1$).
It may be proven in general that the contribution of order $j$ vanishes 
to order $\ge j-1$. Summing up the whole series would yield (up to
$\kappa$-dependent terms)
\BEQ \frac{da_t}{dt}=\sum_{n=2p} \frac{V^{(n+1)}(a_t)}{n!} (\II b_t)^n,
\qquad \frac{db_t}{dt}=-\II \sum_{n=2p+1} \frac{V^{(n+1)}(a_t)}{n!} (\II b_t)^n \EEQ
\BEQ \frac{dc_{t}}{dt}=\sum_{n\ge 1} V^{(n+1)}(a_t) \frac{(\II b_t)^{n-1}}{(n-1)!}.\EEQ
Note that, for $V$ {\em holomorphic}, this is equivalent to the vector field $V'(Z_{t})\partial_{Z_{t}}+V''(Z_{t})$. {\em However}, terms of order $\ge 3$ also produce
{\em polynomials} in $x$ (instead of linear combination of ${\mathfrak{f}}_z(x)$, $z\in\C\setminus\R$). Integrating the generator would yield power series in $x$ which
are maybe controllable, but this would require totally different techniques
with respect to ours.

\bigskip\noindent



\subsection{Main remainder term}


We study in this subsection the $g$-kernels $g_{pot}^{\kappa';\kappa,(3)}$ obtained
by replacing $\chi_R(x)V'(x)f'_{z_T}(x)$ with $\chi_R(x)(x-a_T)^3 W_{a_T}(x-a_T)f'_{z_T}(x)$
in Definition \ref{def:g-kernel}. We want to prove that the operators
\BEQ {\cal H}^{\kappa,(3)}_{pot}: h\mapsto {\cal H}^{\kappa,(3)}_{pot}(h): \left( (a,b)\mapsto\int_{-3R}^{3R} da_T
\int_{-b_{max}}^{b_{max}} db_T\, g_{pot}^{\kappa,(3)}(a,b;a_T,b_T) h(a_T,b_T)  \right)  \label{eq:HV3} \EEQ
and
\BEQ {\cal H}^{\kappa+1;\kappa,(3)}_{pot}: h\mapsto {\cal H}^{\kappa+1;\kappa,(3)}_{pot}(h): \left( (a,b)\mapsto\int_{-3R}^{3R} da_T
\int_{-b_{max}}^{b_{max}} db_T\, g_{pot}^{\kappa+1;\kappa,(3)}(a,b;a_T,b_T) h(a_T,b_T)  \right)  \label{eq:HV3+1} \EEQ
are bounded operators in $L^1(\Pi_{b_{max}})$.

\begin{Definition}
For $f\in C^k(\R,\R)$ and $R>0$, let
\BEQ ||f||_{k,[-R,R]}:=\sup_{0\le j\le k}  \sup_{[-R,R]} |f^{(j)}|.\EEQ
\end{Definition}

\begin{Lemma} \label{lem:S8}
\begin{itemize}
\item[(i)] $|||{\cal H}_{pot}^{\kappa+1;\kappa,(3)}|||_{L^1(\Pi_{b_{max}})\to
L^1(\Pi_{b_{max}})} = O(b_{max}^{-1} ||V'||_{8+\kappa,[-3R,3R]} )$;
\item[(ii)] $|||{\cal H}_{pot}^{\kappa,(3)}|||_{L^1(\Pi_{b_{max}})\to
L^1(\Pi_{b_{max}})} = O( ||V'||_{7+\kappa,[-3R,3R]} )$.
\end{itemize}
\medskip\noindent The same estimates hold for the {\em $L^{\infty}$-operator norms} $||| \ \cdot \
|||_{L^{\infty}(\Pi_{b_{max}})\to L^{\infty}(\Pi_{b_{max}})}$.
\end{Lemma}

\medskip\noindent  As can be checked by
looking at the details of the proof, norms
$|||\ \cdot\ |||_{L^1(\Pi_{b_{max}}) \to L^{\infty}(\Pi_{b_{max}})}$ are deduced
from these by dividing by a volume factor $\Vol([-3R,3R]\times[-b_{max},b_{max}])\approx b_{max}$, which follows from the fact that the kernel $g_{pot}^{\kappa';\kappa,(3)}$ is regular on the diagonal.

\bigskip\noindent We shall actually only prove the
 statement for $|||{\cal H}_{pot}^{\kappa+1;\kappa,(3)}|||_{L^1\to L^1}$. The proof of (ii) reduces
 immediately to (i) by substituting $\kappa\to \kappa-1$ and taking into account
 the extra $b$-prefactor.
 
\medskip \noindent {\bf Remark.} Had we Taylor expanded $V'$ to order 2 instead
of order 3, would we have obtained an {\em unbounded} operator instead of 
${\cal H}^{\kappa+1;\kappa,(3)}$, as testified by the bound (\ref{eq:3.88})
below (with $x^2$ instead of $x^3$ in the numerator, the integral would diverge
in the limit $b_T\to 0$ for $m=3+\kappa$). On the other hand,
for the same reason, it is easy to see by looking at the details of the
 proof that ${\cal H}^{\kappa,(3)}:L^1(\Pi_{b_{max}})\to
  L^{\infty}(\Pi_{b_{max}})$ is
bounded, typically because in (\ref{eq:3.88}) one obtains instead of the
$L^1$-kernel $x\mapsto \frac{b_T}{x^2+b_T^2}$ the bounded function
$x\mapsto \frac{b_T^2}{x^2+b_T^2}$. 
 
\medskip

 \noindent{\bf Proof.}  By definition, 
\BEQ  g^{\kappa+1;\kappa,(3)}_{pot}(a,b;a_T,b_T)= \frac{(-\II b_T) |b_T|^{\kappa}}{(1+\kappa)!}\,   {\cal K}^{\kappa+1}_{b_{max}} \left( 
x\mapsto \chi_R(x) (x-a_T)^3 W_{a_T}(x-a_T) f'_{z_T}(x) \right)(a) \nonumber\\ 
\label{eq:3.72}
 \EEQ

We consider in the computations only the part $g^{\kappa}$ of the kernel $g^{\kappa+1;\kappa,(3)}_{pot}$ obtained by replacing ${\cal K}^{\kappa}$ with the operators
in factor of the bounded operator $1+\underline{\cal K}^{\kappa}_{b_{max}}$ in
 (\ref{eq:calK-kappa-even}, \ref{eq:calK-kappa-odd}). For some
 numerical constants $c_1=c_1(\kappa),c_2=c_2(\kappa)$,
\BEA  g^{\kappa}(a,b;a_T,b_T) &=&(-\II b_T) |b_T|^{\kappa}  \left\{ c_1 [{\cal F}^{-1}
  (|s|^{3+\kappa}) \ast] +c_2 b_{max}^{-(3+\kappa)} \right\}   (F_{a_T+\II b_T})(x) \nonumber\\
&=:& (-\II b_T) |b_T|^{\kappa} (c_1 G^1_{a_T+\II b_T}(x)+c_2 G^2_{a_T+\II b_T}(x)) , \label{eq:H2} \EEA
with $x:=a-a_T$ and 
\BEQ F_{a_T+\II b_T}(y):= \tilde{W}_{a_T}(y)   \frac{y^3}{(y-\II b_T)^2} 
\label{eq:F2} \EEQ
where $\tilde{W}_{a_T}(y)=\chi_R(a_T+y)W_{a_T}(y)$ has support $\subset[-\frac{9}{2}R,\frac{9}{2}R]$ since
$|a_T|\le 3R$.
Thus (\ref{eq:HV3+1}) looks like a convolution formula in the coordinates $a,a_T$
-- but not quite since $\tilde{W}_{a_T}$ depends on $a_T$ --, which
leads us to use the following bound (where $G_{a_T+\II b_T}:=c_1 G^1_{a_T+\II b_T}
+ c_2 G^2_{a_T+\II b_T}$)

\BEA && ||{\cal H}^{\kappa+1;\kappa,(3)}(h)||_{L^1(\Pi_{b_{max}})} \le 
b_{max} \left( \sup_{(a_T,b_T)\in[-3R,3R]\times [-b_{max},b_{max}]} |b_T|^{1+\kappa} \right. \nonumber\\
&& \qquad\qquad\qquad\qquad \left.  \int da\,
|G_{a_T+\II b_T}(a-a_T)|\right)\ ||h||_{L^1(\Pi_{b_{max})}} 
\EEA
whence
\BEQ |||{\cal H}^{\kappa+1;\kappa,(3)}|||_{L^1\to L^1}
 \le  b_{max} \left(\sup_{(a_T,b_T)\in[-3R,3R]\times [-b_{max},b_{max}]} ||\ |b_T|^{1+\kappa} G_{a_T+\II b_T}||_{L^1} \right). 
\EEQ
 
We must therefore bound $||\ |b_T|^{1+\kappa} G_{a_T+\II b_T}||_{L^1}$. {\em The main
issue, on which we shall now concentrate, is to bound 
$|| \ |b_T|^{1+\kappa} G^1_{a_T+\II b_T}||_{L^1([-5R,5R])}$;} as shown later on,   the missing
terms are less
singular and may be bounded similarly. Thus from now on and till below (\ref{eq:3.99}),
$|a_T|\le 3R$, $|x|\le 5R$ and $|a|:=|x+a_T|\le 8R$ are bounded.

\medskip

\noindent If $\kappa$ is {\em odd} then $|s|^{3+\kappa}=s^{3+\kappa}$ is the Fourier
symbol of a differential operator, otherwise $|s|^{3+\kappa}=\sgn(s) s^{3+\kappa}$ involves
a further convolution by a singular kernel. Let us accordingly dinstinguish two cases. But first of all we let, for $\ell\ge 3$ and $|a_T|\le 3R$, 
\BEQ C_{\ell}:=\sup_{3\le \ell'\le \ell} 
 ||\tilde{W}_{a_T}^{(\ell'-3)}||_{\infty} \EEQ
 and note that
\BEQ C_{\ell}=O\left( \sup_{3\le \ell'\le \ell} \sup_{[-3R,3R]}
|V^{(\ell'+1)}| \right)= O( ||V'||_{\ell,[-3R,3R]})\EEQ

More generally, if $\ell\ge \ell'\ge \ell''\ge 3$ and $r\ge 0$,
\BEQ || y\mapsto (\tilde{W}_{a_T}^{(\ell''-3)}(y) y^r)^{(\ell'-\ell'')} ||_{\infty} = O(  C_{\ell}).\EEQ

\begin{itemize}
\item[(i)] ($\kappa$ odd)\qquad 
First 

\BEQ ||\ |b_T|^{1+\kappa} G^1_{a_T+\II b_T}||_{L^1(\R)}=O(C_{6+\kappa} |b_T|^{1+\kappa}) \sum_{m=0}^{3+\kappa} \int_{-5R}^{5R}  dx\, \,  \left|\partial_x^m\left(
\frac{x^3}{(x-\II b_T)^2}\right) \right|. \EEQ

Then 
\BEQ \partial_x^m\left(
\frac{x^3}{(x-\II b_T)^2}\right)= \sum_{p-q=1-m, p\le 3, q\ge 2} C^m_{p,q}
x^p (x-\II b_T)^{-q} \label{eq:mpq} \EEQ
for some coefficients $C^m_{p,q}$. If $3\le m\le 3+\kappa$ then
$|x^p (x-\II b_T)^{-q}|=O(\frac{|x-\II b_T|^{3-m}}{x^2+|b_T|^2})=O(\frac{|b_T|^{3-m}}{x^2+|b_T|^2})$. Multiplying with respect to $|b_T|^{1+\kappa}$ and integrating, one gets (using $b_T\le b_{max}\le 1$)
\BEA  \, ||\ |b_T|^{1+\kappa} \partial_x^m\left(
\frac{x^3}{(x-\II b_T)^2}\right)||_{L^1(\R)} &\le&  C'  |b_T|^{3+\kappa-m} 
\int dx\, \frac{b_T}{x^2+|b_T|^2}=O(1).\nonumber\\   \label{eq:3.88}
\EEA

If $m\le 2$ then $b_T\left|\partial_x^m \left(
\frac{x^3}{(x-\II b_T)^2}\right) \right|=O(|b_T|\frac{|x|^{3-m}}{x^2+|b_T|^2})=O(|x|^{2-m})$. Thus\\
  $ ||\ |b_T|^{1+\kappa} \partial_x^m\left(
\frac{x^3}{(x-\II b_T)^2}\right)||_{L^1([-5R,5R])}=O(|b_T|^{\kappa} )=O(1).$

\medskip

All together:
\BEQ || \ |b_T|^{1+\kappa} G^1_{a_T+\II b_T}||_{L^1([-5R,-5R])}=O(C_{6+\kappa} )
\label{eq:3.89} .\EEQ

\item[(ii)] ($\kappa$ even) Let, for $0\le m\le 3+\kappa$,
\BEQ I(x):=p.v.(\frac{1}{x})\ast \left[ x\mapsto \tilde{W}_{a_T}^{(3+\kappa-m)}(x) \partial_x^m\left(
\frac{x^3}{(x-\II b_T)^2}\right) \right].\EEQ 
We must bound $|b_T|^{1+\kappa} \int dx\,  |I(x)|.$

\noindent We rewrite $I(x)$ as the sum of two contributions (we refer to  section \ref{sec:Stieltjes} without further mention for computations and bounds related to
the principal value integral), $I(x)=:I_{reg}(x)+I_{sing}(x)$, where

\BEQ I_{reg}(x):=\int dy\, \frac{\tilde{W}_{a_T}^{(3+\kappa-m)}(x)-\tilde{W}_{a_T}^{(3+\kappa-m)}(y)}{x-y} \partial_x^m\left(
\frac{x^3}{(x-\II b_T)^2}\right)  \label{eq:Ireg} \EEQ
and
\BEQ I_{sing}(x):=\int dy\, \tilde{W}_{a_T}^{(3+\kappa-m)} (y) \frac{\partial_x^m \left(
\frac{x^3}{(x-\II b_T)^2}\right)-\partial_y^m \left(
\frac{y^3}{(y-\II b_T)^2}\right)}{x-y}. \label{eq:Ising} \EEQ

Using $\left|\frac{\tilde{W}_{a_T}^{(3+\kappa-m)}(x)-\tilde{W}_{a_T}^{(3+\kappa-m)}(y)}{x-y}
\right|\le  ||\tilde{W}_{a_T}^{(4+\kappa-m)}||_{\infty}\le  C_{7+\kappa}$, and
noting that, for $|x|\le 5R$,  the integral $\int dy \, (\cdots)=\int_{B(x,\frac{19}{2}R)} (\cdots)$  (by symmetry) simply produces an extra factor $O(1)$,  one obtains, using (i)
\BEQ |b_T|^{1+\kappa} \int dx\,  |I_{reg}(x)|=O(C_{7+\kappa} ).
\label{eq:3.93} \EEQ

Considering now $I_{sing}(x)$, we expand the numerator of (\ref{eq:Ising}) as 
in (\ref{eq:mpq}) and rewrite
\BEA &&  \frac{x^p(x-\II b_T)^{-q}-y^p(y-\II b_T)^{-q}}{x-y}  =  \frac{x^p-y^p}{x-y} (y-\II b_T)^{-q} + x^p \frac{(y-\II b_T)^q-(x-\II b_T)^q}{(x-y)(x-\II b_T)^q (y-\II b_T)^q} \nonumber\\
&&\qquad\qquad = \sum_{r=0}^{p-1} x^{p-1-r} \frac{y^r}{(y-\II b_T)^q} - \sum_{r'=1}^q \frac{x^p}{(x-\II b_T)^{r'}} 
\frac{1}{(y-\II b_T)^{q+1-r'}}. \nonumber\\ 
\EEA  

The integrals $J_1:=\int dy\,  \frac{\tilde{W}_{a_T}^{(3+\kappa-m)}(y) y^r}{(y-\II b_T)^q}$, $J_2:=\int dy\,  \frac{\tilde{W}_{a_T}^{(3+\kappa-m)}(y)}{(y-\II b_T)^{q+1-r'}}$ may be bounded using  (\ref{eq:phi-ib-n}) since $q\ge 2\ge 1$ and $q+1-r'\ge 1$. Thus
\BEQ |J_1|=O(C_{8+\kappa} ),
 |J_2|=O(C_{8+\kappa} ) \EEQ
Then $\int_{-5R}^{5R} dx\, |x|^{p-1-r}=O(1)$ and
\BEQ |b_T|^{1+\kappa} \int dx\,  \frac{|x|^p}{|x-\II b_T|^{r'}}=O(|b_T|^{(1+\kappa)-(r'-p-1)})=O(b_{max}^{2-r'+\kappa+p}) \qquad (r'\ge p+2)\EEQ
(note that $2-r'+\kappa+p\ge 0$); 
\BEQ |b_T|^{1+\kappa} \int_{-5R}^{5R} dx\, \frac{|x|^p}{|x-\II b_T|^{r'}}\le |b_T|^{1+\kappa}
\int_{-5R}^{5R} \frac{dx}{|x-\II b_T|} =O(b_{max}^{1+\kappa} \ln(1/b_{max}))
\qquad (r'=p+1)
 \EEQ
\BEQ  |b_T|^{1+\kappa} \int_{-5R}^{5R} dx\, \frac{|x|^p}{|x-\II b_T|^{r'}}\le b_{max}^{1+\kappa} \int_{-5R}^{5R}
dx\, |x|^{p-r'} = O( b_{max}^{1+\kappa} )
\qquad (r'\le p) \EEQ

\medskip

All together,  one obtains, summing all terms:
\BEQ  |b_T|^{1+\kappa} \int dx\,  |I_{sing}(x)|=O(C_{8+\kappa} ).
\label{eq:3.99}  \EEQ

\end{itemize}

\bigskip
Let us now quickly deal with the missing terms. First
\BEA || \ |b_T|^{1+\kappa} G^2_{a_T+\II b_T}||_{L^1} &\le&  b_{max}^{-2} \int dx\,
\frac{|x|^3}{|x-\II b_T|^2} |\tilde{W}_{a_T}(x)| \nonumber\\ 
&\le& b_{max}^{-2}\ ||\tilde{W}_{a_T}||_{\infty} \int_{-\frac{9}{2} R}^{\frac{9}{2}R} dx\,  |x| =
O\left( b_{max}^{-2}C_3 \right). \label{eq:3.100}  \EEA
Then one must bound $|b_T|^{1+\kappa} \int_{|x|\ge 5R} dx\,   |G^1_{a_T+\II b_T}(x)|$;
because supp$(\tilde{W}_{a_T})\subset [-\frac{9}{2}R,\frac{9}{2} R]$, this contribution vanishes
except if $\kappa$ is {\em even}, see  (ii) above, in which case (by integration
by parts)
\BEA && |b_T|^{1+\kappa} \int_{|x|\ge 5R} dx\,  |G^1_{a_T+\II b_T}(x)| = |b_T|^{1+\kappa}
\int_{|x|\ge 5R} dx\, \left| \int_{-\frac{9}{2}R}^{\frac{9}{2}R} dy \, p.v.(\frac{1}{(x-y)^2}) 
F_{a_T+\II b_T}^{(2+\kappa)}(y) \right| \nonumber\\
&& \le |b_T| \int_{|x|\ge 5R}  dx\, O(\frac{1}{x^2}) \,  ||\ |b_T|^{\kappa} F_{a_T+\II b_T}^{(2+\kappa)} ||_{L^1([-5R,5R])}.\EEA
The integral $|| \ |b_T|^{\kappa} F_{a_T+\II b_T}^{(2+\kappa)} ||_{L^1([-5R,5R])}$ is
(up to the replacement $\kappa\to \kappa-1$) exactly the one which has been computed
in case (i) above. Hence we find:
\BEQ |b_T|^{1+\kappa} \int_{|x|\ge 5R} dx\,  |G^1_{a_T+\II b_T}(x)| = 
O\left( C_{5+\kappa} \right).
\label{eq:3.102} \EEQ
Since $C_3\le C_{5+\kappa}\le C_{6+\kappa}\le 
C_{8+\kappa}=O(||V'||_{8+\kappa,[-3R,3R]})$ and $b_{max}\le 1$, the sum of  estimates
(\ref{eq:3.89}, \ref{eq:3.93}, \ref{eq:3.99}, \ref{eq:3.100}, \ref{eq:3.102}) is $
O(b_{max}^{-2}  ||V'||_{8+\kappa,[-3R,3R]})$.
 \hfill\eop


\subsection{Away from the support} 


We study in this subsection the $g$-kernels $g_{pot}^{\kappa';\kappa,ext}$ obtained
by assuming that $h=h^{ext}$, whence $|a_T|\ge 2R$. 
We want to prove that the operators
\BEQ {\cal H}^{\kappa,ext}_{pot}:  h\mapsto \left( {\cal H}^{\kappa,ext}(h):  (a,b)\mapsto\int da_T
\int_{-b_{max}}^{b_{max}} db_T\, g_{pot}^{\kappa,ext}(a,b;a_T,b_T) h(a_T,b_T)  \right)  \label{eq:Hext} \EEQ
and
\BEQ {\cal H}^{\kappa+1;\kappa,ext}_{pot}: h\mapsto  \left( {\cal H}^{\kappa+1;\kappa,ext}(h):  (a,b)\mapsto\int da_T
\int_{-b_{max}}^{b_{max}} db_T\, g_{pot}^{\kappa+1;\kappa,ext}(a,b;a_T,b_T) h(a_T,b_T)  \right)   \label{eq:Hext+1} \EEQ
are bounded operators in $L^1(\Pi_{b_{max}})$. We shall actually only prove the
 statement for ${\cal H}^{\kappa+1;\kappa,ext}_{pot}$, and leave the similar proof of the
 statement for ${\cal H}^{\kappa,ext}_{pot}$ to the reader.

\begin{Lemma} \label{lem:S3}
\begin{itemize}
\item[(i)]  $|||{\cal H}^{\kappa+1;\kappa,ext}|||_{L^1(\Pi_{b_{max}})\to
L^1(\Pi_{b_{max}})} =O\left(b_{max}^{-1} ||V'||_{4+\kappa,[-\frac{3}{2}R,
\frac{3}{2}R]} \right); $
\item[(ii)]  $|||{\cal H}^{\kappa,ext}|||_{L^1(\Pi_{b_{max}})\to
L^1(\Pi_{b_{max}})} =O\left( ||V'||_{3+\kappa,[-\frac{3}{2}R,
\frac{3}{2}R]}  \right).$
\end{itemize}

\medskip\noindent The same estimates hold for the $L^{\infty}$-operator norms $|||\ \cdot\ |||_{L^{\infty}(\Pi_{b_{max}})\to L^{\infty}(\Pi_{b_{max}})}$.
\end{Lemma}

\noindent As can be checked by
looking at the details of the proof, norms
$|||\ \cdot\ |||_{L^1(\Pi_{b_{max}}) \to L^{\infty}(\Pi_{b_{max}})}$ are deduced
from these by dividing by a volume factor $\Vol([-3R,3R]\times[-b_{max},b_{max}])\approx b_{max}$, which follows from the fact that the kernel $g_{pot}^{\kappa';\kappa,ext}$ is regular on the diagonal.

\medskip\noindent{\bf Proof.}

\begin{itemize}
\item[(1)] (operator norm of ${\cal H}_{pot}^{\kappa+1;\kappa,ext}$)
 By definition,

\BEQ  g^{\kappa+1;\kappa,ext}_{pot}(a,b;a_T,b_T)= {\bf 1}_{|b|<b_{max}} {\bf 1}_{|a_T|\ge 2R}\,  \frac{1}{(1+\kappa)!}  (-\II b) |b_T|^{\kappa}\,  {\cal K}^{\kappa+1}_{b_{max}}(\chi_R  V' f'_{z_T})(a)  \label{eq:3.43}
 \EEQ
Note that, since $|a_T|\ge 2R$ and supp$(\chi_R)\subset B(0,\frac{3}{2}R)$, the function
$\chi_R V' f'_{z_T}$ is regular and bounded. Furthermore, for every $n=0,1,\ldots,$  $||(\chi_R V' f'_{z_T})^{(n)}||_{\infty}=O(||V'||_{n,[-\frac{3}{2}R,\frac{3}{2}R]}  |a_T|^{-2})$ is integrable at infinity in $a_T$.

\medskip\noindent
We consider in the computations only the part $g^{\kappa}$ of the kernel $g^{\kappa+1;\kappa,ext}_{pot}$ obtained by replacing ${\cal K}^{\kappa}$ with the operators
in factor of the bounded operator $1+\underline{\cal K}^{\kappa}_{b_{max}}$ in
 (\ref{eq:calK-kappa-even}, \ref{eq:calK-kappa-odd}),  and distinguish the cases $\kappa$ odd,  $\kappa$ even. 
Assume first $\kappa$ is {\em odd}.
Then 

\BEA &&  |g^{\kappa}(a,b;a_T,b_T)|  \nonumber\\
&& = {\bf 1}_{|b|<b_{max}}\, {\bf 1}_{|a_T|\ge 
2R} \,  O(|b_T|^{1+\kappa}) \left(  \left|   (\chi_R V' f'_{z_T})^{(3+\kappa)}(a) \right| + b_{max}^{-(3+\kappa)}
|(\chi_R V'f'_{z_T})(a)| \right) \nonumber\\
&& = O\left(  b_{max}^{1+\kappa} ||V'||_{3+\kappa,[-\frac{3}{2}R, \frac{3}{2}R]}
+ b_{max}^{-2}  ||V'||_{0,[-\frac{3}{2}R, \frac{3}{2}R]} \right)\, |a_T|^{-2} \nonumber\\
&&\qquad\qquad\qquad \, {\bf 1}_{(a,b)\in[- \frac{3}{2} R,\frac{3}{2} R]\times [-b_{max},b_{max}]}.
\nonumber\\  \label{eq:3.45}
\EEA

Hence ${\cal H}^{\kappa+1;\kappa,ext}_{pot}$ is
a bounded operator on $L^1(\Pi_{b_{max}})$, with $L^1$-operator norm
\BEA && |||{\cal H}^{\kappa+1;\kappa,ext}_{pot} |||_{L^1\to L^1}  = O\left( b_{max}^{1+\kappa} ||V'||_{3+\kappa,[-\frac{3}{2}R, \frac{3}{2}R]} + |b|^{-2}_{max} ||V'||_{0,[-\frac{3}{2}R, \frac{3}{2}R]} \right)  \ \cdot \nonumber\\
&&\qquad\qquad\qquad \cdot\ \Vol([-3R,3R]\times
[-b_{max},b_{max}])  \nonumber\\
&&\qquad =O\left( b_{max}^{3+\kappa} ||V'||_{3+\kappa,[-\frac{3}{2}R,\frac{3}{2}R]} + 
||V'||_{0,[-\frac{3}{2}R, \frac{3}{2}R]} \right) b_{max}^{-1}. 
\nonumber\\
\label{eq:op-ext1}
\EEA

\medskip\noindent Assume now that $\kappa$ is {\em even}. Then (using (\ref{eq:3.45}) and (\ref{eq:bound-vp2}))\\  $|g^{\kappa}(a,b;a_T,b_T)|\le O(b_{max}^{-2}  ||V'||_{0,[-\frac{3}{2}R, \frac{3}{2}R]}  \, {\bf 1}_{(a,b)\in[-\frac{3}{2}R,\frac{3}{2}R]\times [-b_{max},b_{max}]}) |a_T|^{-2}$, plus
\BEA  &&  {\bf 1}_{|b|<b_{max}} O(|b_T|^{1+\kappa}) \left|  {\bf 1}_{|a_T|\ge 
2R} \,  p.v.(\frac{1}{x})\ast (\chi_R V' f'_{z_T})^{(3+\kappa)}(a) \right|
\nonumber\\
&&=  {\bf 1}_{|b|<b_{max}} O(|b_T|^{1+\kappa}) \left|  {\bf 1}_{|a_T|\ge 
2R} \,  p.v.(\frac{1}{x^2})\ast (\chi_R V' f'_{z_T})^{(2+\kappa)}(a) \right|
\nonumber\\
&&= O(|b_T|^{1+\kappa})\  \left[ {\bf 1}_{(a,b)\in[-3R,3R]\times [-b_{max},b_{max}]} \  || {\bf 1}_{|a_T|\ge 2R} (\chi_R V' f'_{z_T})^{(4+\kappa)}||_{\infty}  \right.\nonumber\\
&& \qquad \qquad \left. + {\bf 1}_{(a,b)\in (\R\setminus[-3R,3R])\times[-b_{max},b_{max}]} 
\ \frac{1}{a^2} \ 
|| {\bf 1}_{|a_T|\ge 2R} (\chi_R V' f'_{z_T})^{(2+\kappa)}||_{\infty} \right].
\nonumber\\ 
\EEA
We conclude  by multiplying by $\bar{\chi}_R+(1-\bar{\chi}_R)$ that 
${\cal H}^{\kappa+1;\kappa,ext}_{pot} h=\tilde{h}^{int}+\tilde{h}^{ext}$, where
\BEQ ||\tilde{h}^{int}||_{L^1}\le O\left(  b_{max}^{3+\kappa} ||V'||_{4+\kappa,[-\frac{3}{2}R, \frac{3}{2}R]}
+
||V'||_{0,[-\frac{3}{2}R, \frac{3}{2}R]}   \right) \  b_{max}^{-1} \ ||h||_{L^1}  \label{eq:op-ext2} \EEQ
and
\BEQ ||\tilde{h}^{ext}||_{L^1}\le O( b_{max}^{2+\kappa} ||V'||_{2+\kappa,[-\frac{3}{2}R, \frac{3}{2}R]})\   ||h||_{L^1} .  \label{eq:op-ext3}
\EEQ

\item[(2)](operator norm of ${\cal H}_{pot}^{\kappa,ext}$)
With respect to (\ref{eq:op-ext1}), (\ref{eq:op-ext2}) and (\ref{eq:op-ext3}),
we remove one order of differentiation and one power of $b_{max}^{-1}$, and
exchange parities. Thus
\BEQ |||{\cal H}^{\kappa,ext}_{pot} |||_{L^1\to L^1}= O\left( b_{max}^{2+\kappa} ||V'||_{2+\kappa, [-\frac{3}{2}R,\frac{3}{2}R]} + 
||V'||_{0,[-\frac{3}{2}R, \frac{3}{2}R]} \right)  \label{eq:op-ext4}
\EEQ
for $\kappa$ even, while for $\kappa$ odd, ${\cal H}^{\kappa,ext}_{pot} h=
\tilde{h}^{int}+\tilde{h}^{ext}$, with
\BEQ ||\tilde{h}^{int}||_{L^1}\le O\left(b_{max}^{2+\kappa} ||V'||_{3+\kappa,[-\frac{3}{2}R, \frac{3}{2}R]}
+
||V'||_{0,[-\frac{3}{2}R, \frac{3}{2}R]}  \right)   \ ||h||_{L^1}  \label{eq:op-ext5} \EEQ
and
\BEQ ||\tilde{h}^{ext}||_{L^1}\le O( b_{max}^{2+\kappa} ||V'||_{1+\kappa,[-\frac{3}{2}R, \frac{3}{2}R]})\    ||h||_{L^1} .  \label{eq:op-ext6}
\EEQ

\end{itemize}

\bigskip\noindent
{\bf Remark.}  When $b_{max}=+\infty$, 
 the adjoint of ${\cal H}^{\kappa,ext}_{pot}$ or ${\cal H}^{\kappa+1;\kappa,ext}_{pot}$ is the generator of a
{\em signed jump Markov process} with good properties. Namely, for all
$a_T,b_T$,
\BEQ \int_{-\infty}^{+\infty} da \int_{-\infty}^{+\infty} db\, g_{pot}^{\kappa,ext}(a,b;a_T,b_T)=0 \EEQ
since $K^{\kappa}_{+\infty}(s=0)=0$, hence ${\cal L}^{\kappa,ext}_{pot}:=
({\cal H}^{\kappa,ext}_{pot})^{\dagger}$ may we written in the following form,
\BEQ {\cal L}^{\kappa,ext}_{pot}(\phi)(a_T,b_T)=\int_{-\infty}^{+\infty} da\,  \int_{-\infty}^{+\infty} db\,  (\phi(a,b)-\phi(a_T,b_T))g^{\kappa,ext}_{pot}(a,b;a_T,b_T). \label{eq:L-ext} \EEQ

For a bona fide Markov process, the function $(a,b)\mapsto -g^{\kappa,ext}_{pot}(a,b;a_T,b_T)\ge 0$ would be the jump rate from to
$(a_T,b_T)$ to $(a,b)$, and one would have $-\int da\, db\, g^{\kappa,ext}_{pot}(a,b;a_T,b_T)=1$. Here
$g$ is a signed kernel, so the probabilistic interpretation fails stricto sensu. However,  
the $L^1$ semi-group generated by ${\cal H}^{\kappa,ext}_{pot}$  has good properties
because ${\cal H}^{\kappa,ext}_{pot}$ is a bounded operator (see section 4).



\subsection{Boundary terms}


Recall the {\em horizontal boundary terms} h-bdry$_0$ (\ref{eq:h-bdry0}),
h-bdry$_{pot}^{(1)}$ (\ref{eq:h-bdry-1}), and the {\em vertical boundary terms}, 
v-bdry$_0$ (\ref{eq:v-bdry0}), v-bdry$_{pot}^{(0)}$ (\ref{eq:v-bdry-0}) and 
v-bdry$_{pot}^{(2)}$ (\ref{eq:v-bdry-2}). We adopt the following notations.
First let 
\BEQ \partial_R \Pi_{b_{max}}:=\partial\left([-3R,3R]\times [-b_{max},b_{max}]\right)= \left([-3R,3R]\times\{\pm b_{max}\}\right) \cup \left( \{\pm 3R\}\times
[-b_{max},b_{max}]\right),\EEQ
and, for $h:[-3R,3R]\times[-b_{max},b_{max}]\to\C$,  let $\partial_R h:=
h\big|_{\partial_R \Pi_{b_{max}}}$ be the restriction of $h$ to the boundary. 
Horizontal boundary terms h-bdry$_{\cdot}^{\cdot}$ (where the dots stand for 
various lower and upper indices) are of the form
h-bdry$_{\cdot}^{\cdot}(b_{max};\partial_R h)$ \ $-$\  h-bdry$_{\cdot}^{\cdot}(-b_{max};\partial_R h)$;
similarly, vertical boundary terms  v-bdry$_{\cdot}^{\cdot}$ are of the form
v-bdry$_{\cdot}^{\cdot}(3R;\partial_Rh)$ \ $-$ \  v-bdry$_{\cdot}^{\cdot}(-3R;\partial_R h)$.

\medskip\noindent  These terms depend
 a priori on the restriction of   $h$ to
 $\partial_R \Pi_{b_{max}}$, which is incompatible with the $L^1$-norms. In order to avoid that, we replace ${\mathfrak{f}}_{a_T\pm \II b_{max}}$
 with ${\cal C}^{\kappa'} \left((a,b)\mapsto  {\cal K}^{\kappa'}_{b_{max}}({\mathfrak{f}}_{a_t\pm\II b_{max}})(a)\right)$, $\kappa'=\kappa,\kappa+1$ so that

\BEA &&
({\mathrm{h\!-\!bdry}}_0+{\mathrm{h\!-\!bdry}}_{pot}^{(1)})(b_{max};\partial_R h) - (b_{max}\leftrightarrow -b_{max}) \nonumber\\
&& \qquad ={\cal C}^{\kappa+1}
({\cal H}^{\kappa+1;\kappa,h-bdry}(\partial_R h))  ={\cal C}^{\kappa}({\cal H}^{\kappa,h-bdry}(\partial_R h)) 
\EEA
 
with ${\cal H}^{\kappa+1;\kappa,h-bdry},{\cal H}^{\kappa,h-bdry}: L^{\infty}(\partial_R \Pi_{b_{max}})\to (L^1 \cap L^{\infty})(\Pi_{b_{max}})$,

\BEQ {\cal H}^{\kappa+1;\kappa}_{\mathrm{h-bdry}}(\partial_R h)(a,b)=\int_{-3R}^{3R} da_T
\, g^{\kappa+1;\kappa,}_{\mathrm{h-bdry}}(a,b;a_T) \partial_R h(a_T,b_{max}) - 
(b_{max}\leftrightarrow -b_{max})  \label{eq:Hhbdry+1}
\EEQ

\BEQ {\cal H}^{\kappa}_{\mathrm{h-bdry}}(\partial_R h)(a,b)=\int_{-3R}^{3R} da_T
\, g^{\kappa}_{\mathrm{h-bdry}}(a,b;a_T) \partial_R h(a_T,b_{max}) - (b_{max}\leftrightarrow -b_{max})  \label{eq:Hhbdry}
\EEQ

and
 
\BEA &&  g^{\kappa';\kappa}_{\mathrm{h-bdry}}(a,b;a_T)= -  \left\{  \frac{\beta}{4} \frac{b_{max}^{1+\kappa}}{(1+\kappa)!} \Im[(M^N_T+M_T)(a_T+\II b_{max})]\, \right.\nonumber\\
&& \left. \qquad + b_{max}^{2+\kappa} V''(a_T)
  \right\}  \  {\cal K}^{\kappa'}_{b_{max}}({\mathfrak{f}}_{a_T+\II b_{max}})(a)
   \label{eq:g-bdry} \EEA
   
Similarly, we replace $f_{\pm 3R+\II b}$ with ${\cal C}^{\kappa'} \left((a,b)\mapsto  {\cal K}^{\kappa'}_{b_{max}}({\mathfrak{f}}_{\pm 3R+\II b})(a)\right)$, $\kappa'=\kappa,\kappa+1$ , so that
\BEA &&  ({\mathrm{v\!-\!bdry}}_0 + {\mathrm{v\!-\!bdry}}_{pot}^{(0)}+ {\mathrm{
v\!-\!bdry}}_{pot}^{(2)})(3R;\partial_R h) - (R\leftrightarrow -R) \nonumber\\
&& \qquad = {\cal C}^{\kappa+1}({\cal H}^{\kappa+1;\kappa}_{v-bdry}(\partial_R h)) = {\cal C}^{\kappa}({\cal H}^{\kappa}_{v-bdry}(\partial_R h)) 
\EEA

with ${\cal H}^{\kappa+1;\kappa}_{\mathrm{v-bdry}},{\cal H}^{\kappa}_{\mathrm{v-bdry}}: L^{\infty}(\partial_R \Pi_{b_{max}})\to (L^1 \cap L^{\infty})(\Pi_{b_{max}})$,

\BEQ {\cal H}^{\kappa+1;\kappa}_{\mathrm{v-bdry}}(\partial_R h)(a,b)=\int_{-b_{max}}^{b_{max}} db_T
\, g^{\kappa+1;\kappa}_{\mathrm{v-bdry}}(a,b;b_T) \partial_R h(3R,b_T) - 
(R\leftrightarrow -R)  \label{eq:Hvbdry+1}
\EEQ

\BEQ {\cal H}^{\kappa}_{\mathrm{v-bdry}}(\partial_R h)(a,b)=\int_{-b_{max}}^{b_{max}} db_T
\, g^{\kappa}_{\mathrm{v-bdry}}(a,b;b_T) \partial_R h(3R,b_T)  - 
(R\leftrightarrow -R)  \label{eq:Hvbdry}
\EEQ

and
 
\BEA &&  g^{\kappa';\kappa}_{\mathrm{v-bdry}}(a,b;b_T)= -(-\II b_T)  \left\{  \frac{\beta}{4}  \frac{|b_T|^{\kappa}}{(1+\kappa)!}   \Re[(M^N_T+M_T)(3R+\II b_T)]\,  \right.\nonumber\\
&& \left. \qquad + (|b_T|^{\kappa} V'(3R)+ |b_T|^{\kappa+2} V'''(3R)
   \right\} \     {\cal K}^{\kappa'}_{b_{max}}(\chi_{|R|}\,  {\mathfrak{f}}_{3R+\II b_T})(a)
    \EEA

\medskip

\begin{Lemma} \label{lem:S9}
Let $\kappa\ge 0$ and $\kappa'=\kappa,\kappa+1$. Then
\begin{itemize}
\item[(i)] $|||{\cal H}^{\kappa';\kappa}_{\mathrm{h-bdry}}|||_{L^{\infty}(\partial_R\Pi_{b_{max}})\to L^{\infty}(\Pi_{b_{max}})}=b_{max}^{-3-(\kappa'-\kappa)} \left\{ O(\frac{1}{b_{max}})+O(b_{max} ||V''||_{0,[-3R,3R]} \right\}$ and\\
 $|||{\cal H}^{\kappa';\kappa}_{\mathrm{h-bdry}}|||_{L^{\infty}(\partial_R\Pi_{b_{max}})\to L^1(\Pi_{b_{max}})}= b_{max}^{-2-(\kappa'-\kappa)} \left\{ O(\frac{1}{b_{max}})+O(b_{max} ||V''||_{0,[-3R,3R]} \right\}$;

\item[(ii)] 
\BEA && |||{\cal H}^{\kappa';\kappa}_{\mathrm{v-bdry}}|||_{L^{\infty}(\partial_R\Pi_{b_{max}})\to L^{\infty}(\Pi_{b_{max}})} \nonumber\\
&&\qquad 
= b_{max}^{-1-(\kappa'-\kappa)}  \left\{ O(\frac{1}{b_{max}})+O(||V'||_{0,[-3R,3R]}) + O(b_{max}^2 ||V'''||_{0,[-3R,3R]}) \right\} \nonumber\\
\EEA
 and
\BEA && |||{\cal H}^{\kappa}_{\mathrm{v-bdry}}|||_{L^{\infty}(\partial_R\Pi_{b_{max}})\to L^1(\Pi_{b_{max}})} \nonumber\\
&&\qquad =b_{max}^{-(\kappa'-\kappa)} \left\{ O(\frac{1}{b_{max}})+O(||V'||_{0,[-3R,3R]}) + O(b_{max}^2 ||V'''||_{0,[-3R,3R]}) \right\}. \nonumber\\
\EEA
\end{itemize}
\end{Lemma}

\noindent As in the two previous subsections,  $L^1$-estimates and $L^{\infty}$-estimates differ
only by a volume factor $\approx b_{max}$.

\medskip\noindent {\bf Proof.} 
\begin{itemize}
\item[(i)]
 Immediate consequence of the bounds
$\left|\Im[(M^N_T+M_T)(a_T+\II b_{max})] \right|\le 2/b_{max}$,
$b_{max} |V''(a_T)|\le b_{max} ||V''||_{0,[-3R,3R]}$  and   
$| {\cal K}^{\kappa'}_{b_{max}}(\chi_R {\mathfrak{f}}_{a_T\pm\II b_{max}})(a)|=O(b_{max}^{-(4+\kappa')}) \ O(\frac{1}{(1+|a|)^2})$ (as seen by using (\ref{eq:calK-kappa-even}) when
$\kappa'$ is even, and (\ref{eq:calK-kappa-odd},\ref{eq:bound-vp2}) when
$\kappa'$ is odd). 

\item[(ii)]  Immediate consequence of the bounds
$\left|\Re[(M^N_T+M_T)(a_T+\II b_{max})] \right|\le 2/b_{max}$,
$|V'(3R)|\le ||V'||_{0,[-3R,3R]}$, $b_T^2 |V'''(3R)|\le b_{max}^2 ||V'''||_{0,[-3R,3R]}$,   and\\   
$| {\cal K}^{\kappa'}_{b_{max}}(\chi_R \, {\mathfrak{f}}_{\pm 3R+\II b_T})(a)|=O(b_{max}^{-(2+\kappa')}) \ O(\frac{1}{(1+|a|)^2})$.

\end{itemize} \hfill \eop


\section{Gaussianity of the fluctuation process}


In this section, we prove our Main Theorem (see \S 1.3), namely, we prove that the finite-$N$ fluctuation process $(Y^N_t)_{t\ge 0}$
converges weakly in $C([0,T],H_{-14})$  to a fluctuation process $(Y_t)_{t\ge 0}$, which is
the unique solution of a martingale problem that we solve explicitly in terms of
the solution of (\ref{eq:PDE-f}). 

\noindent {\em We fix once and for all: $b_{max}=\half$}.

\medskip\noindent Summarizing what we have found up to now and applying  Proposition \ref{prop:Israelsson}, we find for $\kappa=0,1,2,\ldots$:
\BEQ d\langle Y_t^N,{\cal C}^{\kappa} h_t\rangle=\half(1-\frac{\beta}{2}) \langle X_t^N, ({\cal C}^{\kappa} h_t)''\rangle \, dt + \frac{1}{\sqrt{N}} \sum_i  ({\cal C}^{\kappa} h_t)'(\lambda^i_t) dW^i_t \EEQ

where $(h_t)_{t\le T}$ is the solution of the evolution equation
\BEQ \frac{dh_t}{dt}={\cal H}^{\kappa}(t)h(t).\EEQ

On the other hand, the process $(Y^N_t)_{t\ge 0}$ is a solution of the following 
{\em martingale problem}: if $\bar{\phi}=\{\phi_j\}_{1\le j\le k}$ is a family of test
functions in $C_c^{\infty}(\R,\R)$, $F\in C_b^2(\R^k,\R)$, then, letting $F_{\bar{\phi}}(Y_{\cdot}^N):=F(\langle Y_{\cdot}^N,
\phi_1\rangle,\ldots,\langle Y_{\cdot}^N,\phi_k\rangle)$,
\BEQ \Phi^{T,N}_t(Y^N):=F_{\bar{\phi}}(Y^N_T)-F_{\bar{\phi}}(Y^N_t)-\int_t^T ds\, L^N_s F_{\bar{\phi}}(Y^N_s)\EEQ
is a martingale, 
where
\BEA && L^N_s F_{\bar{\phi}}(Y^N_s):=\sum_{j=1}^k \frac{\partial F_{\bar{\phi}}}{\partial
x_j}(Y^N_s) \left( \langle Y^N_s, \frac{\beta}{4} \int \frac{\phi'_j(\cdot)-\phi'_j(y)}{\cdot-y}
(X_t+X_t^N)(dy) - V'(\cdot)\phi'_j(\cdot) \rangle + \right. \nonumber\\
&& \qquad\qquad\left. + \half(1-\frac{\beta}{2}) \langle
X_t^N,\phi''_j\rangle \right) +\half  \sum_{j,l=1}^k  \partial^2_{jl} F_{\bar{\phi}}(Y_s^N)\  \langle X_t^N,\phi'_j\phi'_l\rangle \label{eq:mgle-pb}
\EEA
(see \cite{Isr}, p. 28--29).

\bigskip\noindent
We now use an exponential functional of the process to derive the limit law.

\begin{Definition}
For $\kappa=0,1,2,\ldots$ and $h\in (L^1\cap L^{\infty})(\Pi_{b_{max}})$, let
\BEQ \phi_h(Y_t^N):=e^{\II \langle Y_t^N,{\cal C}^{\kappa} h\rangle}. \EEQ
\end{Definition}

Notice the slight change of notation with respect to (\ref{eq:intro-Gaussian}) in the introduction. 
It\^o's formula implies as in (\cite{Isr}, p. 29)
\BEQ d\phi_{h_t}(Y_t^N)=\phi_{h_t}(Y_t^N) \left(\half \II (1-\frac{\beta}{2}) \langle
X_t^N,({\cal C}^{\kappa} h_t)''\rangle -\half \langle X_t^N, (({\cal C}^{\kappa}h_t)')^2 \rangle \right) dt,\EEQ
where $({\cal C}^{\kappa}h_t)_{0\le t\le T}$ is the solution of (\ref{eq:intro-PDE-f}) for $N$ finite, 
from which for $0\le t\le T$ (letting {\em formally} $N\to\infty$)
\BEQ \esper[\phi_{h_T}(Y_T)\big| {\cal F}_t]=\phi_{h_t}(Y_t)\  \exp\left( \half
\int_t^T  \left[ \II (1-\frac{\beta}{2}) \langle X_s,({\cal C}^{\kappa} h_s)''\rangle -
\langle X_s,(({\cal C}^{\kappa} h_s)')^2\rangle \right] ds  \right)  \label{eq:Gaussian} \EEQ
where  $(h_t)_{0\le t\le T}$ is now  the solution of  the asymptotic
equation (\ref{eq:intro-PDE-f-asymptotic}).
Since $h_s$, $t\le s\le T$ are linear in $h_T$, the term in the exponential in (\ref{eq:Gaussian})
is a sum of a linear and a quadratic term in $h_T$, giving resp. the expectation 
and the variance of a Gaussian process (see Israelsson, \S 2.6 for more details).

\bigskip

The strategy of the proof, following closely the proof in (Israelsson \cite{Isr}, section 2), is the following:

\begin{itemize}
\item[(A)] {\em find bounds} for $\esper[ \sup_{s\le T} |\langle Y_s^N,\phi\rangle|]$, 
$\esper[\sup_{s\le T} |\langle Y_s^N,\int \frac{\phi'(\cdot)-\phi'(y)}{\cdot-y} X_s^N(dy)\rangle|]$
and\\  $\esper[ \sup_{s\le T} |\langle Y_s^N,V'(\cdot)\phi'(\cdot)\rangle|]$
(see first line in (\ref{eq:dY_tnf}) or (\ref{eq:mgle-pb}));

\item[(B)] prove a {\em tightness property} for the family of processes $Y^N$,
implying the existence of a (non necessarily unique) limit in law;

\item[(C)] prove that {\em any} weak limit $Y$ of the $(Y^N)_{n\ge 1}$ {\em satisfies the
limit martingale problem}, i.e. is such that
\BEQ \Phi^{T}_t(Y)=F_{\bar{\phi}}(Y_T)-F_{\bar{\phi}}(Y_t)-\int_t^T ds\, L_s F_{\bar{\phi}}(Y_s)\EEQ
is a martingale, 
where
\BEA && L_s F_{\bar{\phi}}(Y_s):=\sum_{j=1}^k \frac{\partial F_{\bar{\phi}}}{\partial
x_j}(Y_s) \left( \langle Y_s, \frac{\beta}{2} \int \frac{\phi'_j(\cdot)-\phi'_j(y)}{\cdot-y}
X_t(dy) - V'(\cdot)\phi'_j(\cdot) \rangle + \right. \nonumber\\
&& \qquad\qquad\left. + \half(1-\frac{\beta}{2}) \langle
X_t,\phi''_j\rangle \right) +\half  \sum_{j,l=1}^k  \partial^2_{jl} F_{\bar{\phi}}(Y_s) \  \langle X_t,\phi'_j\phi'_l\rangle \label{eq:lim-mgle-pb}
\EEA
(obtained formally from (\ref{eq:mgle-pb}) by letting $N\to\infty$);

\item[(D)] prove that there exists only {\em one} measure with given 
initial measure satisfying (\ref{eq:lim-mgle-pb}), and that it is Gaussian, and
satisfies (\ref{eq:Gaussian}).
\end{itemize}

To prove (C) one must essentially prove that $(\Phi^T_t-\Phi^{T,N}_t)(Y^N)$ is  small,
since $\Phi_t^T$ is continuous (see \cite{Isr}, pp. 47-48).

\bigskip
The proof is borrowed from Israelsson, where it takes up a few pages, see \cite{Isr}.
The main bound, compare with \cite{Isr}, Lemma 15, is the following "$H_8$"-bound,
\BEQ  \esper[ \sup_{s\le T} |\langle Y_s^N,\phi\rangle|]\le C_T ||\phi||_{H_8},
\label{eq:bd-H4} \EEQ

The above $H_8$-bound is obtained by integrating
(\ref{eq:dY_tnf}), namely,
\BEA && \esper\Big[ \sup_{s\le T} |\langle Y_s^N,\phi\rangle| \Big]\le C\Big( \esper[ |\langle Y_0^N,\phi\rangle|]+  \nonumber\\
&& \qquad +  \int_0^T ds \, \esper\Big[|\langle Y_s^N, V'(\cdot) \phi'(\cdot)\rangle|+
|\langle Y_s^N,
\int \frac{\phi'(x)-\phi'(\cdot)}{x-\cdot} (X_s^N(dx)+X_s(dx))\rangle| \Big]  \nonumber\\
&&  \qquad  + \int_0^T  ds\, \esper\Big[|\int \phi''(x)X_s^N(dx)| \Big] + \esper\Big[\sup_{s\le T} |M_s|\Big]  \Big)
\EEA
where $M_s:=\frac{1}{\sqrt{N}} \int_0^s \sum_{i=1}^N \phi'(\lambda_s^i)dW_s^i$ is
a martingale. All these terms are bounded as in \cite{Isr}, Lemma 15, using 
easy, $V$-independent inequalities {\em and} the following {\em fundamental estimate},

\begin{Lemma} \label{lem:fundamental-estimate}
There exists a constant $C$ depending on $T$ such that, for all $t\le T$, 
\BEQ \esper[|\langle Y_t^N,\phi\rangle|^2] \le C||\phi||^2_{H_7},\label{eq:fundamental-estimate}
\EEQ
\end{Lemma}
itself an immediate consequence of
\BEQ \sqrt{\esper[|\langle Y_T^N,{\mathfrak{f}}_z\rangle|^2]} \equiv 
\sqrt{ \esper[|N(M_t^N-M_t)(z)|^2]}  \le C |b|^{-6}, \qquad
b:=\Im z \in [-b_{max},b_{max}]   \label{eq:pre-4.14} \EEQ
see (\ref{eq:4.14}) below.

\medskip\noindent Note the loss of regularity with respect to (\cite{Isr}, Lemma 14), where one has
$||\phi||_{H_2}$ in the r.-h.s.  This changes the bounds in
the course of the proof of (\cite{Isr}, Lemma 15), see in particular p. 43, l. 5,  where estimates are proved
using  Israelsson's Lemma 8 for $q=2$ (with $q-1$ playing the same r\^ole as our $\kappa$): the latter lemma yields (with our notations) a bound on $||{\cal K}^{\kappa}_{b_{max}}(f)||_{L^1}$ in terms of $||f||_{L^1}+||f||_{H_{\kappa+2}}$, resp.
$||f||_{L^1}+||f||_{H_{\kappa+3}}$, depending whether $\kappa$ is even, resp. odd;
because of the loss of regularity one must take $\kappa=5$ ($q=6$, same exponent
as in (\ref{eq:pre-4.14})), hence the "$H_8$"=$H_{\kappa+3}$-bound.

\medskip\noindent The proof of (C), see \cite{Isr}, Lemmas 17 and 20, mainly depends on a
bound for\\ $C(\phi):=\esper\Big[ \sup_{0\le s\le T} \Big| \Big\langle Y_s^N, \int
\frac{\phi'(\cdot)-\phi'(y)}{\cdot - y} X_s^N(dy) \Big\rangle \Big| \Big]$. Assume
$\phi={\cal C}^{\kappa}h$, $\kappa\in\N$; then  $C(\phi)$ is bounded in terms of
the integral against the measure $|b|^{\kappa+1} |h(a,b)| \, da\, db$ on $\Pi_{b_{max}}$ of the random function $N(M_t^N(z)-M_t(z))^2$, averaging to  $O(\frac{1}{N} |b|^{-12})$ by (\ref{eq:pre-4.14}).  Therefore the integral converges if $\kappa\ge 11$, and is
bounded, as recalled in the previous paragraph, by $O(1/N) ||\phi||_{H_{14}}$.

 \medskip\noindent Then the tightness property (B) is proved using a lemma due to Mitoma \cite{Mit}
and the above estimates (see \cite{Isr}, \S 2.4); the Sobolev space
$H_{-14}$ in our Main Theorem is such that there
exists a nuclear mapping $H_{14}\to H_{8}$ (see \ref{eq:bd-H4}), as follows from
Treves \cite{Tre}, which is a requirement in Mitoma's hypotheses. 

\medskip\noindent Finally, letting
$\phi_i={\cal C}^0 h_i$, $i=1,\ldots,k$,  (D) is proved by 
computing\\  $\esper\left[ \exp\II  \left(\langle Y_{t_1},\phi_1\rangle+\ldots+
\langle Y_{t_k},\phi_k\rangle\right)\right]$, $0\le t_k\le \ldots\le t_1\le T$
by induction on $k$ using the assumed martingale property of the limit(s) and solving
in terms of the time-evolved functions $h_i(t)$, $i=1,\ldots,k$. For $k=2$ we
obtain (\ref{eq:Gaussian}).

\bigskip

\noindent So everything boils down to the proof of the above lemma.

\medskip
\noindent {\bf Proof  of Lemma \ref{lem:fundamental-estimate}.}  Let $\phi\in C_c^{\infty}$. Consider
its standard Stieltjes decomposition of order 5, $\phi={\cal C}^5 h$ (take
 $b_{max}=\half$). Then (using  (\ref{eq:Stieltjes-decomposition}) and
 Cauchy-Schwarz's inequaity)
\BEQ |\langle Y_t^N,\phi\rangle|^2 \le ||h||^2_{L^2(\Pi_{b_{max}})} \int_{\Pi_{b_{max}}} da\, db\, b^{12} (N |M_t(z)-M_t^N(z)|)^2 \label{eq:bound-YtNphi2}
\EEQ
and $||h||_{L^2(\Pi_{b_{max}})}=O(||\phi||_{H^7})$ (as follows from the
kernel representation (\ref{eq:calK-kappa-even}) of the standard Stieltjes 
decomposition, together with Parseval-Bessel's formula).
Hence (\ref{eq:fundamental-estimate}) follows if we can show that
\BEQ \esper[|N(M_T^N(z)-M_T(z))|^2]=\esper[|\langle Y_T^N,{\mathfrak{f}}_z\rangle|^2] \le C b^{-12} \label{eq:4.14} \EEQ
for $0<b:=\Im z\le \half$; compare with Israelsson \cite{Isr}, Proposition 1, where a much
better bound in $O((\ln(1+1/b)/b)^2)$ is proved. Note, however, that there is, to the
best of our understanding,
a flaw in Israelsson's proof, see below (\ref{eq:Ii}), whence (despite some efforts) we find
 in fact a bound in $O(b^{-12})$ in the harmonic case too.
Further, introducing the stopping time
\BEQ \tau:=\inf\{0<t\le T: \sup_{1\le i\le N} |\lambda_t^i|>R(T)\},\EEQ
(and letting by convention $\tau\equiv T$ if $\sup_{0\le t\le T} \sup_{1\le i\le N}
|\lambda^i_t|\le R(T)$), see Lemma \ref{lem:support}, we have, using the large deviation bound of  section 5,
and the obvious deterministic bound $M_t^N(z)\le |\Im(z)|^{-1}$,
\BEQ \esper[|N(M_T^N(z)-M_T(z))|^2] \le \esper[|N(M_{\tau}^N(z)-M_{\tau}(z))|^2] + Ce^{-cN} N^2  |\Im z|^{-2}.\EEQ 
So (\ref{eq:4.14}) holds provided we show that
\BEQ  \esper[|N(M_{\tau}^N(z)-M_{\tau}(z))|^2] \le C b^{-12} \label{eq:4.14bis} \EEQ
where now by construction $\sup_i |\lambda_t^{N,i}|\le R(T)$ for all $t\le \tau$, a
support condition which is  essential for the subsequent computations.

\bigskip\noindent Before we can do that, however, we need a long preliminary discussion. 
Indeed, Israelsson's proof of this fact in his Proposition 1 does not carry through immediately to
the case of a general $V$, because it relies in an essential way on
the bounds on characteristics. As explained in the Introduction though, the deterministic
characteristics due to the $\left(\frac{1}{x}\right)$-potential, see \S 3.1, to which
we can safely add the other transport generators without much change, yield the
most singular contribution, so our strategy is to treat the
non-local term  ${\cal H}_{nonlocal}$ as a {\em perturbation} of ${\cal H}_{transport}$, by using a Green function expansion. Note however the {\em twist} here: the operators ${\cal H}_{transport}^{\kappa}$ and ${\cal H}_{nonlocal}^{\kappa}$ are endomorphisms of $L^1(\Pi_{b_{max}})$, but ${\cal H}_{nonlocal}^{\kappa+1;\kappa}$ intertwines in some sense {\em two different copies} of  $L^1(\Pi_{b_{max}})$,
with different mappings, ${\cal C}^{\kappa}$, vs. ${\cal C}^{\kappa+1}$ to
$L^1(\R)$. The intertwining is not trivial, in the sense that $b{\cal H}_{nonlocal}^{\kappa+1;\kappa}\not={\cal H}_{nonlocal}^{\kappa}$. This leads us to introduce the following
operator-valued matrices.

\begin{Definition} \label{def:eps}
\begin{enumerate}
\item Let $L^1[\eps](\Pi_{b_{max}}):=(\R[\eps]/\eps^3)\otimes L^1(\Pi_{b_{max}})\simeq \eps^0\otimes
L^1(\Pi_{b_{max}}) \oplus \eps^1\otimes
L^1(\Pi_{b_{max}}) \oplus \eps^2 \otimes
L^1(\Pi_{b_{max}}).$
\item Let ${\cal H}[\eps]: L^1[\eps](\Pi_{b_{max}})\to L^1[\eps](\Pi_{b_{max}})$ be
represented by
the operator-valued matrix ${\cal H}[\eps]:= \left(\begin{array}{ccc} {\cal H}_{transport}^0
& 0 & 0 \\ {\cal H}_{nonlocal}^{1,0} & {\cal H}_{transport}^1 & 0 \\ 0 & {\cal H}_{nonlocal}^{2,1} & {\cal H}^2 \end{array}\right).$
\item Let {\em Ev:} $L^1[\eps](\Pi_{b_{max}}) \to L^1(\Pi_{b_{max}})$ be the {\em evaluation mapping},
\BEQ {\mathrm{Ev}}(\eps^0 \otimes h^0 + \eps^1\otimes h^1 + \eps^2 \otimes h^2)(a,b)=
h^0(a,b)+bh^1(a,b)+b^2 h^2(a,b).\EEQ
\end{enumerate}
\end{Definition}

The $(2,2)$-coefficient ${\cal H}^2$ --  the sum
${\cal H}^2_{transport}+{\cal H}^2_{nonlocal}$ -- is coherent with the truncation.
Another possibility (also coherent with Lemma \ref{lem:eps} below, but introducing pointless
complications) would be to consider the un-truncated infinite-dimensional matrix
$\tilde{\cal H}[\eps]:=\left(\begin{array}{ccccc} {\cal H}_{transport}^0
& 0 & 0 & 0 & \cdots \\ {\cal H}_{nonlocal}^{1,0} & {\cal H}_{transport}^1 & 0 & 0 & \cdots \\ 0 & {\cal H}_{nonlocal}^{2,1} & {\cal H}^2 & 0 & \cdots \\ 0 & 0 & {\cal H}^{3,2}_{nonlocal} & {\cal H}^3_{transport} & 0 \\ \vdots & \vdots & \vdots & \vdots & \vdots \end{array}\right)$ acting on $\tilde{L}^1[\eps](\Pi_{b_{max}})\equiv \R[\eps]\otimes
L^1(\Pi_{b_{max}})$, with evaluation mapping Ev$(\sum_{j\ge 0} \eps^j\otimes h^j)(a,b)=\sum_{j\ge 0} b^j h^j$.

\begin{Lemma} \label{lem:eps}
Let $(h_t)_{0\le t\le T}\in L^1[\eps](\Pi_{b_{max}})$ be the solution of the time-evolution
problem $\frac{dh_t}{dt}={\cal H}[\eps](t)h_t$ with terminal condition $h_T\equiv \eps^0\otimes h_T$. 

\noindent Then $f_t:={\cal C}^0\circ {\mathrm{Ev}}(h_t)$ solves (\ref{eq:PDE-f}) with
initial condition ${\cal C}^0(h_T)$.
\end{Lemma}

\noindent {\bf Proof.} By definition.

\hfill \eop

\medskip\noindent Thus {\em our time-evolution operator is ${\cal H}[\eps]$}. Let ${\cal H}_{transport}[\eps]:=\left(\begin{array}{ccc}  {\cal H}_{transport}^0
& 0 & 0 \\ 0 & {\cal H}_{transport}^1 & 0 \\ 0 & 0 & {\cal H}_{transport}^2 \end{array}\right)$ and
${\cal H}_{nonlocal}[\eps]:=\left(\begin{array}{ccc} 0
& 0 & 0 \\ {\cal H}_{nonlocal}^{1,0} & 0 & 0 \\ 0 & {\cal H}_{nonlocal}^{2,1} &{\cal H}^2_{nonlocal} \end{array}\right)$.
The Green function first-order expansion then reads as follows,

\BEQ U[\eps](t,T)=U_{transport}[\eps](t,T)-\int_t^T ds\, U[\eps](t,s) {\cal H}_{nonlocal}[\eps](s)
 U_{transport}[\eps](s,T), \label{eq:Green-function-expansion} \EEQ
 
$U[\eps](t,T)$, resp. $U_{transport}[\eps](t,T)$ being the Green kernels (or evolution operators)
obtained by integrating the time-inhomogeneous evolution systems generated by ${\cal H}[\eps]$,
resp. ${\cal H}_{transport}[\eps]$, i.e. $u(t)=U[\eps](t,T)u(T)$, resp. $u_{transport}(t)=U_{transport}[\eps](t,T)u_{transport}(T)$, 
 solves the linear equation $\frac{du(t)}{dt}={\cal H}[\eps](t)u(t)$, resp. $\frac{du_{transport}}{dt}={\cal H}_{transport}[\eps](t)u_{transport}(t)$.  We shall actually require
  a {\em second-order} expansion of  the Green kernel, obtained by iterating
 (\ref{eq:Green-function-expansion}),
\BEA &&   U[\eps](t,T)=U_{transport}[\eps](t,T)-\int_t^T ds\, U_{transport}[\eps](t,s) {\cal H}_{nonlocal}[\eps](s)
 U_{transport}[\eps](s,T) \nonumber\\
 && \qquad + \int_t^T ds \int_t^s ds'\, U[\eps](t,s'){\cal H}_{nonlocal}[\eps](s')
 U_{transport}[\eps](s,'s){\cal H}_{nonlocal}[\eps](s)U_{transport}[\eps](s,T) \label{eq:2Nd-Green-function-expansion-eps} \nonumber\\
 \EEA
 
Thus (considering a terminal condition $h_T\equiv \eps^0\otimes h_T$)
\BEA &&   ({\mathrm{Ev}}\circ U[\eps](t,T))(h_T)(a,b)=U^0
_{transport}(t,T)h_T(a,b) \nonumber\\
 && \qquad  - b \int_t^T ds\, U_{transport}^1(t,s) {\cal H}_{nonlocal}^{1,0}(s)
 U_{transport}^0(s,T) h_T(a,b) \nonumber\\
 && \qquad + b^2\int_t^T ds \int_t^s ds'\, U^2(t,s'){\cal H}_{nonlocal}^{2,1}(s')
 U_{transport}^1(s,'s){\cal H}_{nonlocal}^{1,0}(s)U_{transport}^0(s,T)h_T(a,b). \label{eq:2Nd-Green-function-expansion} \nonumber\\
 \EEA
 
\medskip\noindent {\em Define
\BEQ ||f||_{(L^1\cap L^{\infty})(\Pi_{b_{max})}}:=\sup\left( ||f||_{L^1(\Pi_{b_{max}})} ,
 ||f||_{L^{\infty}(\Pi_{b_{max}})} \right)  \EEQ
and, for an operator ${\cal H}:(L^1\cap L^{\infty})(\Pi_{b_{max}}) \to
(L^1\cap L^{\infty})(\Pi_{b_{max}})$,
\BEQ |||{\cal H}|||_{(L^1\cap L^{\infty})(\Pi_{b_{max}})}:= \sup_{ ||f||_{(L^1\cap
L^{\infty})(\Pi_{b_{max}})}=1}   ||{\cal H}f|| _{(L^1\cap
L^{\infty})(\Pi_{b_{max}})}.\EEQ }

 From the estimates proved in section 3, that is,
from Lemma  \ref{lem:U-transport} on the one hand, and Lemmas \ref{lem:S8}, \ref{lem:S3},  \ref{lem:S9} on the other, 
 we know that, for all $\kappa\in\N$ and $0\le t\le s$:
\BEQ |||U^{\kappa}_{transport}(t,s)|||_{L^1(\Pi_{b_{max}})\to L^1(\Pi_{b_{max}})},
 |||U^{\kappa}_{transport}(t,s)|||_{L^{\infty}(\Pi_{b_{max}})\to L^{\infty}(\Pi_{b_{max}})}\le 1; \EEQ
 
\BEQ   |||{\cal H}_{nonlocal}^{\kappa+1;\kappa}|||_{
(L^1\cap L^{\infty})(\Pi_{b_{max}}) \to (L^1\cap L^{\infty})(\Pi_{b_{max}})} = O( ||V'||_{8+\kappa,[-3R,3R]} ),\EEQ

\BEQ    |||{\cal H}_{nonlocal}^{\kappa}|||_{(L^1\cap L^{\infty})(\Pi_{b_{max}}) \to
(L^1\cap L^{\infty})(\Pi_{b_{max}})} = O(||V'||_{7+\kappa,[-3R,3R]} )
\label{eq:Green-estimates-H-kappa} \EEQ

\noindent Hence

\begin{Lemma}  \label{lem:L1Linfty-norm}
Let $T>0$ fixed, and $0\le t\le s\le T$. Then 
\BEQ 
|||U^{\kappa}(t,s)|||_{(L^1\cap L^{\infty})(\Pi_{b_{max}}) \to
 (L^1\cap L^{\infty})(\Pi_{b_{max}})}\le e^{c||V'||_{7+\kappa,[-3R,3R]} t}
\EEQ
  for some constant $c>0$.  
\end{Lemma}

\noindent{\bf Proof.} Results from (\ref{eq:Green-estimates-H-kappa}),  Tanabe \cite{Tan}, Theorem 4.4.1 (construction of fundamental
solutions of temporally inhomogeneous equations) and Proposition 4.3.3 (bounded perturbations
of generators of "stable" strongly continuous semi-groups, here of ${\cal H}^{\kappa}_{transport}$
by ${\cal H}^{\kappa}_{nonlocal}$). \hfill \eop
 
\bigskip 
\noindent
{\em Proof of (\ref{eq:4.14bis}).}  We start from It\^o's formula (\ref{eq:Ito}),
\BEQ d\langle Y_t^N,f_t\rangle=\half(1-\frac{\beta}{2}) \langle X_t^N,f''_t\rangle dt +
\frac{1}{\sqrt{N}} \sum_{i=1}^N f'_t(\lambda_t^i)dW_t^i\EEQ
where 
\BEQ f_{\tau}(x)=\frac{\chi_R(x)}{x-z_{\tau}}, \qquad z_{\tau}\equiv z, 
\EEQ represented as $({\cal C}^{\kappa} h_{\tau})(x)$
 for some $\kappa$ (to be chosen later), with $h_{\tau}\ge 0 $ defined
as in (\ref{eq:fz=Ch}), and $f_t$ satisfies the finite-$N$ evolution equation
(\ref{eq:intro-PDE-f}). Recall (see (\ref{eq:hL1b-1},\ref{eq:hLinftyb-3})) that $||h_{\tau}||_{L^1(\Pi_{b_{max}})}\approx 1/b_{\tau}^{1+\kappa}$,
$||h_{\tau}||_{L^{\infty}(\Pi_{b_{max}})}\approx 1/b_{\tau}^{3+\kappa}$.

\medskip\noindent Integrating, we must bound three terms:

\begin{itemize}
\item[(1)] (initial condition) $\esper|\langle Y_0^N,f_0\rangle|^2$;
\item[(2)] (drift term) $\esper \left(\int_0^{\tau} dt\, |\langle X_t^N, f''_t\rangle|
\right)^{2}$;
\item[(3)] ("martingale term") $\esper\left( \frac{1}{\sqrt{N}} \int_0^{\tau} dt\, \sum_{i=1}^N f'_t(\lambda_t^i)
dW^i_t \right)^2$.
\end{itemize} 
 
where $f_t={\cal C}^{\kappa} h_t$. 
 
\bigskip 
\noindent 
Bounding (1) is  easy.   We use the $0$-th order Stieltjes decomposition, \\ $f_0(x)=({\cal C}^0 h_0)(x)=\int
da \int_{-b_{max}}^{b_{max}} (-\II b)\, db\,   {\mathfrak{f}}_z(x) \ h_0(a,b)$, together
with the Cauchy-Schwarz inequality, and obtain as in (\ref{eq:bound-YtNphi2})

\BEQ  |\langle Y_0^N,f_0\rangle|^2 \le ||h_0||^2_{L^2(\Pi_{b_{max}})} \int_{\Pi_{b_{max}}} da\, db\, b^2 (N |M_0(z)-M_0^N(z)|)^2  \label{eq:bound-(1)}
 \EEQ

We use the obvious $L^1-L^{\infty}$-bound, $||h_0||^2_{L^2(\Pi_{b_{max}})}\le 
||h_0||_{L^1(\Pi_{b_{max}})} ||h_0||_{L^{\infty}(\Pi_{b_{max}})}$,  and 
Lemma \ref{lem:L1Linfty-norm},
\BEQ  ||h_0||_{L^1(\Pi_{b_{max}})} ||h_0||_{L^{\infty}(\Pi_{b_{max}})}
\le C ||h_{\tau}||^2_{(L^1\cap L^{\infty})(\Pi_{b_{max}})} 
\le  C'/b_{\tau}^{6} \label{eq:b8}
\EEQ
with $C,C'$ are constants depending on $T$, $R\equiv R(T)$ and  $||V'||_{7,[-3R,3R]}$.

 \medskip\noindent There remains to bound the integral in the r.-h.s. of 
(\ref{eq:bound-(1)}), using our Assumption (iii) on the initial measure, see
(\ref{hyp:initial-measure}).  We split $\int_{\Pi_{b_{max}}}(\cdots)$  into
$\int_{[-2R,2R]\times[-b_{max},b_{max}]}(\cdots)+\int_{(\R\setminus[-2R,2R])\times
[-b_{max},b_{max}]}(\cdots)$. The integral over
 $[-2R,2R]\times[-b_{max},b_{max}]$ is $O(1)$.
As for the integral over $z\in(\R\setminus[-2R,2R])\times[-b_{max},b_{max}]$, we first remark that 
\BEQ N|M_0^N(z)-M_0(z)|= \langle Y_0^N,\chi_R {\mathfrak{f}}_z\rangle \EEQ
and
\BEQ \chi_R(x){\mathfrak{f}}_z(x)=\int da' \int_{-b_{max}}^{b_{max}} (-\II b')\, db'\, 
\frac{1}{x-z'} {\cal K}^0_{b_{max}}(\chi_R {\mathfrak{f}}_z)(a') \EEQ
with (see (\ref{eq:calK-kappa-even}))  ${\cal K}^0_{b_{max}}(\chi_R {\mathfrak{f}}_z)(a')=O(\frac{1}{|a|})
O(\frac{1}{1+
|a'|^2})$. Hence (using Cauchy-Schwarz's inequality and (\ref{hyp:initial-measure}) once
again)

\BEA  && \esper\  \int_{|a|>2R} da \, \int_{|b|<b_{max}} db\, b^2\  \Big( N|M_0^N(z)-M_0(z)| \Big)^2  \nonumber\\
&& \le C\   \esper \int_{(\R\setminus[-2R,2R])\times[-b_{max},b_{max}]} \frac{b^2\,  da\, db}{|a|^2} 
\left|  \int da' \int_{-b_{max}}^{b_{max}} db'\, 
\, \frac{N|b'|\ |M_0^N(z')-M_0(z')|}{1+|a'|^2}  \right|^2 \nonumber\\
&&\le C'  \int \frac{da'}{1+|a'|^2} \int_{-b_{max}}^{b_{max}}  db'\, 
\esper[( N|b'| \ |M_0^N(z')-M_0(z')|)^2] = O(1).\EEA

All together we have proved: $\esper|\langle Y_0^N,f_0\rangle|^2=O(b_T^{-6}).$

\bigskip\noindent
The bound for (2) is essentially pathwise, but more subtle and relies on our perturbative expansion for the
Green kernel, which yields the optimal exponent of $1/b_{\tau}$.  {\em Main term} is obtained by replacing $f''_t$ in (2) with
$({\cal C}^0 u_t)''$, where $u_t:=U_{transport}(t,\tau)h_{\tau}$. By assumption $h_{\tau}\ge 0$, so (see (\ref{eq:utilde})) $\tilde{u}_t:=\tilde{U}_{transport}(t,{\tau})h_{\tau}=
|u_t|\ge 0$
for $0\le t\le {\tau}$, and  these terms
may be bounded as in Israelsson in a probabilistic way, by using the characteristic estimates
proved in section 3.

\medskip\noindent   First (using $|({\mathfrak{f}}_z)''(x)|\le \frac{2}{|x-z|^3}\le \frac{2}{|b|}
\frac{1}{|x-z|^2}$), we have pathwise
\BEA &&  \int_0^{\tau} dt\, |\langle X_t^N, ({\cal C}^0 u_t)''\rangle| \le \int_0^{\tau} dt\, 
\left| \Big\langle X_t^N, \int da \int_{-b_{max}}^{b_{max}} db\, |b|  |({\mathfrak{f}}_z)''(x)| \tilde{u}_t(a,b) \Big\rangle \right|
\nonumber\\
&&\qquad\qquad\le  2 \int_0^{\tau} dt \int_{-3R}^{3R} da \int_{-b_{max}}^{b_{max}} db\, \frac{1}{b} \Im( M_t^N(z)) \, \tilde{u}_t(a,b)
\nonumber\\
&&\qquad\qquad\le  2 \int_0^{\tau} dt \int_{-3R}^{3R} da \int_{-b_{max}}^{b_{max}} db\, \frac{1}{b} \Im \Big[(M_t^N+M_t)(z)\Big] \tilde{u}_t(a,b)
\nonumber\\
&&\qquad\qquad= 2  \int_0^{\tau} dt \int_{-3R}^{3R} da \int_{-b_{max}}^{b_{max}} db\, \partial_b\Big[ \Im \Big[(M_t^N+M_t)(z)\big] \tilde{u}_t(a,b)
\Big] \ln(1/|b|) + {\mathrm{bdry}}_1
\nonumber\\
 \label{eq:(2)}\EEA
where 
\BEQ {\mathrm{bdry}}_1:=-2\int_0^{\tau} dt\int_{-3R}^{3R} da \, \left[\Im(M_t^N+M_t)(a+\II b_{max}) \tilde{u}_t(a,b_{max}) - (b_{max}\leftrightarrow -b_{max}) \right]
\ln(1/|b_{max}|)\EEQ
is a boundary term.

\medskip\noindent Now, we compare the first term in the r.-h.s. of (\ref{eq:(2)})  to 
\BEQ \frac{8}{\beta}\int_{-3R}^{3R} da\int_{-b_{max}}^{b_{max}} db\, (\tilde{u}_{\tau}(a,b)-\tilde{u}_0(a,b)) \ln(1/|b|) 
\equiv  \frac{8}{\beta}\int_0^{\tau} dt\, \frac{d}{dt} \left( \int_{-3R}^{3R} da\int_{-b_{max}}^{b_{max}}
db\, \tilde{u}_t(a,b) \ln(1/|b|) \right).  \label{eq:8/beta1} \EEQ
The main terms in $\frac{8}{\beta}\frac{d}{dt} \tilde{u}_t(a,b)$ are
those due to the $(\frac{1}{x})$-kernel, 
 $$2\left[ \partial_a\left(  \Re(M^N_t+M_t)(z)\tilde{u}_t(a,b) \right)  + \partial_b\left( \Im(M^N_t+M_t)(z)
\tilde{u}_t(a,b) \right) \right].$$ The horizontal drift term 
$\partial_a \, \Re(M^N_t+M_t)$ vanishes by integration by parts up to a boundary term,
\BEQ  {\mathrm{bdry}}_2:= 2\int_0^{\tau} dt \int_{-b_{max}}^{b_{max}}  db\, \ln(1/|b|)
\left[\Re(M_t^N+M_t)(3R+\II b)\tilde{u}_t(3R+\II b) - (R\leftrightarrow -R) \right] 
\EEQ
(note that $\Re(M_t^N+M_t)(\pm 3R+\II b)=O(1/R)$ is bounded), 
 while the vertical drift
term is identical to (\ref{eq:(2)}). The other terms in ${\cal H}_{transport}(t)$, see
(\ref{eq:da(tau)/dtau}), (\ref{eq:db(tau)/dtau}), 
\BEQ \partial_a \left( (V'(a)-\half V'''(a)b^2)\tilde{u}_t(a,b)\right)+ \left(\partial_b \left(
V''(a)b \tilde{u}_t(a,b)\right)\right) + {\mathrm{bdry}} 
\label{eq:4.24} \EEQ
contribute respectively: yet another boundary term,
\BEQ  {\mathrm{bdry}}_3:=\frac{8}{\beta} \int_0^{\tau} dt \int_{-b_{max}}^{b_{max}}
db\, \ln(1/|b|) \left[ (V'(3R)-\half V'''(3R)b^2)\tilde{u}_t(3R,b) - (R\leftrightarrow -R) \right];\EEQ
and
 \BEA && \frac{8}{\beta} \int_0^{\tau} dt\, 
\int_{-3R}^{3R} da \, V''(a) \int_{-b_{max}}^{b_{max}} db\,  \tilde{u}_t(a,b)
\nonumber\\
&& \qquad \qquad\qquad=O\left(\sup_{[-3R,3R]} |V''|\right) \ e^{c||V'||_{7,[-3R,3R]}  T} ||u_T||_{L^1(\Pi_{b_{max}})} \EEA
by Lemma \ref{lem:L1Linfty-norm}, plus a boundary term,
\BEQ  {\mathrm{bdry}}_4:=\frac{8}{\beta} \int_0^{\tau} dt \int_{-3R}^{3R} da\, 
V''(a)  \left[ b_{max} \ln(1/|b_{max}|) \tilde{u}_t(a,b_{max})  - (b_{max}\leftrightarrow -b_{max}) \right],\EEQ 
and also the integral over the domain $[-3R,3R]
\times [-b_{max},b_{max}]$  of  ($\frac{8}{\beta} \ln(1/|b|)$ times the boundary terms of \S 3.9.  Finally, the contribution of the $\tilde{c}$-characteristic is
known from the contraction property to be of the form
 $\frac{8}{\beta}\int_0^{\tau} dt\, 
\int_{3R}^{3R} da \int_{-b_{max}}^{b_{max}} db\, \sigma_t(a,b) \tilde{u}_t(a,b) \ln(1/|b|)\ge 0$ with $\sigma_t(\cdot,\cdot)\ge 0$, hence positive.

\medskip\noindent  Using the $(L^1\cap L^{\infty})$-bound of $\tilde{u}_t$, one sees
that all boundary terms are\\ $O(1) ||u_{\tau}||_{(L^1\cap L^{\infty})(\Pi_{b_{max}})}=
O(|b_{\tau}|^{-3})$,  times some derivative of $V$, $||V^{(j)}||_{0,[-3R,3R]}$, $j=1,2,3$, times 
possibly $\int_{-b_{max}}^{b_{max}} db\, \ln(1/|b|)=O(1)$.  But, actually, we have
a much better bound {\em for $T$ small enough}: because  $h_{\tau}(z)$ is $O(1)$, independent of
$b_{\tau}$, far from the support $[-R,R]\times\{0\}$, say, on $\Pi_{b_{max}}\setminus\Big(
[-2R,2R]\times[-\half b_{max},\half b_{max}]\Big)$, we shall have 
\BEQ ||u_t||_{L^{\infty}(\partial_R
\Pi_{b_{max}})}=O(1) \label{eq:bd-bdry}
\EEQ
 for
all $t\in[0,\tau]$, as explained in the side remark before Lemma  \ref{lem:U-transport}.
Anticipating on the next terms featuring in the second-order expansion of the Green
kernel (see below (\ref{eq:4.43})), it is easy to see that ${\cal H}^{1,0}_{nonlocal}
(s)u_s$, whence $U^1_{transport}(t,s)(|{\cal H}^{1,0}_{nonlocal}
(s)u_s|)$ and $U^2(t,s'){\cal H}_{nonlocal}^{2,1}(s')
(U^1_{transport}(s',s){\cal H}_{nonlocal}^{1,0}(s) u_s$ too, enjoy the same property (\ref{eq:bd-bdry}). Incidentally, this implies $||h_0||_{L^1(\Pi_{b_{max}})}
||h_0||_{L^{\infty}(\Pi_{b_{max}})} \le C ||h_{\tau}||_{L^1(\Pi_{b_{max}})} 
||h_{\tau}||_{L^{\infty}(\Pi_{b_{max}})} \le C'/b_{\tau}^4$ instead of $C'/b_{\tau}^6$
in (\ref{eq:b8}).

\medskip\noindent Consider now the {\em left-hand side} of (\ref{eq:8/beta1}). Considering the {\em adjoint evolution}, we get a time-reversed sub-Markov process 
$(\tilde{A}_t,\tilde{B}_t)$ with kernel $p(s,\tilde{a}_s;t,\tilde{z}_t)$, $s\le t$
(see \S 3.6).
 Since $t\mapsto |\tilde{B}_t|$ {\em decreases},  we
obtain
\BEA && \left| \int da\int_{-b_{max}}^{b_{max}} db\,  \tilde{u}_0(a,b) \ln(1/|b|)\right| =
\int da_{\tau} \int db_{\tau}\, u_{\tau}(a_{\tau},b_{\tau})  \nonumber\\
&& \qquad\qquad\qquad\qquad \int da_0 \int db_0\,  p(0,\tilde{a}_0,\tilde{b}_0;{\tau},a_{\tau},b_{\tau}) \ln(1/|\tilde{b}_0|) \nonumber\\
&&\le  \int da_{\tau} \int db_{\tau}\,   \ln(1/|b_{\tau}|)  u_{\tau}(a_{\tau},b_{\tau})
= O(\ln(1/|b_{\tau}|)b_{\tau}^{-1}) \nonumber\\  \label{eq:4.43}
\EEA
by the log-estimate (\ref{eq:log-estimate}). 
So much for the contribution of $U^0_{transport}$ to (2), which we have
shown to be overall $O((b_{\tau}^{-3})^2)=O(b_{\tau}^{-6})$, and even 
$O((\ln(1+1/b_{\tau})/b_{\tau})^2)$ for $T$ small enough, as in \cite{Isr}. 

\bigskip\noindent We now use the second-order 
expansion of the Green kernel (\ref{eq:2Nd-Green-function-expansion}). The {\em second
term} in the expansion,
\BEQ v(a,b):= \int_t^{\tau} ds\, U_{transport}^1(t,s) {\cal H}_{nonlocal}^{1,0}(s)
 U_{transport}^0(s,\tau) h_{\tau}(a,b) \EEQ
leads to a development similar to (\ref{eq:(2)}):

\BEA &&\int_0^{\tau} dt\, |\langle X_t^N, ({\cal C}^0 v_t)''\rangle|  \nonumber\\
 &&\qquad \le \int_0^{\tau} dt\, 
\left| \langle X_t^N, \int_{-3R}^{3R} da \int_{-b_{max}}^{b_{max}} db\, b^2  |({\mathfrak{f}}_z)''(x)|\int_t^{\tau} ds
(U^1_{transport}(t,s)|{\cal H}_{nonlocal}^{1,0}(s) u_s|)(a,b)\rangle \right|
\nonumber\\
&& \qquad \le  2 \int_0^{\tau} ds \left[ \int_0^s dt \int_{-3R}^{3R} da \int_{-b_{max}}^{b_{max}} db\,  |\Im (M_t^N+M_t)(z)| \tilde{u}^s_t(a,b) \right]
\nonumber\\
&&\qquad = -2 \int_0^{\tau} ds \left[ \int_0^s dt \int_{-3R}^{3R} da \int_{-b_{max}}^{b_{max}} db\, \partial_b\left[ \Im (M_t^N+M_t)(z) \tilde{u}^s_t(a,b)
\right] b \right] + {\mathrm{bdry}}'_1
\nonumber\\
 \label{eq:(2)bis}\EEA
where $\tilde{u}^s_t(a,b):=\tilde{U}^1_{transport}(t,s) (|{\cal H}^{1,0}_{nonlocal}(s)u_s|)(a,b)$ $(\ge 0)$,
which we compare to
\BEA && \int_0^{\tau} ds \left[ \frac{8}{\beta}\int_{-3R}^{3R} da\int_{-b_{max}}^{b_{max}} db\, (\tilde{u}_s^s(a,b)-\tilde{u}_0^s(a,b)) b  \right]
\nonumber\\
&& \qquad 
\equiv  \int_0^{\tau} ds \left[ \frac{8}{\beta}\int_0^s dt\, \frac{d}{dt} \left( \int_{-3R}^{3R} da\int_{-b_{max}}^{b_{max}}
db\, \tilde{u}_t^s(a,b) b \right) \right].  \label{eq:8/beta1bis} \nonumber\\
 \EEA
The right-hand side in (\ref{eq:8/beta1bis}) decomposes  in the same way as explained below
(\ref{eq:4.24}) -- but with $\kappa=1$ now --. Compared to the main term studied in the
previous two pages, $\ln(b)$ has been replaced with $b$ (which may simply be
bounded by a constant, $b\le \half$), so logarithms disappear in the estimates,
while the replacement of $u_s$ by $|{\cal H}^{1,0}_{nonlocal}(s)u_s|$ produces
 the supplementary factor $|||{\cal H}^{1,0}_{nonlocal}(s)|||_{(L^1\cap L^{\infty})(\Pi_{b_{max}})} =O(1)$. The total contribution to (2) is therefore  $O(b_{\tau}^{-6})$
 or even $O((\ln(1+1/b)/b)^2)$  as for the main term. 

\medskip\noindent
  The {\em last} term in the Green kernel expansion,

 \BEQ w(a,b):= \int_t^{\tau} ds \int_t^s ds'\, U^2(t,s'){\cal H}_{nonlocal}^{2,1}(s')
 U_{transport}^1(s,'s){\cal H}_{nonlocal}^{1,0}(s)U_{transport}^0(s,{\tau})h_{\tau}(a,b),\EEQ
leads now to a third contribution which is bounded in a very simple way,

\BEA &&\int_0^{\tau} dt\, \Big| \langle X_t^N, ({\cal C}^0 w_t)''\rangle \Big| \le \int_0^{\tau} dt\, 
\left| \langle X_t^N, \int_{-3R}^{3R} da \int_{-b_{max}}^{b_{max}} db\, |b|^3  |( {\mathfrak{f}}_z)''(x)| \right. \nonumber\\
 && \qquad\qquad\qquad \left. \int_t^{\tau} ds
 \int_t^s ds'\, \left|U^2(t,s'){\cal H}_{nonlocal}^{2,1}(s')
(U^1_{transport}(s',s){\cal H}_{nonlocal}^{1,0}(s) u_s)(a,b)\right| \big{\rangle} \right|  \label{eq:4.41} \nonumber\\
\EEA
Since $|b|^3 |({\mathfrak{f}}_z)''(x)|=O(1)$,  (\ref{eq:4.41}) is  simply bounded in the end for arbitrary $T$
by\\  $||u_{\tau}||_{(L^1\cap L^{\infty})(\Pi_{b_{max}})}=O(|b_{\tau}^{-3}|)$, times the product of the $(L^1\cap L^{\infty})(\Pi_{b_{max}})$- operator
norms $|||U^i_{transport}(\cdot,\cdot)|||, |||{\cal H}^{i+1,i}_{nonlocal}(\cdot)|||$ ($i=0,1$) figuring in the integral, yielding
once again a total contribution  $O(b_{\tau}^{-6})$, or  (for $T$ short enough) $O((\ln(1+1/b)/b)^2)$ to (2). 

\vskip 2cm
\bigskip\noindent
 We finally proceed to 
bound the "martingale term" (3). 

\medskip\noindent A caveat is required here: for finite $N$, $f_t(\cdot)$ is {\em not}
${\cal F}_t$-measurable, since it is obtained by integrating the ordinary
differential equation with {\em random coefficients}  (\ref{eq:intro-PDE-f})
{\em backwards} from time $\tau$ to time $t$. Hence 
\BEQ  \esper \Big[ \Big(\frac{1}{\sqrt{N}} \sum_{i=1}^N I_i(\tau)\Big)^2\Big], \qquad I_i(\tau):=\int_0^{\tau} dt \, f'_t(\lambda_t^i)
dW_t^i \label{eq:Ii}
\EEQ
 {\em cannot} be bounded by $\esper \int_0^{\tau} dt\, (f'_t(\lambda_t^i))^2
$ using standard tools of stochastic calculus. By the way, this points out to a 
mistake in the proof of the estimate for $\esper[|N(M_T^N(z)-M_T(z))|^2]$ given
in Proposition 1 of Israelsson's article. Our arguments below yield  a bound in $O(1/b^{12})$ independently of $V$ -- in particular in the harmonic case, instead of the bound in $O((\ln(1+1/b)/b)^2)$ found by Israelsson. See remark at the end of this section for 
some after-thoughts.

\medskip\noindent The correct way to cope with the stochastic integral
$I_i(\tau)$  (\ref{eq:Ii}) is
the following. Since $s\mapsto f'_s(\lambda)$ is $C^1$ for $\lambda$ {\em fixed},  it can be considered as the
{\em non-adapted}
finite variation part of a semi-martingale, in the extended
definition briefly mentioned below Definition (1.17) of
\cite{RY}, Chapter 4. Hence the  integration-by-parts lemma of standard differential calculus holds,   and we can rewrite $I_i(\tau)$ as $\int_0^{\tau} ds\, f'_0(\lambda_s^i) dW_s^i -\int_0^{\tau} dt \, J_i(t)$, where
\BEQ J_i(t):=\int_t^{\tau} \frac{\partial f'_t}{\partial t}(\lambda_s^i) dW_s^i.\EEQ
Then, {\em considering the standard Stieltjes decomposition of $h_t$ of order $\kappa=3$ this time} (which turns out in the end of the ensuing computations to be the minimum possible order
yielding finite results in the neighbourhood of the real axis),
\BEQ f'_t(\lambda_s^i)= \frac{\partial}{\partial x} {\cal C}^3(h_t)(x)=\int da
\int_{-b_{max}}^{b_{max}} \,  db\, (-\II b)\, |b|^3 ({\mathfrak{f}}_z)'(\lambda_s^i) h_t(a,b)
\EEQ
and
\BEQ J_i(t)=\int da \int_{-b_{max}}^{b_{max}} db\, (-\II b)\, |b|^3  ({\cal H}(t)
h_t)(a,b) J_i^z(t),  \label{eq:Ji} 
\EEQ
where the stochastic integral 
\BEQ J_i^z(t):=\int_t^{\tau} ( {\mathfrak{f}}_z)'(\lambda_s^i)
dW_s^i
\EEQ
 is now  a standard (i.e. non-anticipative) It\^o integral, whence (using primed
 integration variables $t',a',b'$ for  $I_{i'}(\tau)$ in the
 averaged squared quantity $\esper\Big[ (\int_0^{\tau} dt\, \sum_i J_i(t))
 (\int_0^{\tau} dt'\, \sum_{i'} J_{i'}(t')) \Big]$) 
\BEQ \esper[ J_i^z(t) J_j^{z'}(t')]=\del_{i,j} \int_{\max(t,t')}^{\tau}  ds\,  ( {\mathfrak{f}}_z)'(\lambda_s^i) ( {\mathfrak{f}}_{z'})'(\lambda_s^i) . \label{eq:JJ}
\EEQ

\medskip\noindent
The first step consists in transferring to the $\mathfrak{f}'_z$-factors the derivatives $\partial_a,\partial_b$ coming
from the action of ${\cal H}_{transport}(t)$ on $h_t$. We concentrate on the
most singular terms coming from ${\cal H}_0^{\kappa}(t)$, namely, $({\cal H}_0^{\kappa}(t)h_t)(a,b)=
\frac{\beta}{4} \Big[ \partial_a (\Re (M_t+M_t^N)(z))h(t;a,b)) +  \partial_b (\Im (M_t+M_t^N)(z))h(t;a,b)) \Big]+\cdots$, where the missing order 0 part $(\cdots)$ is as in
(\ref{eq:H0}).  Integrating these two terms w.r. to the measure $\int da\, db\, 
(-\II b) |b|^3 \, {\mathfrak{f}}'_z(\lambda_s^i)$ yields by integration by parts
$-\frac{\beta}{4} \int da\, (-\II db)\, h_t(a,b)   \Big( \pm b^4 \Re (M_t+M_t^N)(z) {\mathfrak{f}}''_z(\lambda_s^i) \pm \II \Im(M_t+M_t^N)(z) \partial_b
\Big(  b^4 {\mathfrak{f}}'_z(\lambda_s^i) \Big) \Big)$.

\medskip\noindent For finite $N$, ${\cal H}(t)h_t(\cdot)$ is random and not ${\cal F}_t$-measurable, hence\\  $\esper\Big[ ({\cal H}(t)h_t(a,b) {\cal H}(t')h_{t'}(a',b')) 
(J_i^z(t) J_{i'}^{z'}(t')) \Big]$ may not directly be bounded using (\ref{eq:JJ}) 
(see Remark below). Instead, we use the bounds
\BEQ |h_t(a,b)|,|h_{t'}(a',b')|=O(||h_{\tau}||_{L^{\infty}(\Pi_{b_{max}})})=O(b_{\tau}^{-3-\kappa})=O(b_{\tau}^{-6}) \EEQ
\BEQ |b^4 \Re(M_t+M_t^N)(z)|=O(|b|^{3}) , \qquad ||{\mathfrak{f}}''_z||_{\infty}\le |b|^{-3}
\EEQ
and
\BEA && \esper \Big|  \sum_{i,i'} \int_t^{\tau} {\mathfrak{f}}''_z(\lambda_s^i) dW_s^i
\ \int_{t'}^{\tau} {\mathfrak{f}}''_{z'}(\lambda_s^{i'}) dW_s^{i'} \Big| \nonumber\\
&&\le  \Big[\esper \Big( \sum_i  \int_t^{\tau} {\mathfrak{f}}''_z(\lambda_s^i) dW_s^i 
\Big)^2  \Big]^{1/2}  
 \Big[\esper \Big( \sum_{i'}  \int_{t'}^{\tau} {\mathfrak{f}}''_{z'}(\lambda_s^{i'}) dW_s^{i'} 
\Big)^2  \Big]^{1/2}  \nonumber\\
&&=   \Big[ \sum_i \esper   \int_t^{\tau} ds\,  |{\mathfrak{f}}''_z(\lambda_s^i)|^2  
  \Big]^{1/2}  
 \Big[\sum_{i'} \esper   \int_{t'}^{\tau}ds\,  |{\mathfrak{f}}''_{z'}(\lambda_s^{i'}|^2 
  \Big]^{1/2}  \nonumber\\
&&\le \half \esper\Big[ \left(\frac{b}{b'}\right)^3  \sum_i \int_t^{\tau} ds\, |{\mathfrak{f}}''_z(\lambda_s^i)|^2  +  \left(\frac{b'}{b}\right)^3  \sum_i \int_{t'}^{\tau} ds\, |{\mathfrak{f}}''_{z'}(\lambda_s^i)|^2 \Big] \nonumber\\
&&= \frac{N}{2} \esper\Big[ \left(\frac{b}{b'}\right)^3 \int_t^{\tau} ds\,  
\langle X_s^N,({\mathfrak{f}}''_z)^2\rangle +  \left(\frac{b'}{b}\right)^3 \int_t^{\tau} ds\,  
\langle X_s^N,({\mathfrak{f}}''_{z'})^2\rangle \Big] \nonumber\\
&&= O(N (|bb'|)^{-3})
\EEA
Considering instead the terms of the type $4b^3 {\mathfrak{f}}'_z(\lambda_s^i)$ coming from
$\partial_b(b^4  {\mathfrak{f}}'_z(\lambda_s^i))$,  or those coming from the missing order 0 part $(\cdots)$
above, leads to the same scaling in
$b,b'$ when $b,b'\to 0$, as can easily be seen, while terms coming from the
bounded operator ${\cal H}_{nonlocal}$ or from the time $0$ contribution $\int_0^{\tau}
ds\, f'_0(\lambda^i_s)dW^i_s$ are less singular.
Thus we finally find, as expected:
\BEQ  \esper \Big[ \Big(\frac{1}{\sqrt{N}} \sum_{i=1}^N I_i(\tau)\Big)^2\Big] \le 
 \Big(\int_{-3R}^{3R} da\,  \int_{-b_{max}}^{b_{max}} db\, O(b_{\tau}^{-6}) \Big)^2 =
 O(b_{\tau}^{-12}),  \label{eq:poor-estimates}
\EEQ
plus an $O(1)$-contribution coming from $({\cal H}(t)h_t)_{ext}$.
\hfill \eop

\bigskip\noindent  {\bf Remark.} Israelsson's bounds in $O((\ln(1+1/b)/b)^2)$ are
recovered if one (somewhat carelessly and out of the blue!) {\em replaces} 
$f_t$, solution of the finite-N evolution equation (\ref{eq:intro-PDE-f}), with
the {\em deterministic} solution $f^{\infty}_t$ of the asymptotic evolution equation (\ref{eq:intro-PDE-f-asymptotic}).    Namely, in that case, It\^o's formula applies, see
(\ref{eq:Ii}).
 Using $|({\mathfrak{f}}_z)'(x)|\le \frac{1}{|b|} \frac{1}{|x-z|}$ and (\ref{eq:Im-rho-infty}), we get, letting $f^{\infty}_t={\cal C}^0 (h^{\infty}_t)$ and 
 $I_i^{\infty}(\tau):=\int_0^{\tau} dt\, (f_t^{\infty})'(\lambda_t^i) dW_t^i$ : 
\BEA &&  \esper \Big[ \Big(\frac{1}{\sqrt{N}} \sum_{i=1}^N I_i^{\infty} (\tau)\Big)^2\Big]= \int_0^{\tau} dt |\langle X_t^N, (({\cal C}^0 h^{\infty}_t)')^2\rangle| dt \nonumber\\
 &&\le \int_0^{\tau} dt \,
 \int da \int_{-b_{max}}^{b_{max}} db\, \ \cdot\ \int da'\int_{-b_{max}}^{b_{max}} db'\,  |h_t(a,b)|
 |h_t(a',b')|\  |\langle X_t^N, |b f'_z(\cdot)|\  |b'f'_{z'}|(\cdot)\rangle | \nonumber\\
 &&\le  \int_0^{\tau} dt \,
 \int da \int_{-b_{max}}^{b_{max}} db\, \ \cdot\ \int da'\int_{-b_{max}}^{b_{max}} db'\,  |h_t(a,b)|
 |h_t(a',b')| \langle X_t^N, b^2 |(f'_z(\cdot))|^2\rangle \nonumber\\
 &&\le \sup_{0\le t\le {\tau}}  \left( \int da'\int_{-b_{max}}^{b_{max}} db'\, |h_t(a',b')|\right)
 \ \cdot\  \int_0^{\tau} dt\, \left( \int da \int_{-b_{max}}^{b_{max}} db\, \frac{1}{b} \Im M_t^N(z) |h_t(a,b)| \right)
. \nonumber\\  \label{eq:4.23} \EEA
The second factor in (\ref{eq:4.23}) is bounded exactly like
the drift term (2), while the first one
is just $||h_t||_{L^1(\Pi_{b_{max}})}$. All together,  $I^{\infty}(\tau)$
is bounded by $ O(b_{\tau}^{-6})$, or even by $O(\ln(1+1/b_{\tau})/b_{\tau}^2)$ for
$T$ small enough.

\medskip\noindent A way to improve our poor estimates (\ref{eq:poor-estimates}) would  be to separate
$h_t$ into $h_t^{\infty}$ plus a $O(1/N)$ fluctuation $\del h_t^{\infty}$, whose
contribution to (\ref{eq:4.14bis}) would be hopefully $O(1/N)$ times some inverse
power of $b$, and would therefore vanish when $N\to\infty$. However, the time-evolution of
$\del h_t^{\infty}$  is a priori governed by the Jacobian of (\ref{eq:intro-PDE-f-asymptotic}) around $h_t^{\infty}$, whose characteristics are obtained by linearizing
those of the transport operator ${\cal H}_{transport}$. Alas, the linearization of the 
already singular characteristics of ${\cal H}_0$, see Proposition \ref{prop:Israelsson},
leads to an exponential factor of the type $\exp\Big(c\int dt\, ((M_t^N)''(Z_t)+(M''_t(Z_t))
\Big)$, which is exponentially large for small $|b|$ near the points $x$ of the real axis
at which $M_t(x\pm\II 0)$ is not differentiable, e.g. near the end points of the support
for a standard density of the semi-circle type $\frac{1}{\pi} \sqrt{2-x^2}$ 
($a>0$), with associated Stieltjes transform $M_t(z)=-z+\sqrt{z^2-2}$, see (\ref{eq:cut-eq2}).


\section{Large deviation bound for the support of the measure} \label{sec:entropy}


As a key technical argument required for the convergence of our scheme, we prove
in this section the following bound for the probability that the support of the measure is large. Since the number $N$ of eigenvalues varies in this
section, we emphasize the $N$-dependence of the process when we judge it
necessary by writing $\lambda_t^{N,i}$ instead of $\lambda_t^i$.

\begin{Lemma}[large deviation bound]  \label{lem:support}
Assume the large deviation estimate (\ref{eq:support-bound0}) for the initial support holds, namely, $\proba[ \max_{i=1,\ldots,N} |\lambda_0^{N,i}| > R_0] \le 
C_0 e^{-c_0 N} $ for some constants $R_0,c_0,C_0>0$.
 Let $T>0$. There exists some radius $R=R(T)$ and constant $c$, depending on $V$ and $R_0,c_0$ but uniform in $N$,  such that
\BEQ \proba\left[ \sup_{0\le t\le T} \sup_{i=1,\ldots,N} |\lambda_t^{N,i}| > R\right]
\le C e^{-cN}. 
\EEQ
\end{Lemma}

\noindent The principle of the proof was obligingly provided by a referee. It relies on uniform-in-time
moment bounds for the empirical measure, and on a comparison principle for sde's.

\medskip\noindent
First, we use as an input moment bounds  proved {\em in the case $V=0$}  by induction on $p=1,2,\ldots,\eps N$ in Anderson-Guionnet-Zeitouni \cite{And}. Let $(\tilde{\lambda}_t^{N,i})_{t\ge 0}, i=1,\ldots,N$ be the solution of the modified system of coupled
stochastic differential equations with {\em zero potential}, 
\BEQ d\tilde{\lambda}_t^{N,i}=\frac{1}{\sqrt{N}} dW_t^i+ \frac{\beta}{2N}
\sum_{j\not=i} \frac{dt}{\tilde{\lambda}^{N,i}_t-\tilde{\lambda}^{N,j}_t}, \qquad i=1,\ldots,N  \label{eq:SDE0} \EEQ
with initial condition $\tilde{\lambda}^{N,i}_0\equiv \lambda^{N,i}_0$ coinciding
with that of (\ref{eq:SDE}), and $\tilde{X}_t^N:=\frac{1}{N} \sum_{i=1}^N \del_{\tilde{\lambda}^{N,i}_t}$ be the corresponding random point process. Under  $\Omega_0  : \left( \max_{i=1,\ldots,N} |\lambda_0^{N,i}| \le  R_0 \right)$,  an event of proability $1-C_0 e^{-c_0 N}$,  
eq. (4.3.45) in \cite{And} holds,  namely,
\BEQ \esper\left[ {\bf 1}_{\Omega_0}\  \sup_{0\le t\le T} \int \tilde{X}^N_t(dx)\, |x|^p \right] \le R_1(T)^p. \EEQ
(An explicit expression for the constants $\eps$ and $R_1(T)$, depending on $R_0$, can be obtained by following
computations on p. 274, as  a consequence of   Lemma 4.3.17.)  The above bound implies
$\esper\left[{\bf 1}_{\Omega_0}\ \sup_{0\le t\le T} \sup_{i=1,\ldots,N} |\tilde{\lambda}_t^{N,i}|^p\right]\le N R_1(T)^p$
and then
(by Markov's inequality), letting $R_2(t):=eR_1(t)$,
\BEQ \proba[{\bf 1}_{\Omega_0}\ \sup_{0\le t\le T} \sup_{i=1,\ldots,N} |\tilde{\lambda}_t^{N,i}|> R_2(T)]\le Ne^{-\eps N}.  \label{eq:LDB-R2} \EEQ 
 
\medskip\noindent Then, one compares the two eigenvalue processes $(\lambda^i_t)_{0\le t\le T}$ and  $(\bar{\lambda}^i_t)_{0\le t\le T}$, adapting the argument given in \cite{And}, Lemma 4.3.6. Let $E_t^{N,i}:=\lambda^{N,i}_t-
\bar{\lambda}^{N,i}_t-\alpha t$ $(\alpha >0)$. Subtracting the  sde's for the two processes, one gets
\BEQ \frac{dE_t^{N,i}}{dt}=-\frac{\beta}{2N} \sum_{j\not=i} \frac{E_t^{N,i}-E_t^{N,j}}{(\lambda^{N,i}_t-\lambda^{N,j}_t)(\bar{\lambda}^{N,i}_t-\bar{\lambda}^{N,j}_t)} - (V'(\lambda^{N,i}_t)-
V'(\bar{\lambda}^{N,i}_t)) - V'(\bar{\lambda}^{N,i}_t)-\alpha .  \label{eq:LDB-E} \EEQ
Whatever the ordering chosen for the eigenvalues, the denominator  in (\ref{eq:LDB-E})
is always $>0$ because eigenvalues never cross.
We assume that the event $\Omega: \left( \max_{i=1,\ldots,N} |\lambda_0^{N,i}| \le  R_0 \right) \cap \left(\sup_{0\le t\le T} \sup_{i=1,\ldots,N} |\bar{\lambda}_t^{N,i}|\le  R_2(T) \right)$ is realized, an event of probability $1- e^{-c_1 N}$; then
$|V'(\bar{\lambda}^{N,i}_t)|$ is bounded uniformly in $N$, $i$, and $t\le T$ by some
constant $C_2$ depending  on $V$; we assume $\alpha>C_2$.   Initially,
$E_0^{N,i}\le R_0+R_2(T)$, $i=1,\ldots,N$ by construction. {\em Assume} that there exists some $t<T$ and $i$ 
such that $E_t^{N,i}\ge R_3(T):=R_0+R_2(T)+1$, and let $t_{min}>0$ be the first time at
which one such inequality holds, so that $E_{t_{min}}^{N,i}=R_3(T)$ for some $i$, while
$E_t^{N,j}<R_3(T)$ for $t<t_{min}$ and $j=1,\ldots,N$. But then $E_{t_{min}}^{N,i}-E_{t_{min}}^{N,j}\ge 0$ for all $j\not=i$, and (by convexity of $V$) $V'(\lambda_{t_{min}}^{N,i})-
V'(\bar{\lambda}_{t_{min}}^{N,i})\ge 0$. Hence $\frac{dE_{t_{min}}^{N,i}}{dt}<0$: a contradiction.  Reversing the signs of the inequalities, one proves similarly that $\tilde{\lambda}_t^{N,i}-\lambda_t^{N,i}-\alpha t\le R_3(T)$. Concluding: with high probability, $\sup_{0\le t\le T} \sup_{i=1,\ldots,N} |\lambda_t^{N,i}|\le R_0+2R_2(T)+
C_2 T+1$. 

\hfill \eop


\section{Appendix. Generalized transport operators} \label{sec:transport}


Many operators in this article are of the following type,
\BEQ {\cal H}_tf(x)=\sum_k v_k(t,x)\partial_{x_k}f(x)+ \tau(t,x) f(x) \EEQ
with $f:\Omega\to \R$, where $\Omega$ is a domain in $\R^d$ (in practice, we need
only consider $\Omega=\Pi^{\pm}$), and $\vec{v}(t,\cdot)$ a vector field, resp.
$\tau(t,\cdot)$ a function,
on $\Omega$. Let us call such operators {\em generalized
transport operators}. 

It is well-known how to solve PDEs generated by generalized transport operators, i.e.
of the type
\BEQ \frac{\partial f_t}{\partial t}(x)={\cal H}_t f_t(x) \label{eq:transport-PDE} \EEQ
with terminal condition $f_T\equiv f$. 
Namely, let $y_t\equiv \Phi_t^T(y)$ (called: {\em characteristics of (\ref{eq:transport-PDE})}) be the solution of the ode $\frac{dy_t}{dt}=\vec{v}(t,y_t)$ with terminal condition $y_T=y$.
One checks immediately  that
\BEQ f_t(y)=c_t f(\Phi_t^T(y)), \qquad c_t:=\exp\left(-\int_t^T \tau(y_s)\, ds\right) \EEQ
is a solution. In particular, $\supp(f_t), t\le T$ is the inverse image of 
$\supp(f)$ by $\Phi_t^T$; so, if $\vec{v}\big|_{\partial\Omega}$ is {\em inward on the boundary of some domain $\Omega$ containing the support of $f_T$, then
$\supp(f_T)\subset\Omega$ for all $t\le T$.} 
In the article we actually refer either to the basis trajectory
$y_{\cdot}=a_{\cdot}+\II b_{\cdot}$ or to the "extended"
trajectory $(a_{\cdot}+\II b_{\cdot},c_{\cdot})$ as characteristics.

\medskip 

The {\em Jacobian} of the ode, $J_t:=\frac{dy_t}{dy}$, solves the linearized ode 
$\frac{dJ_t}{dt}=\nabla \vec{v}(t,y_t)J_t$ with terminal condition $J_T=\Id$. In particular
(letting $|\, \cdot\, |$ denote the determinant),
$\frac{d}{dt}\big|_{t=T} |J_t|=\Tr (\nabla \vec{v}(T,y))=\nabla\cdot \vec{v}(T,y).$
The time-variation of the $L^1$-norm of $f_t$ is
\BEA   \frac{d}{dt}\big|_{t=T} \int dy\, |f_t(y)| &=& \int dy\, \left(\Re \frac{d}{dt}
 c_t(y) \right)|f(y)| - \int dy\,  \left(\frac{d}{dt}\big|_{t=T} |J_t|\right)
 |f(y)| \nonumber\\
 &=& \int dy \left(\Re \tau(T,y)-\Tr (\nabla \vec{v}(T,y))\right) |f(y)|; \label{eq:time-variation-L1-norm}
\EEA

it vanishes when $\Re c=\Tr (\nabla \vec{v})$, in particular when
\BEQ {\cal H}_t=\left( \sum_k v_k(t,x)\partial_{x_k}\right)^{\dagger}=-\sum_k
v_k(t,x)\partial_{x_k}-\nabla\cdot \vec{v}(t,x) \EEQ
is in {\em divergence form}, i.e. is the {\em adjoint} of a transport operator. Thus ${\cal H}_t$ is the generator of a
strongly continuous semi-group of contractions of $L^1(\R)$, see e.g. \cite{Paz}, chapter 1. The latter observation extends to the case when 
${\cal H}_t=\left( \sum_k v_k(t,x)\partial_{x_k}\right)^{\dagger} - \tau(t,x)
$  with $\tau(t,\cdot)\le 0$, in the sense that $\int dy\, |f_t(y)|\le \int dy\, |f_T(y)|$ for $0\le t\le T$.

\bigskip


\section{Appendix. Stieltjes transforms} \label{sec:Stieltjes}


We collect in this section some definitions and elementary properties concerning
Stieltjes transforms. We make use of the Fourier transform normalized as follows,
\BEQ {\cal F}(f)(s)=\int_{-\infty}^{+\infty} f(x) e^{-\II xs}\, dx \label{Fourier} \EEQ
with inverse ${\cal F}^{-1}(g)(x)=\frac{1}{2\pi} \int g(s) e^{\II xs}\, ds$.

\medskip

\noindent Let, for $z=a+\II b\in\C\setminus\R$ 
\BEQ {\mathfrak{f}}_z(x)=\frac{1}{x-z}, \qquad x\in\R.\EEQ
For fixed $b\not=0$, ${\mathfrak{f}}_z(x)$ may be seen as  a convolution kernel
$K_b$, 
\BEQ K_b(x-a)=\frac{1}{(x-a)-\II b}. \EEQ

Note that 
\BEQ \Im({\mathfrak{f}}_z)(x)=\frac{b}{|x-z|^2}=\frac{b}{(x-a)^2+b^2},  \qquad \Re({\mathfrak{f}}_z)(x)=\frac{x-a}{(x-a)^2+b^2}.  \label{eq:8.4} \EEQ
In particular, 
\BEQ \Im({\mathfrak{f}}_z)(x)\ge 0 \qquad (b>0). \label{eq:Im-fz>0} \EEQ

Many estimates are based on the simple remark that
\BEQ \int \frac{b}{(x-a)^2+b^2} \, dx=\pi \qquad (b>0) \EEQ
is a constant. The Plemelj
formula,
\BEQ \frac{1}{x-\II 0}=p.v.\left(\frac{1}{x}\right)+\II \pi \del_0 \EEQ
implies the  following boundary value equations for $K_b$,
\BEQ \lim_{b\to 0^+} \int dy\, K_b(x-y) \phi(y) \ -\ \lim_{b\to 0^-} \int dy\, K_b(x-y) \phi(y)=
2\II \pi \phi(x)  \label{eq:bv-} \EEQ
\BEQ \lim_{b\to 0^+} \int dy\, K_b(x-y) \phi(y) \ +\ \lim_{b\to 0^-} \int dy\, K_b(x-y) \phi(y)=
2\  p.v.\ \int \frac{dy}{x-y} \phi(y)  \label{eq:bv+} \EEQ

\medskip

\noindent Then:
\BEQ {\cal F} {\mathfrak{f}}_z(s)=2\II\pi  e^{-b|s|-\II as} {\bf 1}_{s<0} \qquad (b>0), \qquad
-2\II\pi e^{b|s|-\II as} {\bf 1}_{s>0} \qquad (b<0)  \label{eq:8.10} \EEQ
hence (for $b>0$)
\BEQ {\cal F}(\Im({\mathfrak{f}}_z))(s)=\pi e^{-b|s|-\II as}, \qquad {\cal F}(\Re({\mathfrak{f}}_z))=-\II\pi \ \sgn(s) e^{-b|s|-\II as}. \label{eq:8.11} \EEQ

\bigskip

\noindent{\bf Properties of the Stieltjes transform of $\rho_t$.}

\bigskip
\noindent 
Let $M_t(z):=\langle X_t,{\mathfrak{f}}_z\rangle$ ($b:=\Im z>0$). Then:

\BEQ \Im(M_t(z))=\langle X_t,\frac{b}{(x-a)^2+b^2}\rangle; \label{eq:Im-rho-infty} \EEQ
\BEQ |M_t(z)|=|\langle X_t,\frac{1}{(x-a)-\II b} \rangle|\le 1/b; \label{eq:1st-bound-Mt} \EEQ
\BEQ |M'_t(z)|=|\langle X_t,\frac{1}{((x-a)-\II b)^2}\rangle| \le \frac{1}{b}
\langle X_t,\frac{b}{(x-a)^2+b^2}\rangle =\frac{1}{b} \Im(M_t(z)). \label{eq:M'M} \EEQ

When $|a|\gg R$,  we get much better estimates, e.g.
\BEQ |M_t(z)|\le 2/|a|, \qquad |a|\ge 2R. \label{eq:2Nd-bound-Mt} \EEQ

On the other hand, if $b\to 0$ and $a\in \supp(X_t)$, then  $M_t(z)$
 may diverge in general. In particular,
\BEQ |\Re (M_t(z))|\le  C||\rho_t||_{\infty} \ln(R/b) \qquad (b\le \frac{R}{2}, |a|\le 2R),
\label{eq:log-bound} \EEQ
However, if $\rho_t$ is bounded then $\Im( M_t(z)) \in [0,\pi||\rho_t||_{\infty}]$; and
$\Re( M_t(z))=O(1)$ if the space derivative of the density, $\rho'_t$,  is bounded.

\bigskip

\noindent{\bf Some distributions.}

\medskip\noindent Let $\phi\in C_c^{\infty}$ be a  smooth function supported on $[-r,r]$, and $b>0$. Let
\BEQ \langle {\mathfrak{f}}_{\II b},\phi\rangle:=\int dy\, \frac{\phi(y)}{y-\II b} \EEQ
Then  
\BEQ \langle {\mathfrak{f}}_{\II b},\phi\rangle=\phi(0) \int_{-r}^r \frac{dy}{y-\II b} + \II
\int_{-r}^r dy\, (\phi(y)-\phi(0)) \frac{b}{y^2+|b_T|^2} + \int_{-r}^r dy\, 
\frac{y(\phi(y)-\phi(0))}{y^2+b^2} \EEQ
is $O(||\phi||_{\infty}+r||\phi'||_{\infty})$ since: $|\int_{-r}^r \frac{dy}{y-\II b}|\le 
\int dy\, \frac{b}{y^2+b^2}=O(1)$, and $|\frac{y(\phi(y)-\phi(0))}{y^2+b^2}|\le ||\phi'||_{\infty}$. Hence (as seen by integration by parts), $y\mapsto (y-\II b)^{-n}$ ($n\ge 1$) is
a distribution of order $n$, namely,
\BEQ \left|\int dy  \frac{\phi(y)}{(y-\II b)^n} \right|=O(||\phi^{(n-1)}||_{\infty}+
r||\phi^{(n)}||_{\infty}). \label{eq:phi-ib-n} \EEQ

\bigskip
\noindent
{\bf The $\left(\frac{1}{x}\right)$-kernel and its family.}

\medskip\noindent It is known that 
\BEQ {\cal F}^{-1}(s\mapsto \sgn(s){\cal F}f(s))(x)=\II Hf(x):=\frac{\II}{\pi} p.v. \int_{-\infty}^{+\infty} \frac{1}{x-y} f(y)\, dy \label{eq:C1} \EEQ
defined for a compactly supported $f\in C^1$ either as $\frac{\II}{\pi}
\lim_{\eps\to 0^+} \int_{|x-y|>\eps} \frac{1}{x-y}f(y)\, dy$ or as
$\frac{\II}{\pi} \int \frac{f(y)-f(x)}{x-y} \, dy$,  
from which by differentiating
\BEQ {\cal F}^{-1}(s\mapsto |s|{\cal F}f(s))(x)=-\frac{1}{\pi} p.v. \int_{-\infty}^{+\infty} \frac{1}{(x-y)^2} f(y)\, dy \label{eq:C2} \EEQ

For a function $f$  supported on $[-R,R]$, we have the following bounds:
\BEA \left|p.v. \int \frac{1}{x-y}f(y)\, dy\right| &=& {\bf 1}_{|x|\le 2R} 
\left| \int_{x-3R}^{x+3R} \frac{f(y)-f(x)}{x-y}\, dy \right| + {\bf 1}_{|x|> 2R}
\left| \int_{-R}^R \frac{f(y)}{x-y}\, dy\right| \nonumber\\
&=& {\bf 1}_{|x|\le 2R}\  O(R||f'||_{\infty}) +  {\bf 1}_{|x|>2R}\  O(||f||_{\infty} R/|x|), 
\label{eq:bound-vp} \EEA
and similarly
\BEA \left|p.v. \int \frac{1}{(x-y)^2}f(y)\, dy\right| &=& {\bf 1}_{|x|\le 2R} 
\left| \int_{x-3R}^{x+3R} \frac{f'(y)-f'(x)}{x-y}\, dy \right| + {\bf 1}_{|x|> 2R}
\left| \int_{-R}^R \frac{f(y)}{(x-y)^2}\, dy\right| \nonumber\\
&=& {\bf 1}_{|x|\le 2R}\  O(R||f''||_{\infty}) +  {\bf 1}_{|x|>2R}\  O(||f||_{\infty}R/x^2).  \label{eq:bound-vp2}
 \EEA





\bigskip

\noindent { \bf {\em Acknowledgements.}} The author had the opportunity and the pleasure to discuss
 this
project at several occasions and places with N. Simm (Univ. of Warwick), to whom
he therefore wishes to express his gratitude. Referees helped greatly simplify section 5,
pointed out the difficulty with the "martingale term" of section 4, and contributed to
the readability of the manuscript.



\bigskip

\begin{thebibliography}{99}

\bibitem{And} G. W. Anderson, A. Guionnet, O. Zeitouni. {\em An introduction to random
matrices}, Cambridge studies in advanced mathematics {\bf 118}, Cambridge University
Press (2010).

\bibitem{Ben} M. Bender. {\em Global fluctuations in general $\beta$ Dyson
Brownian motion}, Stoch. Proc. Appl. {\bf 118} (6), 1022-1042 (2008).

\bibitem{BenCag} D. Benedetto, E. Cagliotti, J. A. Carrillo, M. Pulvirenti.
{\em A non-Maxwell steady distribution for one-dimensional granular media},
J. Stat. Phys. {\bf 91}, 979--990 (1998).

\bibitem{Car} J. Carrillo, R. McCann, C. Villani. {\em Kinetic equilibration rates for
granular media and related equations: entropy dissipation and mass transportation estimates}, Rev. Mat. Iberoamericanan {\bf 19}, 971--1018 (2003). 

\bibitem{DemZei} A. Dembo, O. Zeitouni. {\em Large deviations techniques
and applications}, Applications of Mathematics {\bf 38}, Springer (2010).


\bibitem{Dui}  M. Duits. {\em On global fluctuations for non-colliding processes},
arxiv.org/abs/1510.08248.

\bibitem{Dys} F. J. Dyson. {\em A Brownian-motion model of the eigenvalues of a random
matrix}, J. Math. Phys. {\bf 3}, 1191--1198 (1962).

\bibitem{EthKur} S. N. Ethier, T. G. Kurtz. {\em Markov proceses: characterization and
convergence}, Wiley, New York (1986). 

\bibitem{For} P. J. Forrester, {\em Log-gases and random matrices}, Princeton University
Press (2010).

\bibitem{ForNag} T. Nagao, P. J. Forrester. {\em Multilevel dynamical correlation functions for
Dyson's Brownian model of random matrices}, Phys. Lett. {\bf A 247}, 42--46 (1998).

\bibitem{FreWen}  M. Freidlin,  A. Wentzell. {\em Random perturbations of
dynamical systems},  Grundlehren der mathematischen Wissenschaften {\bf 260},
Springer (2012).

\bibitem{FyoKhoSim} Y. V. Fyodorov, B. A. Khoruzhenko, N. J. Simm. {\em Fractional
Brownian motion with Hurst index $H=0$ and the Gaussian Unitary Ensemble}, 
Ann. Probab. {\bf 44 (4)}, 2980-3031 (2016).

\bibitem{GolRut} R. J. Goldston, P. H.
Rutherford. {\em An introduction to
plasma physics}, IOP Publishing (1995).


\bibitem{Isr} S. Israelsson. {\em Asymptotic fluctuations of a particle system with singular
interaction}, Stoch. Proc. Appl. {\bf 93}, 25--56 (2001).

\bibitem{Joh} K. Johansson. {\em On fluctuations of eigenvalues of random Hermitian
matrices}, Duke Math J {\bf 91} (1), 15166204 (1998).

\bibitem{KS} I. Karatzas, S. E. Shreve. {\em Brownian motion and stochastic
calculus}, Springer (1991).

\bibitem{Kat} T. Kato. {\em Integration of the equation of evolution in a Banach
space}, J. Math. Soc. Japan {\bf 5}, 208--234 (1953).

\bibitem{KipLan}  C. Kipnis, C. Landim. {\em Scaling limits of interacting particle
systems}, Grundlehren der mathematischen Wissenschaften {\bf 320}, Springer (1999).

\bibitem{LiLiXie} S. Li, X.-D; Li, Y.-X. Xie. {\em Generalized Dyson Brownian
motion, McKean-Vlasov equation and eigenvalues of random matrices}, arxiv.org/abs/1303.1240. {\em On the Law of Large Numbers for the empirical measure process of Generalized Dyson Brownian motion},  arxiv.org/abs/1407.7234.

\bibitem{LodSim} A. Lodhia, N. Simm. {\em Mesoscopic linear statistics of Wigner
matrices}, arxiv.org/abs/1503.03533. 

\bibitem{MacMac} A. F. Macedo, A. M. S. Mac\^edo. {\em Brownian motion ensembles of random matrix
theory: A classification scheme and an integral transform method}, Nucl. Phys. {\bf B 752},
439-475 (2006).

\bibitem{Meh} M. L. Mehta. {\em Random matrices}, Academic Press (1991).

\bibitem{Mit} I. Mitoma. {\em Tightness of probabilities on $C([0,1];{\cal S}')$ and
$D([0,1];{\cal S}')$}, Ann. Prob. {\bf 11 (4)}, 989--999 (1983).

\bibitem{Ott} F. Otto. {\em The geometry of dissipative evolution equations: the porous
medium equation}, Comm. Part. Diff. Eq. {\bf 26}, 101--174 (2001).

\bibitem{OttVil} F. Otto, C. Villani. {\em Generalization of an inequality by Talagrand and
links with the logarithmic Sobolev inequality}, J. Funct. Anal. {\bf 173}, 361--400 (2000).

\bibitem{Pas} L. Pastur, M. Shcherbina. {\em Eigenvalue distribution of large
random matrices}, Mathematical surveys and monographs {\bf 171}, American Mathematical
Society (2011).

\bibitem{Paz} A. Pazy. {\em Semigroups of linear operators and applications to partial
differential equations}, Applied Mathematical Sciences {\bf 44}, Springer-Verlag (1983).

\bibitem{RY}  D. Revuz, M. Yor. {\em Continuous martingales
and Brownian motion}, Grundlehren der  mathematischen Wissenschaften {\bf 293}, Springer (1991).

\bibitem{RogShi} L. Rogers, Z. Shi. {\em Interacting Brownian particles and the 
Wigner law}, Prob. Th. Rel. Fields {\bf 95} (4), 555--570 (1993).

\bibitem{Tan} H. Tanabe. {\em Equations of evolution}, Monographs and studies in
mathematics {\bf 6}, Pitman (1979).

\bibitem{Tao} T. Tao's blog on Brownian motion, https://terrytao.wordpress.com/2010/01/18/254a-notes-3b-brownian-motion-and-dyson-brownian-motion/.

\bibitem{Tre} F. Treves. {\em Topological vector spaces, distributions and
kernels}, Academic Press, New York (1967).

\bibitem{Unt-law} J. Unterberger. {\em Global fluctuations for 1D log-gas
dynamics: the stationary regime} (in preparation).

\bibitem{ValVir} B. Valko, B. Virag. {\em Continuum limits of random matrices and the Brownian carousel}, Invent. math. {\bf 177}, 463--508 (2009).

\bibitem{Vil} C. Villani. {\em Optimal transport, Old and New}, Springer-Verlag, Berlin
(2009).

\end{thebibliography}
\end{document}